\documentclass[12pt]{amsart}

\usepackage{amssymb, amsmath, fancyhdr, amsthm, bm,graphicx}

\usepackage[mathscr]{eucal}

\usepackage{amscd}

\usepackage{verbatim}

\usepackage{graphics}

\input epsf.tex

\renewcommand{\phi}{\varphi}

\newcommand{\Rm}{\mathbb{R}}

\newcommand{\B}{\mathcal{B}}

\newcommand{\fL}{\mathcal B}

\newcommand{\ve}{\varepsilon}

\newcommand{\sB}{\mathcal B}
\newcommand{\sS}{\mathcal B}

\def\hmu{{\hat \mu}}

\def\home{{\hat \omega}}

\def\cF{{\mathcal F}}

\def\cB{\mathcal {B}}
\def\sB{{\mathcal {B}}^{2\odot}}

\def\cN{{\mathcal N}}

\theoremstyle{plain}
\newtheorem{theorem}{Theorem}
\newtheorem{lemma}[theorem]{Lemma}
\newtheorem{prop}[theorem]{Proposition}

\theoremstyle{definition}
\newtheorem{defn}[theorem]{Definition}
\newtheorem{ex}[theorem]{Example}

\numberwithin{theorem}{section}

\title[Ergodic Theory: Nonsingular Transformations]{Ergodic Theory: Nonsingular Transformations\\ $\text{\rm{(This survey is an update of the 2008 version)}}$}

\author{Alexandre~I.~Danilenko}
\address{B. I. Verkin Institute for Low Temperature Physics \& Engineering of
Ukrai\-nian National Academy of Sciences, 47  Nauky Ave.,
Kharkiv, 61164, UKRAINE
}
\email{danilenko@gmail.com}
\author{Cesar E. Silva}
\address{Department of Mathematics, Williams College, Williamstown, MA 01267}
\email{csilva@williams.edu}

\begin{document}

\maketitle

\tableofcontents
\addcontentsline{toc}{section}{Glossary}

\centerline{{\large{G}}{\small{LOSSARY}}}

\bigskip

\noindent{\bf Nonsingular dynamical system:} Let $(X,\fL, \mu)$ be a
standard Borel
space equipped with  a $\sigma$-finite  measure. A Borel map $T:X\to X$ is a {\it nonsingular
transformation}\index{nonsingular transformation} of $X$ if
   for any $N\in \fL$, $\mu(T^{-1}N)=0 \text {  if
  and only if } \mu(N)=0.$
In this case the measure  $\mu$ is called {\it
quasi-invariant\,}\index{quasi-invariant} for $T$; and the quadruple $(X,\sS,\mu,T)$ is called a {\it nonsingular dynamical system}. If $\mu(A)=\mu(T^{-1}A)$ for all $A\in\fL$ then $\mu$ is said to be
{\it invariant } under $T$ or, equivalently, $T$ is
{\it measure-preserving}.
  \smallskip

   \noindent{\bf Conservativity:}
$T$ is {\it
conservative}\index{conservative}  if for all sets $A$ of positive
measure there exists an
integer $n>0$ such that  $\mu(A\cap T^{-n}A)>0$.
  \smallskip

\noindent{\bf Ergodicity:} $T$ is {\it ergodic} if every measurable  subset $A$ of
$X$ that is invariant under $T$ (i.e., $T^{-1}A= A$)  is either $\mu$-null or $\mu$-conull.
Equivalently, every Borel function $f:X\to \mathbb R$ such that $f\circ
T=f$ is constant a.e.
  \smallskip

\noindent{\bf Types II, II$_1$, II$_\infty$ and III:} Suppose that $\mu$ is non-atomic and $T$ is invertible and ergodic  (and hence conservative). If there exists a
$\sigma$-finite measure $\nu$ on $\sS$  which is equivalent to $\mu$
and invariant under $T$ then $T$ is said {\it to be of type II}. It
is easy to see that $\nu$ is unique  up to scaling. If $\nu$ is finite
then $T$ is {\it of type $II_1$}. If $\nu$ is infinite then $T$ is of
type $II_\infty$.
 If $T$ is not of type $II$ then $T$ is said {\it  to
be of type $III$}.

\section{Definition of the subject}
\label{S:Introduction}
An abstract measurable dynamical system consists of a set $X$ (phase space) with a transformation $T:X\to X$ (evolution law or time)
 and a finite or $\sigma$-finite measure $\mu$ on $X$ that specifies a class of  negligible subsets. Nonsingular ergodic theory studies systems where $T$ respects $\mu$ in a weak sense: the transformation preserves only the class of negligible subsets but it may not preserve $\mu$.
This survey is  about dynamics and invariants of nonsingular systems.
Such systems model `non-equilibrium' situations in which events that are impossible at some time remain impossible at any other time.
Of course, the first question that arises is whether it is possible to find an equivalent invariant measure, i.e., pass to a hidden  equilibrium without changing the negligible subsets? It turns out that there exist systems which do not admit an equivalent invariant finite or even $\sigma$-finite measure. They are of our primary interest here.
 In a way (Baire category) most of systems are like that.

Nonsingular dynamical systems arise naturally in various fields of mathematics: topological and smooth dynamics,  probability theory, random walks, theory of numbers, von Neumann algebras, unitary representations of groups, mathematical physics and so on. They also can appear  in the study of probability preserving systems: some criteria of mild mixing and distality, a problem of Furstenberg on disjointness, etc. We briefly discuss this in \S\,\ref{applic}. Nonsingular ergodic theory studies all of them from a general point of  view:
\begin{itemize}
\item[---] What is the qualitative nature of the dynamics?
\item[---] What are the orbits?
\item[---]  Which properties are typical  within a class of systems?
\item[---] How do we  find computable invariants to compare or distinguish
    various systems?
\end{itemize}
Typically there are two kinds of results: some are extensions to
nonsingular systems of theorems for finite measure-preserving
transformations (for instance,  \S\,2, \S\,4, \S\,12) and the other are about
new properly `nonsingular' phenomena (see \S\,5--\S\,9). 
Philosophically
speaking, the dynamics of nonsingular systems is more diverse comparatively
with their finite measure-preserving counterparts. 
That is why it is
usually easier to construct counterexamples than to develop a general
theory. While infinite measure preserving transformations are not the main subject of this survey, we cover them partially as they are also nonsingular systems and arise often as natural examples or counterexamples in the nonsingular setting.
     Because of shortage of space we concentrate mainly on
invertible transformations, and we have not included as many references as we had wished. General group or semigroup actions are practically not considered here (with some exceptions in \S\,\ref{applic} devoted to applications).  A number of open problems are scattered through the entire text.

We thank J. Aaronson,  J. R. Choksi, V. Ya. Golodets, M. Lema\'nczyk, F. Parreau, E. Roy for useful remarks to the first edition of this survey.
Many new results related to nonsingular dynamical systems have appeared since the
release of the first edition.
 The second edition is enlarged essentially to cover (partially) this
progress. 
In particular, we added new \S7 and \S9 and totally rewrote \S8.
More than 100 new references have been added.

\section{Introduction and Basic Results}

This section includes the basic results involving conservativity and ergodicity  as well as some direct nonsingular counterparts of the basic machinery from classic  ergodic theory: mean and pointwise ergodic theorems, Rokhlin lemma, ergodic decomposition, generators, Glimm-Effros theorem and special representation of nonsingular flows. The historically first example of a transformation of type $III$
(due to Ornstein)  is also given here with  full proof.

\subsection{Nonsingular transformations}

In this paper we will consider mainly {\it invertible} nonsingular
  transformations, i.e., those which are bijections when restricted to
an invariant Borel subset of full measure. Thus when we  refer to a nonsingular dynamical system
$(X,\fL,\mu,  T)$ we shall assume that  $T$ is an invertible nonsingular transformation (unless the contrary is specified explicitly).  Of course, each measure
$\nu$ on $\sS$ which is {\it equivalent} to $\mu$, i.e., $\mu$ and
$\nu$ have the same null sets, is also quasi-invariant under $T$. In
particular, since $\mu$ is $\sigma$-finite, $T$ admits an equivalent
quasi-invariant probability measure.
For each $i\in\mathbb Z$, we denote by $\omega_i^\mu$ or  $\omega_i$
the Radon-Nikodym derivative ${d(\mu\circ T^i)}/{d\mu}\in L^1(X,\mu)$.
The derivatives satisfy the cocycle equation
$\omega_{i+j}(x)=\omega_i(x) \omega_j(T^ix)$ for a.e. $x$ and all
$i,j\in\mathbb Z$.

\subsection {Basic properties of conservativity and ergodicity}
\label{S:Recurrence}
  A measurable set $W$ is said to be
{\it wandering}\index{ wandering}  if for all $i, j \geq 0$ with $i\neq j$,
$T^{-i}W\cap
T^{-j}W=\emptyset$.
 Clearly, if $T$ has a wandering set of positive measure then
it cannot be conservative.  A nonsingular transformation $T$ is {\it
incompressible}\index{incompressible}  if whenever $T^{-1}C\subset C$,
then $\mu(C\setminus
T^{-1}C)=0$. 

\begin{prop} {\rm(see e.g. \cite{K85})}{\label{recurrentequivalences}}
  Let $(X,\sS,\mu,T)$ be a
nonsingular dynamical system.  The following are equivalent:
\begin{enumerate}
\item[(i)] $T$ is conservative.

\item[(ii)] For every measurable set $A$,
$
\mu(A\setminus
\bigcup_{n=1}^{\infty}T^{-n}A)=0.
$

\item[(iii)] $T$ is incompressible.

\item[(iv)]  Every wandering set for $T$ is null.

\item[(v)] $\sum_{i=0}^{+\infty}\omega_i(x)=\infty$ at a.e. $x$ (provided that $\mu(X)<\infty$).

\end{enumerate}
\end{prop}

   Since any finite measure-preserving transformation is
incompressible, we deduce  that  it is conservative. This is the
statement of the  classical Poincar\'e recurrence lemma.
If $T$ is a conservative nonsingular transformation of $(X,\sS,\mu)$
and $A\in\sS$ a subset of positive measure, we can define an {\it
induced transformation} $T_A$ of the space $(A,\sS\cap
A,\mu\restriction A)$ by setting $T_Ax:=T^nx$ if $n=n(x)$ is the
smallest natural number such that $T^nx\in A$. $T_A$ is also
conservative. As shown in \cite[5.2]{ST1}, if $\mu(X)=1$ and $T$ is conservative and ergodic, 
 $\int_A \sum_{i=0}^{n(x)-1} \omega_i(x)\ d\mu(x)=1$, which is a nonsingular version of the 
 well-known Ka\c{c}s formula.

\begin{theorem}[Hopf Decomposition, see e.g. \cite{Aa}] \label{T:hopf}
Let $T$ be a nonsingular transformation.
Then there exist disjoint invariant sets $C,D\in\sS$ such that
$X=C\sqcup D$,  $T$ restricted to $C$ is conservative, and $D=
\bigsqcup_{n=-\infty}^\infty T^nW$, where $W$ is a wandering set.
If $f\in L^1(X,\mu)$, $f>0$,  then
$
C=\{x: \sum_{i=0}^{+\infty}  f(T^ix) \omega_i(x) =\infty \text{ a.e.}\}$
and $D=\{x: \sum_{i=0}^{+\infty}  f(T^ix) \omega_i(x) <\infty \text{ a.e.}\}.
$
\end{theorem}

The set $C$ is called the {\it conservative part} of $T$ and $D$ is
called the   {\it dissipative part} of $T$.
If $D$ is of positive measure we call $T$ {\it dissipative}.
If $D$ is of full measure we call $T$ {\it totally dissipative}.

 If $T$ is ergodic and $\mu$ is non-atomic  then
$T$ is automatically conservative. The translation by 1 on the group
$\mathbb Z$ furnished with the counting measure is an example of an
ergodic non-conservative (infinite measure-preserving) transformation.

\begin{prop}\label{fmpergodic}  Let $(X,\sS,\mu,T)$ be a nonsingular
  dynamical
system.
The following are equivalent:
\begin{enumerate}
\item[(i)] $T$ is  conservative and ergodic.

\item[(ii)]  For every set $A$ of positive measure, $\mu(X\setminus
\bigcup_{n=1}^{\infty}T^{-n}A)=0$.  (In this case we will say $A$
sweeps out.)

\item[(iii)] For every measurable set $A$ of positive measure and for a.e.
$x\in X$ there exists an integer $n>0$ such that
$
T^{n}x\in A.
$
\item[(iv)] For all sets $A$ and $B$ of positive measure there exists an
integer $n>0$ such that $\mu(T^{-n}A\cap B)>0$.

\item[(v)]  If $A$ is such that $T^{-1}A\subset A$, then $\mu(A)=0$ or $\mu(A^c)=0$.
\end{enumerate}
\end{prop}

A set $W$ of positive measure is said to be {\it weakly wandering} if there is a sequence
$n_i\to\infty$ such that $T^{n_i}W\cap T^{n_j}W=\emptyset$ for all $i\neq j$.
Clearly, a finite measure-preserving transformation cannot have a weakly wandering set.
Hajian and Kakutani \cite{HK} showed that a nonsingular transformation $T$
is of type $II_1$ if and only if $T$ does not have a weakly wandering set.
This survey is mainly about systems of type $III$.
For some time it was not quite obvious  whether such systems exist at
all. The historically first example was constructed by Ornstein in
1960.

\begin{ex} (Ornstein  \cite{O60})\label{OrnEx}
Let $A_n=\{0,1,\dots,n\}$, $\nu_n(0)=0.5$ and $\nu_n(i)=1/(2n)$ for
$0<i\le n$ and all $n\in\mathbb N$. Denote by $(X,\mu)$ the infinite
product probability space $\bigotimes_{n=1}^\infty(A_n,\nu_n)$. Of
course, $\mu$  is non-atomic. A point of $X$ is an infinite sequence
$x=(x_n)_{n=1}^\infty$ with $x_n\in A_n$ for all $n$.  Given $a_1\in
A_1,\dots,a_n\in A_n$, we denote  the cylinder
$\{x=(x_i)_{i=1}^\infty\in X : x_1=a_1,\dots, x_n=a_n\}$  by $[a_1,\dots,a_n]$.
Define a Borel map $T:X\to X$ by setting

\begin{eqnarray}
\label{odometer}
(Tx)_i=
\begin{cases}
0, & \text{if }i<l(x)\\
x_i+1, &\text{if }i=l(x)\\
x_i, &\text{if }i>l(x),
\end{cases}
\end{eqnarray}
where $l(x)$ is the smallest number $l$ such that $x_l\ne l$. It is
easy to verify that $T$ is a nonsingular transformation of $(X,\mu)$
and
$$
\omega^\mu_1(x)=\prod_{n=1}^\infty\frac{\nu_n((Tx)_n)}{\nu_n(x_n)}=
\begin{cases}
(l(x)-1)!/l(x), &\text{if } x_{l(x)}=0\\
(l(x)-1)!, &\text{if } x_{l(x)}\ne 0.
\end{cases}
$$
We prove that $T$ is of type III by contradiction. Suppose that
there exists a $T$-invariant $\sigma$-finite measure $\nu$ equivalent
to $\mu$.
Let $\phi:=d\mu/d\nu$. Then
\begin{eqnarray} \label{formula}
\omega^\mu_i(x)=\phi(x)\phi(T^ix)^{-1}\text{ \ for a.a. $x\in X$ and
all $i\in\mathbb Z$}.
\end{eqnarray}
Fix a real $C>1$ such that the set $E_C:=\phi^{-1}([C^{-1},C])\subset
X$ is of positive measure.  By a standard approximation argument,  for
each sufficiently large $n$, there is a cylinder $[a_1,\dots,a_n]$
such that $\mu(E_C\cap[a_1,\dots,a_n])>0.9\mu([a_1,\dots,a_n])$. Since
$\nu_{n+1}(0)=0.5$, it follows that
$\mu(E_C\cap[a_1,\dots,a_n,0])>0.8\mu([a_1,\dots,a_n,0])$. Moreover,
by the pigeon hole principle there is $0<i\le n+1$ with
$\mu(E_C\cap[a_1,\dots,a_n,i])>0.8\mu([a_1,\dots,a_n,i])$. Find
$N_n>0$ such that $T^{N_n}[a_1,\dots,a_n,0]=[a_1,\dots,a_n,i]$. Since
$\omega^{\mu}_{N_n}$ is constant on $[a_1,\dots,a_n,0]$, there is a
subset $E_0\subset E_C\cap[a_1,\dots,a_n,0]$ of positive measure such
that $T^{N_n}E_0\subset E_C\cap
[a_1,\dots,a_n,i]$. Moreover, $\omega^\mu_{N_n}(x)=\nu_{n+1}(i)/\nu_{n+1}(0)=(n+1)^{-1}$ for
a.a. $x\in [a_1,\dots,a_n,0]$. On the other hand, we deduce
from~(\ref{formula})
that $\omega^\mu_{N_n}(x)\ge C^{-2}$ for all $x\in E_0$, a contradiction.
\end{ex}

\subsection {Mean and pointwise  ergodic theorems. Rokhlin lemma}
\label{S:ErgodicTheorem}
Let $(X,\sS,\mu, T)$  be a nonsingular dynamical system.
Define a unitary operator $U_T$ of $L^2(X,\mu)$  by setting
\begin{align}\label{koopman oper}
U_Tf:=\sqrt{\omega_1}\cdot f\circ T.
\end{align}
We note that $U_T$ preserves the cone of positive functions $L^2_+(X,\mu)$. Conversely,
every positive unitary operator in $L^2(X,\mu)$ that preserves $L^2_+(X,\mu)$  equals $U_T$ for a $\mu$-nonsingular transformation $T$.
We call $U_T$ the {\it Koopman operator} generated by $T$.

\begin{theorem}
[von Neumann mean Ergodic Theorem, see e.g. \cite{Aa}]  $T$ has no $\mu$-absolutely continuous $T$-invariant probability if and only if  $n^{-1}\sum_{i=0}^{n-1}U_T^i\to 0$ in the strong operator topology.
\end{theorem}

\begin{proof} 
Let $P$ denote the orthogonal projector in $L^2(X,\mu)$ onto the subspace of $U_T$-invariant vectors.
By the well-known fact from the theory of Hilbert spaces,
$n^{-1}\sum_{i=0}^{n-1}U_T^i\to P$ in the strong operator topology.
Then $P\ne 0$ if and only if there is $f\in L^2(X,\mu)$ such that $f\ne 0$ and $U_Tf=f$.
Of course, $U_T|f|=|f|$.
We now define a non-trivial finite measure $\lambda\prec\mu$ by setting $\frac{d\lambda}{d\mu}:=|f|^2$.
It is straightforward  to verify that $\lambda$ is invariant under $T$.
\end{proof}

 Denote by ${\mathcal I}$  the sub-$\sigma$-algebra of $T$-invariant sets. Let  ${\mathbb E}_\mu[.|{\mathcal I}]$  stand for the conditional expectation with respect to ${\mathcal I}$.
 Note that if   $T$ is ergodic, then ${\mathbb E}_\mu[f|{\mathcal I}]=\int f\,d\mu$.
Now we state a nonsingular analogue of Birkhoff's pointwise ergodic theorem, due to Hurewicz \cite{Hur} and in
the form stated by Halmos \cite{Hal46}.

\begin{theorem}
[Hurewicz pointwise Ergodic Theorem]\label{HpET}   If  $T$ is conservative,
$\mu(X)=1$,
 $f , g\in\ L^1(X,\mu)$ and  $g>0$,   then
$$
\frac{\sum_{i=0}^{n-1}  f(T^ix) \omega_i(x) } {\sum_{i=0}^{n-1}  g(T^ix) \omega_i(x)} \to  \frac{{\mathbb E}_\mu[f|{\mathcal I}](x)}
{{\mathbb E}_\mu[g|{\mathcal I}](x)} \ \text{ as }n\to\infty\text{ for a.e. }x.
$$
 \end{theorem}

 There is a nonsingular version of the subadditive ergodic theorem of Kingman \cite{ST1}. 
 Let $T$ be a nonsingular transformation and let $(\omega_n)$ be its sequence of Radon-Nikodym derivatives. 
 A sequence of functions $(f_n)$ is said to be  {\it{subadditive}} if $f_{n+m}\leq f_m+f_n\circ T^m \omega_m$ for all $n,m\ge 0$.
Since $(f_n)$ is a subadditive, one can verify that the following limit 
\[
\mathbb S_\mu[(f_n)](x):=\lim_{n\to\infty}\frac{1}{n} {\mathbb E}_\mu [f_n|{\mathcal I}](x)
\] 
exists almost everywhere.

\begin{theorem}[Nonsingular subadditive Ergodic Theorem]\label{sat}
   If  $T$ is conservative,
$\mu(X)=1$, $(f_n)$ is a subadditive  sequence of integrable functions,
 $ g\in\ L^1(X,\mu)$ and  $g>0$,   then
$$
\frac{f_n(x) } {\sum_{i=0}^{n-1}  g(T^ix) \omega_i(x)} \to  \frac{{\mathbb S}_\mu[(f_n)](x)}
{{\mathbb E}_\mu[g|{\mathcal I}](x)} \ \text{ as }n\to\infty\text{ for a.e. }x.
$$
 \end{theorem}

Of course, Theorem~\ref{HpET} follows from Theorem~\ref{sat} if we set $f_n(x):=
\sum_{i=0}^{n-1}  f(T^ix) \omega_i(x)$ for all $n>0 $ and $x\in X$.

A transformation $T$ is {\it aperiodic } if the  $T$-orbit of a.e. point from $X$ is infinite. The following classical statement can be deduced easily from Proposition \ref{recurrentequivalences}.

\begin{lemma}[Rokhlin's lemma \cite{F70}]\label{Rol} Let $T$ be an aperiodic nonsingular  transformation of a standard probability space $(X,\mu)$.
 For each $\ve>0$ and
integer $N>1$ there exists a measurable set $A$ such that the sets $A, TA,\dots,T^{N-1}A$
are disjoint and $\mu(A\cup TA\cup\cdots\cup T^{N-1}A)>1-\ve$.
\end{lemma}

This lemma was refined later (for ergodic transformations) by Lehrer and Weiss
as follows.

\begin{theorem}[$\epsilon$-free Rokhlin lemma \cite{LeW}] \label{epsilon free} Let $T$ be ergodic and $\mu$ non-atomic. Then for a subset $B\subset X$ and any $N$ for which $\bigcup_{k=0}^\infty T^{-kN}(X\setminus B)=X$, there is a set $A$ such that
the sets $A, TA,\dots,T^{N-1}A$
are disjoint and $A\cup TA\cup\cdots\cup T^{N-1}A\supset B$.
\end{theorem}

The condition $\bigcup_{k=0}^\infty T^{-kN}(X\setminus B)=X$ holds of course for each $B\ne X$ if  $T$ is {\it totally ergodic}, i.e., $T^p$ is ergodic for any $p$, or if $N$ is prime.
We now state a  nonsingular version of Alpern's lemma which is a generalization of 
Lemma~\ref{Rol}.

\begin{theorem}[Alpern's lemma \cite{AP}] Let $T$ be an aperiodic nonsingular transformation of a standard probability space $(X,\mu)$.
Let $\pi=(\pi_1,\pi_2,\dots)$ be a probability vector
such that
$\{k\mid \pi_k>0\}$ is a relatively prime set of integers.
Then there is a measurable partition $P=\{P_{k,i}\mid k>0, i=1,\dots,k\}$ of $X$ satisfying
\begin{enumerate}
\item[(a)] $TP_{k,i}=P_{k,i+1}$ for each $k$ and every $i<k$ and 
\item[(b)]  $\sum_{i=1}^k\mu(P_{k,i})=\pi_k$
 for each $k$.
\end{enumerate}
\end{theorem}

\subsection{Ergodic decomposition}

A proof of the following theorem may be found in \cite[2.2.8]{Aa} and \cite[\S 6]{Sc}.

\begin{theorem}[Ergodic Decomposition Theorem]\label{T:ergodicdecomp}
Let $T$ be a conservative  nonsingular transformation on a standard probability space
$(X,\B,\mu)$. There there exists a standard probability space $(Y,\nu,\mathcal A)$
and a family of probability measures $\mu_y$ on $(X,\B)$, for $y\in Y$, such that
\begin{enumerate}
\item[(i)]  For each $A\in\B$ the map $y\mapsto \mu_y(A)$ is Borel and for each $A\in\B$
\[
\mu(A)=\int \mu_y(A) d\nu(y).
\]
\item[(ii)] For $y,y^\prime\in Y$ the measures $\mu_y$ and $\mu_{y^\prime}$ are mutually singular.
\item[(iii)] For each $y\in Y$ the transformation $T$ is nonsingular and conservative, ergodic
on $(X,\B,\mu_y)$.
\item[(iv)] For each $y\in Y$,
$\omega^{\mu_y}_1=\omega^\mu_1 \
 \mu_y\text{-a.e.}
$
\item[(v)] (Uniqueness) If there exists another probability space $(Y^\prime,\nu^\prime,\mathcal A^\prime)$
and a family of probability measures $\mu^\prime_{y'}$ on $(X,\B)$, for $y'\in Y'$, satisfying
\rm{(i)-(iv)}, then there exists a measure-preserving isomorphism $\theta:Y\to Y^\prime$ such that
$\mu_y = \mu^\prime_{\theta y}$ for $\nu$-a.e. $y$.
\end{enumerate}

\end{theorem}

It follows that if $T$ preserves an equivalent $\sigma$-finite measure
then the system $(X,\B,\mu_y,T)$ is of type $II$ for a.a. $y$.
The space $(Y,\nu,\mathcal A)$ is called {\it the space of  $T$-ergodic components}.

\subsection{Generators}\label{generators} 
A $\sigma$-algebra $\mathcal F$ is called a {\it generator} for a nonsingular  transformation $T$ on a standard probability space $(X,\mathcal B,\mu)$, if  $\bigvee_{n=-\infty}^\infty T^n{\mathcal F}=\mathcal B$. 
It was shown in  \cite{R65} and  \cite{wP66}  that $T$ has a {\it countable   generator}, i.e., a  countable partition $ P$ of $X$ so that the the $\sigma$-algebra of $P$-measurable
sets is a generator for $T$.
It was refined by Krengel \cite{Kr70}: if $T$ is of type $II_\infty$ or $III$ then there exists  $P$ consisting of two sets only. 
Moreover, given a sub-$\sigma$-algebra $\mathcal F\subset\mathcal B$ such that $\mathcal F\subset T\mathcal F$ and $\bigcup_{k>0}T^k\mathcal F=\mathcal B$,  the set $\{A\in \mathcal F\mid (A,X\setminus A)\text{ is a generator of }T\}$ is dense in $\mathcal F$. It follows, in particular, that $T$ is isomorphic to the shift on $\{0,1\}^\mathbb Z$ equipped with a quasi-invariant probability measure.

\subsection{The Glimm-Effros Theorem}\label{S:Glimm}
The classical Bogolyubov-Krylov theorem states that each homeomorphism of a compact space admits an ergodic invariant probability measure~\cite{CFS}. The following statement by Glimm \cite{Gli} and Effros \cite{Eff}  is a ``nonsingular'' analogue of that theorem. (We consider here only a particular  case of $\mathbb Z$-actions.)

\begin{theorem}\label{T:Glimm} Let $X$ be a Polish space and $T:X\to X$ an aperiodic  homeomorphism. Then the following are equivalent:
\begin{itemize}
\item[(i)] $T$ has a recurrent point $x$, i.e., $x=\lim_{n\to\infty}T^{n_i}x$ for a sequence $n_1<n_2<\cdots$.
\item[(ii)]
There is an orbit of $T$ which is not locally closed.
\item[(iii)] There is no a Borel set which intersects each orbit of $T$ exactly once.
\item[(iv)]
 There is a continuous probability Borel measure $\mu$ on $X$ such that $(X,\mu,T)$ is an ergodic nonsingular system.
\end{itemize}
\end{theorem}

A natural question arises: under the conditions of the theorem how many  such $\mu$ can exists? It turns out that there is a wealth of such measures. To state a corresponding result we first write an  important definition.

\begin{defn}\label{orbit equivalent}
Two
nonsingular
systems $(X,\B,\mu,T)$ and
$(X,\B',\mu',T')$ are called  {\it orbit equivalent} if there is a
one-to-one bi-measurable map $\phi:X\to X$ with $\mu'\circ\phi\sim\mu$ and
such that $\phi$ maps the $T$-orbit of $x$ onto the $T'$-orbit of $\phi(x)$ for a.a. $x\in X$.
\end{defn}

The following theorem was proved in \cite{KW0}, \cite{Sc4} and \cite{Kr7}.

\begin{theorem}
Let $(X,T)$ be as in Theorem~\ref{T:Glimm}. Then for each ergodic dynamical system $(Y,\mathcal C,\nu, S)$ of type II$_\infty$ or III, there exist uncountably many mutually disjoint Borel measures $\mu$ on $X$ such that
$(X,T,\B,\mu)$ is orbit equivalent to $(Y,\mathcal C,\nu,S)$.
\end{theorem}

On the other hand,  $T$ may not have any  finite invariant measure.
The first such example appeared in \cite{ehw}.
We present a simpler one.

\begin{ex}
 Let $T$ be an irrational rotation on the circle $\mathbb T$ and let $K$ be a nowhere dense closed subset of $\Bbb T$ of positive Lebesgue measure.
 Let $X$ be the complement of the $T$-orbit $\bigcup_{n\in\Bbb Z}T^nK$ of $K$.
Then  $X$  is a $T$-invariant $G_\delta$-subset of zero Lebesgue  measure.
Hence $X$ is Polish in the induced topology and $T\restriction X$ is an aperiodic homeomorphism of $X$.
Since $T$ is minimal,  $X$ is dense in $\mathbb T$ and the $(T\restriction X)$-orbit of each point of $X$ is dense in $X$.
Hence every point is recurrent. 
By Theorem~2.10, there exists a continuous ergodic nonsingular probability Borel measure $\lambda$ on $X$.
If it is invariant under $T\restriction X$ then $\lambda$ can be considered also as a finite $T$-invariant measure on $\mathbb T$.
Since $T$ is uniquely ergodic,  $\lambda$ is  the Lebesgue measure.
However $X$ is of zero Lebesgue measure, a contradiction.
 \end{ex}

Let $T$ be an aperiodic Borel transformation of a standard Borel space $X$. Denote by $\mathcal M(T)$ the set of all ergodic $T$-nonsingular continuous measures on $X$.
Given $\mu\in\mathcal M(T)$, let $N(\mu)$ denote  the family of all Borel $\mu$-null subsets. Shelah and Weiss showed \cite{SWe} that $\bigcap_{\mu\in\mathcal M(T)}N(\mu)$ coincides with the collection of all Borel $T$-wandering sets.

\subsection{Minimal Radon uniquely ergodic models for infinite measure preserving transformations}
We first note that there is only one, up to a homeomorphism, locally compact non-compact Cantor (i.e., zero-dimensional, perfect, metrizable) set. 
Denote it by $C$.
We recall that a Borel measure on $C$ is called {\it Radon} if it is finite on every compact subset of $C$.
The following is an infinite version of the well known  Jewett-Krieger theorem.

\begin{theorem}[Yuasa, \cite{Yuasa1}]
Let $T$ be an ergodic  measure preserving transformation of the standard infinite $\sigma$-finite
 measure space $(X,\mu)$.
Then there exists a minimal homeomorphism $R$ of $C$ that admits a unique, up to scaling, $R$-invariant Radon measure $\nu$ such that $(X,\mu,T)$ is isomorphic to $(C,\nu,R)$.
\end{theorem}

 Yuasa proved a relative version of this theorem in  \cite{Yuasa2}.
Similar  results on {\it strictly ergodic models} for ergodic systems of type $III$ are unknown yet.

\subsection{Special representations of ergodic flows}
Nonsingular flows (=$\mathbb R$-actions) appear naturally in the study of orbit equivalence for systems of type $III$ (see Section~\ref{orbits}).
Here we record some basic notions related to nonsingular flows.
Let $(X,\B,\mu)$ be a standard Borel space with a $\sigma$-finite measure $\mu$ on $\B$.
A nonsingular {\it flow} on $(X,\mu)$ is a Borel map $S:X\times\Rm\ni(x,t)\mapsto S_tx\in X$ such that $S_tS_s=S_{t+s}$
for all $s,t\in\mathbb R$ and each $S_t$ is a nonsingular transformation of $(X,\mu)$.  Conservativity  and ergodicity for flows are defined in a similar way as for transformations.

A very useful example of a flow is a flow built under a function.  Let $(X,\B,\mu,T)$ be a nonsingular dynamical system
and $f$ a positive Borel function on $X$ such
that $\sum_{i=0}^\infty f(T^ix)=\sum_{i=0}^\infty f(T^{-i}x)=\infty$  for all $x\in X$.  Set
$X^{f}:=\{(x,s): x\in X, 0\leq s <f(x)\}.$
Define $\mu^{f}$ to be the restriction of the product measure $\mu\times\text{Leb}$ on $X\times\Rm$ to $X^{f}$ and
define, for 
$t \geq 0$,

\[S^{f}_{t(x,s)} : = (T^n x,s+t-\sum_{i=0}^{n-1}f(T^ix),\]

where $n$ is the unique integer that satisfies
\[
 \sum_{i=0}^{n-1}f(T^ix)<s+t\leq \sum_{i=0}^{n}f(T^ix).
\]
A similar definition applies when $t<0$. In particular, when $0<s+t<\phi(x) $, $S^f_t(x,s)=(x,s+t)$,
so that the flow moves the point $(x,s)$ up $t$ units, and when it reaches $(x,\phi(x))$ it is sent
to $(Tx,0)$. It can be shown  that $S^f=(S^f_t)_{t\in\mathbb R}$ is a free $\mu^{f}$-nonsingular flow   and that it preserves $\mu^{f}$ if and only if  $T$ preserves $\mu$ \cite{Nad}.
It is called the {\it flow built under the function $\phi$ with the base transformation $T$}. Of course, $S^f$ is conservative or ergodic if and only if so is $T$.

Two flows $S=(S_t)_{t\in\mathbb R}$ on $(X,\B,\mu)$ and $V=(V_t)_{t\in\mathbb R}$ on $(Y,\mathcal C,\nu)$
are said to be {\it isomorphic} if there exist invariant co-null  sets $X^\prime\subset X$ and $Y^\prime\subset Y$ and  an invertible nonsingular map  $\rho:X^\prime\to Y^\prime$  that interwines the actions of
the flows: $\rho\circ S_t=V_t\circ\rho$ on $X^\prime$ for all $t$.  The following nonsingular version of Ambrose--Kakutani representation theorem was proved by Krengel \cite{Kre69} and  Kubo \cite{Ku}.

\begin{theorem}\label{ambrose} Let $S$ be a free nonsingular flow. Then it is isomorphic to a flow built under a function.
\end{theorem}

Rudolph showed that in the Ambrose-Kakutani  theorem one can choose the function $\phi$ to take two values.  Krengel \cite{Kre76} showed that this can also be assumed in the nonsingular case.

\section{Panorama of  Examples}
\label{S:BasicExample}
 This section is devoted entirely to examples of nonsingular systems. We describe here the most popular (and simple) constructions of nonsingular systems: product odometers, nonsingular Markov odometers, tower transformations, rank-one and finite rank systems,  nonsingular Bernoulli  and Markov shifts.

\subsection{Nonsingular product odometers}\label{SS:odometerdef}

Given a sequence $m_n$ of natural numbers, we let
$A_n:=\{0,1,\dots,m_n-1\}$. Let $\nu_n$ be a probability on $A_n$ and $\nu_n(a)>0$ for all $a\in A_n$. Consider now the infinite product probability space
$(X,\mu):=\bigotimes_{n=1}^\infty(A_n,\nu_n)$.
Assume that $\prod_{n=1}^\infty\max\{\nu_n(a)\mid a\in A_n\}=0$. Then $\mu$
is non-atomic. Given $a_1\in A_1,\dots,a_n\in A_n$, we denote by $[a_1,\dots,a_n]$ the cylinder ${x=(x_i)_{i>0}\mid x_1=a_1,\dots,x_n=a_n}$. If $x\ne (0,0,\dots)$, we let $l(x)$ be the smallest number $l$ such that the $l$-th coordinate of $x$ is not $m_l-1$.
 We define a Borel map $T:X\to X$ by (\ref{odometer}) if $x\ne (m_1,m_2,\dots)$ and  put $Tx:=(0,0,\dots)$ if $x=(m_1,m_2,\dots)$. Of course, $T$ is isomorphic to a rotation on a compact monothetic totally disconnected Abelian group. It is easy to check that $T$ is $\mu$-nonsingular and
$$
\omega_1^\mu(x)=\prod_{n=1}^\infty\frac{\nu_n((Tx)_n)}{\nu_n(x_n)}=
\frac{\nu_{l(x)}(x_{l(x)}+1)}
{\nu_{l(x)}(x_{l(x)})}\prod_{n=1}^{l(x)-1}\frac{\nu_n(0)}{\nu_n(m_n-1)}
$$ for a.a. $x=(x_n)_{n>0}\in X$.
It is also easy to verify that $T$ is ergodic. It is called the {\it
nonsingular product odometer} associated to $(m_n,\nu_n)_{n=1}^ \infty$. We note
that Ornstein's transformation (Example \ref{OrnEx}) is a nonsingular product 
odometer.

\subsection{Markov odometers}\label{MarkOd}

We define Markov odometers as in \cite{DoH1}. An ordered Bratteli
diagram $B$ \cite{HPS} consists of
\begin{enumerate}
\item[(i)]
a vertex set $V$ which is a disjoint union of finite sets $V^{(n)}$, $n\ge
0$, $V_0$ is a singleton;
\item[(ii)]
an edge set $E$ which is a disjoint union of finite sets $E^{(n)}$, $n>0$;
\item[(iii)]
source mappings $s_n:E^{(n)}\to V^{(n-1)}$ and range mappings
$r_n:E^{(n)}\to V^{(n)}$ such that $s_n^{-1}(v)\ne\emptyset$ for all $v\in
V^{(n-1)}$ and $r_n^{-1}(v)\ne\emptyset$ for all $v\in V^{(n)}$, $n>0$;
\item[(iv)]
a partial order on $E$ so that $e,e'\in E$ are comparable if and only if
$e,e'\in E^{(n)}$ for some $n$ and $r_n(e)=r_n(e')$.
\end{enumerate}
A {\it Bratteli compactum} $X_B$ of the diagram $B$ is the space of
infinite paths
$$
\{x=(x_n)_{n>0}\mid x_n\in E^{(n)} \text{ and }r(x_n)=s(x_{n+1})\}
$$
on $B$. $X_B$ is equipped with the natural topology induced by the product
topology on $\prod_{n>0}E^{(n)}$. We will assume always that the diagram is
{\it essentially simple}, i.e., there is only one infinite path
$x_{\max}=(x_n)_{n>0}$ with $x_n$ maximal for all $n$ and only one
$x_{\min}=(x_n)_{n>0}$ with $x_n$ minimal for all $n$. The {\it
Bratteli-Vershik} map $T_B:X_B\to X_B$ is defined as follows:
$Tx_{\max}=x_{\min}$. If $x=(x_n)_{n>0}\ne x_{\max}$ then let $k$ be the
smallest number such that $x_k$ is not maximal. Let $y_k$ be a successor of
$x_k$. Let $(y_1,\dots,y_k)$ be the unique path such that
$y_1,\dots, y_{k-1}$ are all minimal. Then we let
$T_Bx:=(y_1,\dots,y_k,x_{k+1},x_{k+2},\dots)$. It is easy to see that $T_B$
is a homeomorphism of $X_B$. Suppose that we are given a sequence
$P^{(n)}=(P^{(n)}_{(v,e)\in V^{n-1}\times E^{(n)}})$ of stochastic
matrices, i.e.,
\begin{enumerate}
\item[(i)]
$P^{(n)}_{v,e}>0$ if and only if $v=s_n(e)$ and
\item[(ii)] $\sum_{\{e\in E^{(n)}\mid s_n(e)=v\}}P^{(n)}_{v,e}=1$ for each
$v\in V^{(n-1)}$.
\end{enumerate}
For $e_1\in E^{(1)},\dots, e_n\in E^{(n)}$, let $[e_1,\dots,e_n]$ denote
the cylinder $\{x=(x_j)_{j>0}\mid x_1=e_1,\dots, x_n=e_n\}$. Then we define
a {\it Markov measure} on $X_B$ by setting
$$
\mu_{P}([e_1,\dots,e_n])=P^{1}_{s_1(e_1),e_1}P^{2}_{s_2(e_2),e_2}\cdots
P^{n}_{s_n(e_n),e_n}
$$
for each cylinder $[e_1,\dots,e_n]$. The dynamical system $(X_B,\mu_{P},
T_B)$ is called a {\it Markov odometer}. It is easy to see that every
nonsingular product odometer is a Markov odometer where the corresponding $V^{(n)}$
are all singletons.

\subsection{Tower transformations}\label{tower} This construction is a discrete analogue of flow under a function. Given a nonsingular dynamical system
$(X,\mu,T)$ and a measurable map $f:X\to\mathbb N$, we define a new dynamical system $(X^f,\mu^f,T^f)$ by setting
$$
\begin{aligned}
X^f &:=\{(x,i)\in X\times\mathbb Z_+\mid 0\le i< f(x)\},\\
d\mu^f(x,i) &:=d\mu(x) \ \text{and}\\
T^f(x,i) &:=
\begin{cases}
(x,i+1), & \text{if }i+1<f(x)\\
(Tx,0),  & \text{ otherwise.}
\end{cases}
\end{aligned}
$$
Then $T^f$ is $\mu^f$-nonsingular and $(d\mu^f\circ T^f/d\mu^f)(x,i)=(d\mu\circ T/d\mu)(x)$ for a.a. $(x,i)\in X^f$.
This transformation is called the (Kakutani) {\it tower over $T$ with  height function} $f$. It is easy to check that $T^f$ is conservative if and only if $T$ is conservative; $T^f$ is ergodic if and only if $T$ is ergodic; $T^f$ is of type $III$ if and only if $T$ is of type $III$. Moreover, the induced transformation  $(T^f)_{X\times \{0\}}$ is isomorphic to $T$. Given a subset $A\subset X$ of positive measure, $T$ is the tower over the induced transformation $T_A$ with the first return time to $A$ as the height function.

\subsection{Rank-one transformations. Chac\'on maps. Finite rank}\label{S:rankone}
The definition
 uses the process of  ``cutting and stacking.''
We construct by induction a sequence of columns $C_n$. A {\it column}  $C_n$ consists of a
finite sequence of bounded intervals (left-closed, right-open)  $C_n=\{I_{n,0},\dots,I_{n,h_n-1}\}$ of
{\it height} $h_n$.  A column $C_n$
determines a {\it column map} $T_{C_n}$ that sends each interval $I_{n,i}$ to the interval above it
$I_{n,i+1}$ by the unique orientation-preserving affine map between the intervals. $T_{C_n}$ remains undefined on the top interval $I_{n,h_n-1}$. Set $C_0=\{[0,1)\}$ and let $\{r_n \geq 2\}$ be a sequence of positive integers, let $\{s_{n}\}$ be a sequence of functions $s_n:\{0,\ldots, r_n-1\}\to \mathbb N_0$,  and let $\{w_{n}\}$ be a  sequence of probability vectors on $\{0,\ldots, r_n-1\}$.  If   $C_n$ has been defined, column $C_{n+1}$ is defined as follows. First  ``cut'' (i.e., subdivide) each interval $I_{n,i}$ in $C_n$ into $r_n$ subintervals $I_{n,i}[j], j=0,\ldots,r_n-1$,
whose lengths are in the proportions $w_n(0) : w_n(1) : \cdots : w_n(r_n-1)$.  Next  place, for each
$j=0,\ldots r_n-1$, $s_{n}(j)$
 new subintervals above $I_{n,h_n-1}[j]$, all of the same  length as $I_{n,h_n-1}[j]$.   Denote these intervals, called {\it spacers}, by
$S_{n,0}[j], \ldots S_{n,s_{n}(j)-1}[j]$. This yields, for each $j\in\{0,\ldots,r_n-1\}$, $r_n$ subcolums
each consisting of the subintervals
\[I_{n,0}[j],\ldots I_{n,h_n-1}[j]\text{\  followed by the spacers\ }  S_{n,0}[j], \ldots S_{n,s_{n}(j)-1}[j].
\]
Finally each subcolumn is stacked from left to right so that the top subinterval in subcolumn $j$ is sent to the bottom subinterval in subcolumn $j+1$, for $j=0,\ldots,r_n-2$ (by the unique orientation-preserving affine map between the intervals).
For example, $S_{n,s_{n}(0)-1}[0]$ is sent to $I_{n,0}[1]$. This defines a new column $C_{n+1}$ and new column map $T_{C_{n+1}}$, which remains undefined on its top subinterval. Let $X$ be the union of all intervals in all columns and let $\mu$ be Lebesgue measure restricted to $X$.
We assume that as $n\to\infty$ the maximal length of the intervals in $C_n$
converges to $0$, so we may define a transformation $T$ of $(X,\mu)$ by $Tx:=\lim_{n\to\infty} T_{C_{n}}x$.
One can verify that $T$ is well-defined a.e. and that it is  nonsingular and ergodic. $T$ is said to be the {\it rank-one} transformation associated with $(r_n,w_n,s_{n})_{n=1}^\infty$. If all the probability vectors $w_n$ are uniform the resulting transformation is measure-preserving. The measure is infinite ($\sigma$-finite) if and only if the total mass of the spacers is infinite. In the case $r_n=3$ and $s_n(0)=s_n(2)=0$, $s_n(1)=1$ for all $n\geq 0$, the associated rank-one transformation is called a {\it nonsingular Chac\'on map}.

It is easy to see that every nonsingular product odometer is of rank-one (the corresponding maps $s_n$ are all trivial).
Each rank-one map $T$ is a tower over a nonsingular product odometer (to obtain such an odometer reduce $T$ to a column $C_n$).

A rank $N$ transformation is defined in a similar way.  A nonsingular transformation $T$ is said to be of
{\it rank $N$ or less} if at each stage of its construction there exits $N$ disjoint columns,  the levels of the columns generate the $\sigma$-algebra and the Radon-Nikodym derivative of $T$ is constant on each non-top level of every column. $T$ is said to be of {\it rank $N$} if it is of rank $N$ or less and not of
rank $N-1$ or less. A rank $N$ transformation, $N\geq 2$, need not be ergodic.

\subsection{Nonsingular Bernoulli shifts}\label{S:hamachi}
  A {\it nonsingular Bernoulli} transformation is  a  transformation $T$ such that there exists a  generating $\sigma$-algebra  $\mathcal P$ for $T$ (see \S\,\ref{generators}) such that the $\sigma$-algebras  $T^n{\mathcal P}$, $n\in\mathbb Z$, are mutually independent.
Thus we may think that $T$ is the left shift on the probability space $(X,\mu)=(A^\Bbb Z,\bigotimes_{n\in\Bbb Z}\mu_n)$, where $A$ is a standard Borel space and $(\mu_n)_{n\in\Bbb Z}$ is a sequence
of probability measures on $A$.
We will always assume that 
$\mu$ is nonatomic. 
It follows from Kakutani's criterion for equivalence of infinite  product measures \cite{Kak}  that $\mu$ is nonsingular  if and only if $\mu_n\sim\mu_{n+1}$ for each $n$ and
\begin{align}\label{ffff}
\sum_{n\in\Bbb Z} H^2(\mu_n,\mu_{n+1})<\infty,
\end{align}
where $H^2(\mu,\nu)$ denotes the {\it Hellinger distance} defined by
$H^2(\mu,\nu):=1-\int_X\sqrt{\frac{d\mu}{d\xi}\frac{d\nu}{d\xi}}d\xi$, where $\xi$ is a probability
such that $\mu\prec\xi$ and  $\nu\prec\xi$.
If (\ref{ffff}) is satisfied then for a.e. $x=(x_n)_{n\in\Bbb Z}\in X$,
$$
\omega^\mu_1(x)=\prod_{n\in\Bbb Z}\frac{\mu_{n-1}(x_n)}{\mu_n(x_n)}.
$$

 \subsection{Nonsingular Markov shifts}\label{S:MarkSh} Let $A$ be a finite set and let $M=(M(a,b))_{a,b\in A}$ be a 0-1-valued  $A\times A$-matrix.
 We let $X_M:=\{x=(x_i)_{i\in\Bbb Z}\in A^{\Bbb Z}\mid  M(x_i,x_{i+1})>0\text{ for all }i\in\Bbb Z\}$.
 Let $T$ denote the restriction of the left shift to $X_M$.
 Then $T$ is a  
 shift of finite type.
 Given two integers $i\le j$, and a finite sequence $a=(a_l)_{l=i}^j$ of elements from $A$ such that $M(a_l,a_{l+1})=1$ for $l=i,\dots,j-1$, we define by $[a]_i^j$ the cylinder $\{x\in X_M\mid x_l=a_l\text{ for all }l=i,\dots,j\}$.
 Suppose that there is a sequence of probability measures $(\pi_n)_{n\in\Bbb Z}$ on $A$  and a sequence $(P_n)_{n\in\Bbb Z}$ of row-stochastic $A\times A$-matrices such that
 $\pi_nP_n=\pi_{n+1}$ and $P_n(a,b)>0$ if and only if $M(a,b)>0$ for each $n\in\Bbb Z$.
 Then there is a unique probability measure $\mu$ on $X_M$ such that for every cylinder
 $[a]_i^j$ in $X_M$.
 $$
 \mu([a]_i^j)=\pi_i(a_i)P_i(a_i,a_{i+1})\cdots P_{j-1}(a_{j-1},a_j).
 $$
 It is called Markov measure on $X_M$ generated by $(\pi_n,P_n)_{n\in\Bbb Z}$.
 By analogy with \cite{Kak}, using \cite{KLS} one can find necessary and sufficient conditions for $T$ to be $\mu$-nonsingular.
 Then the system $(X_M,T,\mu)$ is called a {\it nonsingular Markov shift}.
It is a natural generalization of nonsingular Bernoulli shifts.
By a standard computation,
$$
\omega^\mu_1(x)=\lim_{n\to+\infty}\frac{\pi_{-n-1}(x_{-n})}{\pi_{-n}(x_{-n})}\prod_{j=-n}^n\frac{P_{j-1}(x_j,x_{j+1})}{P_j(x_j,x_{j+1})}.
$$

 We also mention here so-called {\it infinite Markov shifts}, i.e., Markov transformations preserving an infinite  $\sigma$-finite measure (see \cite{KP} and  \S4.5 from \cite{Aa}).
 Let $A=\Bbb Z$ and let $P=(P(a,b))_{a,b\in A}$ be a row stochastic  $A\times A$-matrix.
 Suppose that $P$ is irreducible, i.e., for each pair $a,b\in A$, there is $n>0$ with $P^n(a,b)>0$.
 Suppose that there is   a strictly positive function $\pi:A\to\Bbb R_+^*$ such that
 $\sum_{a\in A}\pi(a)=\infty$ and
 $\sum_{a\in A}\pi(a)P(a,b)=\pi(b)$ for all $b\in A$.
 Let $T$ denote the restriction of the left shift to $X_P$.
 Define a measure $\mu$ on $X_P$ by setting
 $$
 \mu([a]_i^j)=\pi(a_i)P(a_i,a_{i+1})\cdots P(a_{j-1},a_j).
 $$
 Then $\mu$ is infinite and $\sigma$-finite and $T$ preserves $\mu$.
 By \cite{KP}, if $\sum_{n=1}^\infty P^k(0,0)=\infty$ then $T$ is ergodic.
 We call the system $(X_P,\mu,T)$ {\it the infinite Markov shift associated with $(P,\pi)$}.

 \subsection{Nonsingular Poisson suspentions}\label{NPT}
 Let $(X,\mathcal B,\mu)$ be a $\sigma$-finite Lebesgue space with an infinite non-atomic measure.
 Let $\mathcal B_0\subset\mathcal B$ be the collection of subsets of finite measure.
 Denote by $X^*$ the space of measures $\omega$ of the form  
 $\omega=\sum_{i\in I}{\delta_{x_i}}$, where $I$ is a countable set.
 Endow $X^*$ with the smallest $\sigma$-algebra $\mathcal B^*$ such that the maps 
 $$
 N_A:X^*\ni\omega\mapsto\omega(A)\in\Bbb Z_+\sqcup\{\infty\}
 $$
  are all measurable, $A\in\mathcal B_0.$
 Let $\mu^*$ be the (only) probability measure on $\mathcal B^*$
 such that 
 
 \begin{itemize}
 \item
 $\mu^*\circ N_A^{-1}$ has the Poisson distribution with parameter $\mu(A)$, $A\in\mathcal B_0$ and
 \item
 if $A,B\in\mathcal B_0$ and $A\cap B=\emptyset$ then the maps
 $N_A$ and $N_B$ are independent.
 \end{itemize} 
 Let $T$ be a nonsingular invertible transformation of $(X,\mathcal B,\mu)$.
Define a transformation $T_*:X^*\to X^*$ by setting: $T_*\omega:=\omega\circ T^{-1}$.
If $T_*$ is $\mu^*$-nonsingular then $(X^*,\mathcal B^*,\mu^*, T_*)$ is called {\it the nonsingular Poisson suspension} of $(X,\mathcal B,\mu, T)$ \cite{DaKoRo1}.

 \begin{theorem}[\cite{DaKoRo1}]\label{nonsing-P} 
 $T_*$ is $\mu^*$-nonsingular if and only if
 $\sqrt{\frac{d\mu\circ T}{d\mu}}-1\in L^2(\mu)$.
 If  $\frac{d\mu\circ T}{d\mu}-1\in L^1(\mu)$ then   $T_*$ is $\mu^*$-nonsingular  and
 $$
 \frac{d\mu^*\circ T_*}{d\mu^*}(\omega)=e^{-\int_X(\frac{d\mu\circ T}{d\mu}-1)d\mu}
 \prod_{\omega(\{x\})=1}\frac{d\mu\circ T}{d\mu}(x)
 $$
 at a.e. $\omega\in X^*$.
 \end{theorem}

 \subsection{Nonsingular Gaussian transformations}\label{gaus}
 Let $\gamma$ stand for the normalized Gaussian measure on $\Bbb R$:
 $d\gamma(t)=\frac1{\sqrt{2\pi}}e^{-\frac{t^2}2}dt$.
 Let $(X,\mu):=(\Bbb R,\gamma)^{\Bbb N}$.
 Denote by $\mathcal O$ the group of orthogonal operators in a separable infinite dimensional real Hilbert space $\mathcal H$.
 Fix an orthonormal base in $\mathcal H$.
 Then we can identify $\mathcal H$ with $\ell^2(\Bbb N)$ and every operator from $\mathcal O$ with
 an infinite $\Bbb N\times\Bbb N$-matrix.
Every $O\in\mathcal O$ determines a Borel transformation $T_O:X\to X$ by the formula
$$
T_Ox:=((T_Ox)_n)_{n\in\Bbb N}, \text{ where } (T_Ox)_n:=\sum_{m=1}^\infty O_{n,m}x_m,\  n\in\Bbb N. 
$$
More precisely,  $T_O$ is defined $\mu$-almost everywhere and it is invertible (mod 0).
Moreover, $T_O$ preserves $\mu$. 
The set $\{T_O\mid O\in\mathcal O\}$
 is  exactly the family of classical (probability preserving) Gaussian transformations, which are
well studied in ergodic theory.
A wider class of  nonsingular Gaussian transformations is defined by Arano, Isono and Marrakchi in \cite{AIM} (see also \cite{DanL22} for a different but equivalent presentation of their concepts).
For each $y\in X$, consider a transformation $S_y:X\ni x\mapsto x+y\in X$.
By the Cameron-Martin theorem,  $S_y$ is $\mu$-nonsingular if and only if $y\in\ell^2(\Bbb N)$ and
$$
\frac{d\mu\circ S_{y}^{-1}}{d\mu}(x)=e^{-\frac12\|y\|_2^2+\sum_{n=1}^\infty x_ny_n}\text{\ at $\mu$-a.e. $x=(x_n)_{n=1}^\infty\in X$ and $y=(y_n)_{n=1}^\infty\in\ell^2$.}
$$
It is easy to see that $S_y$ is totally dissipative.
It is straightforward to verify that $T_OS_yT_O^{-1}=S_{Oy}$.
Denote the affine group  $ \mathcal H\rtimes \mathcal O$ of $\mathcal H$ by Aff$\,\mathcal H$.
Then for each $(h,O)\in \text{Aff}\,\mathcal H$, consider a transformation
 $G_{h,O}:=S_hT_O$ of $(X,\mu)$.
It is called {\it a nonsingular Gaussian transformation}.
Thus, each nonsingular Gaussian transformation is a composition of a classical probability preserving Gaussian transformation and a  nonsingular  totally dissipative ``rotation''. 
It is obvious  that 
$$
\frac{d\mu\circ G_{h,O}^{-1}}{d\mu}(x)=\frac{d\mu\circ S_{h}^{-1}}{d\mu}(x)=e^{-\|h\|_2^2+\sum_{n=1}^\infty x_nh_n}\text{\ at $\mu$-a.e. $x=(x_n)_{n=1}^\infty\in X$.}
$$

\subsection{IDPFT transformations}

\begin{defn}[\cite{DanL18}] 
Let $T_n$ be an ergodic nonsingular invertible transformation of a standard probability space $(X_n,\frak B_n,\mu_n)$ for each $n\in\Bbb N$. 
Denote by $T$ the infinite direct product of $T_n$, $n\in\Bbb N$, acting on the infinite product space $(X,\frak B,\mu):=\bigotimes_{n\in\Bbb N}(X_n,\frak B_n,\mu_n)$.
If $T$ is $\mu$-nonsingular and each $T_n$ is of finite type, i.e that there exists a 
$\mu_n$-equivalent probability measure $\nu_n$ which is invariant under $T_n$ for each $n\in\Bbb N$
then $T$ is called an {\it infinite direct product of finite types (IDPFT).}
\end{defn}

Let $\phi_n:=\frac{d\mu_n}{d\nu_n}$.
It follows from the Kakutani criterion \cite{Kak} that $T$ is $\mu$-nonsingular if and only if
$\prod_{n=1}^\infty\int_{X_n}\sqrt{\phi_n\cdot\phi_n\circ T_n}>0$.
Moreover, $\mu$ is mutually singular with the $T$-invariant probability $ \nu:=\bigotimes_{n\in\Bbb Z}\nu_n$ if and only if $\prod_{n=1}^\infty\int_{X_n}\sqrt{\phi_n}\,d\nu_n=0$.
If $T$ is $\mu$-nonsingular then
$$
\frac{d\mu\circ T}{d\mu}(x)=\prod_{n=1}^\infty\frac{d\mu_n\circ T_n}{d\mu_n}(x_n)\ \text{\ at a.e. }x=(x_n)_{n=1}^\infty\in X
.$$

 \subsection{Natural extensions of nonsingular endomorphisms}
 Let $(X,\mathcal B,\mu)$ be a $\sigma$-finite standard measure space.
 A {\it nonsingular endomorphism} is a measurable map $R:X\to X$ such that $\mu(A)=0$ if and only if $\mu(R^{-1}A)=0$.
 Suppose that $\mu$ is $\sigma$-finite on $R^{-1}\mathcal B$.
 We define the {\it Radon-Nikodym derivative} $\omega_1^\mu$ of $R$ by setting
 $\omega_1^\mu=\frac{d\mu}{d\mu\circ R^{-1}}\circ R$. 
 It was shown in \cite{Si} and \cite{ST2}
 that there exists a $\sigma$-finite standard measure space $(X^*,\mathcal B^*,\mu^*)$, an invertible $\mu^*$-nonsingular transformation $R^*$ and a Borel map $\pi:X^*\to X$
 such that the following hold: $\mu^*\circ\pi^{-1}=\mu$, $\pi R^*=R\pi$, $\omega_1^{\mu^*}$ is $\pi^{-1}(\mathcal B)$-measurable and $\bigvee_{n>0}R^n\pi^{-1}(\mathcal B)=\mathcal B^*$.
 The dynamical system $(X^*,\mathcal B^*,R^*,\mu^*)$ is defined uniquely (up to a natural isomorphism) and called {\it the natural extension of $R$}.
 It coincides with the standard Rokhlin definition of the natural extension  in the case where $R$ preserves $\mu$ and $\mu$ is finite.
 
\begin{theorem}[\cite{Si}, \cite{ST2}] $R^*$ is conservative if and only if  $R$ is $\mu$-recurrent, i.e.,
$$
\sum_{i\ge 0}h\circ R^i\omega_i=+\infty\quad \text{a.e., where $\omega_i=\prod_{j=0}^{i-1}\omega_1^{\mu^*}\circ R^j$}
$$
for each integrable function $h>0$.
Moreover, if $R$ is $\mu$-recurrent then $R^*$ is ergodic if and only if $R$ is ergodic.
\end{theorem}

 Let $R$ be a nonsingular one-sided Bernoulli shift $(X,\mu)=\bigotimes_{n=1}^\infty (A,\mu_n)$.
 Then the natural extension of $R$ is isomorphic to the two-sided nonsingular Bernoulli shift $T$ on
$(X^*,\mu^*)=\bigotimes_{n=-\infty}^\infty (A,\mu_n^*)$, where $\mu_n^*=\mu_n$ if $n>0$
and $\mu_n^*=\mu_1$ if $n\le 0$.
The corresponding projection $\pi:X^*\to X$ is the natural projection, i.e., $\pi(\dots,a_{-1},a_0,a_1,a_2,\dots):=(a_1,a_2,\dots)$.

 \section{Topological groups Aut$(X,\mu)$, Aut$_2(X,\mu)$ and Aut$_1(X,\mu)$}\label{topological group}

 \subsection{} Let  $(X,\mathcal B,\mu)$ be a standard probability space and let  Aut$(X,\mu)$ denote the group of all  nonsingular transformations
    of $X$. Let $\nu$ be a finite or $\sigma$-finite measure equivalent to $\mu$; the subgroup of the $\nu$-preserving transformations is denoted by Aut$_0(X,\nu)$. Then Aut$(X,\mu)$ is a simple group \cite{Eig} and it has no outer automorphisms \cite{Eig2}. Ryzhikov showed \cite{Ry} that every element of this group is a product of three involutions (i.e., transformations of order~2). Moreover, a nonsingular transformation is a product of two involutions if and only if it is conjugate to its inverse by an involution.

Inspired by \cite{Hal}, Ionescu Tulcea \cite{IT65} and Chacon and Friedman \cite{ChF} introduced the {\it weak} and the {\it uniform} topologies  respectively on  Aut$(X,\mu)$. The weak one---we denote it by $d_w$---is induced from the weak  operator topology on the group of unitary operators in $L^2(X,\mu)$ by the embedding $T\mapsto U_T$ (see \S\,\ref{S:ErgodicTheorem}). Then $(\text{Aut$(X,\mu)$}, d_w)$ is a Polish topological group and Aut$_0(X,\nu)$ is a closed subgroup of Aut$(X,\mu)$.
This topology will not be affected if we replace $\mu$ with any equivalent measure. We note that $T_n$ weakly converges to $T$ if and only if $\mu(T_{n}^{-1}A\bigtriangleup T^{-1}A)\to 0$ for each $A\in\sS$ and
$d(\mu\circ {T_n})/d\mu\to d(\mu\circ T)/d\mu$ in $L^1(X,\mu)$.
For each $p\ge 1$, one can also embed Aut$(X,\mu)$ into the isometry
group of $L^p(X,\mu)$ via a formula similar to (\ref{koopman oper}) but
with another power of the Radon-Nikodym derivative in it. The strong
operator topology on the isometry group induces  the very same weak
topology on Aut$(X,\mu)$ for all $p\ge 1$ \cite{CK79}.
 Danilenko showed in \cite{Dan} that $(\text{Aut$(X,\mu)$}, d_w)$ is contractible. It follows easily from the  Rokhlin lemma that periodic transformations are dense in Aut$(X,\mu)$.

It is natural to ask which properties of nonsingular transformations are typical in the sense of Baire category.  The following technical lemma
(see see \cite{F70}, \cite{CK79}) is an indispensable  tool when considering such problems.

  \begin{lemma} \label{conjugacy}The conjugacy class of each  aperiodic transformation $T$ is dense in {\rm Aut}$(X,\mu)$ endowed with the weak topology.
    \end{lemma}

Using this lemma  and the  Hurewicz ergodic theorem Choksi and  Kakutani \cite{CK79} proved that
    the  ergodic transformations form a  dense
    $ G_\delta$ in Aut$(X,\mu)$.  The same holds for the subgroup Aut$_0(X,\nu)$ (\cite{S71} and \cite{CK79}). Combined with \cite{IT65}
the above implies that the ergodic transformations of type $III$ is a dense $G_\delta$ in Aut$(X,\mu)$. For further refinement of this statement we refer to Section~\ref{orbits}.

Since the map $T\mapsto T\times \cdots\times T$\,($p$ times) from Aut$(X,\mu)$ to Aut$(X^p,\mu^{\otimes p})$ is continuous for each $p>0$, we deduce that the set $\mathcal E_\infty$ of transformations with infinite ergodic index (which means that $T\times \cdots\times T$\,($p$ times) is ergodic for each $p>0$)
is a  $G_\delta$ in Aut$(X,\mu)$. It is non-empty by \cite{KP}. Since this $\mathcal E_\infty$ is invariant under conjugacy, it is dense in Aut$(X,\mu)$ by Lemma~\ref{conjugacy}. Thus we obtain that $\mathcal E_\infty$ is a dense $G_\delta$. In a similar way one can show that
$\mathcal E_\infty\cap\text{Aut}_0(X,\nu)$ is a dense $G_\delta$ in Aut$_0(X,\nu)$ (see also \cite{S71},  \cite{CK79}, \cite{CN00} for original proofs of these claims).

 A nonsingular transformation $T$ is called {\it rigid} if $T^{n_i}\to\text{Id}$ weakly for some sequence $n_k\to\infty$.
 The rigid transformations form a dense $G_\delta$ in Aut$(X,\mu)$.
  It follows  that the set of multiply recurrent nonsingular transformations   is residual \cite{AgSi}. 
 A finer result was established in \cite{DanS}: the set of polynomially recurrent transformations  in 
 Aut$_0(X,\nu)$ is  residual  in Aut$_0(X,\nu)$.
 For the definition of multiple and polynomial recurrence we refer to \S\,\ref{MPR} below.

Given $T\in\text{Aut}(X,\mu)$, we denote the {\it centralizer}  $\{S\in \text{Aut}(X,\mu)\mid ST=TS\}$ of $T$ by $C(T)$. Of course, $C(T)$ is a closed subgroup of Aut$(X,\mu)$ and  $C(T)\supset\{T^n\mid n\in\mathbb Z\}$.
 In a similar way, if
 $T\in \text{Aut$_0(X,\nu)$}$,  the {\it measure preserving centralizer}
$C_0(T):=\text{Aut$_0(X,\nu)$}\cap C(T)$
 of $T$
 is a weakly closed subgroup of \text{Aut$_0(X,\nu)$}.
The following problems solved  (by   several authors) for probability preserving systems are still open for the nonsingular case. Are the  properties:
\begin{itemize}
\item[(i)] $T$ has square root;
\item[(ii)] $T$ embeds into a flow;
\item[(iii)] $T$ has non-trivial invariant sub-$\sigma$-algebra;
\item[(iv)] $C(T)$ contains a torus of arbitrary dimension
\end{itemize}
 typical (residual) in Aut$(X,\mu)$ or  Aut$_0(X,\nu)$?

    The {\it uniform} topology on Aut$(X,\mu)$, finer than $d_w$, is defined by the metric
    \[
    d_u(T,S) =\mu(\{x: Tx\neq Sx\}) +\mu(\{x: T^{-1}x\neq S^{-1}x\}).
    \]
    This topology is  also complete metric.
    It depends only on the measure class of $\mu$. However the uniform topology is not separable and that is why it is of less importance in ergodic theory.
We refer to \cite{ChF}, \cite{F70}, \cite{CK79} and
\cite{CP83} for the properties of $d_u$.

\subsection{}
Suppose now that $\mu(X)=\infty$ but $\mu$ is $\sigma$-finite.
We now let
$$
\begin{aligned}
\text{Aut}_2(X,\mu)&:=\bigg\{T\in\text{Aut}(X,\mu)\mid \sqrt{\frac{d\mu\circ T}{d\mu}}-1\in L^2(\mu)\bigg\}\quad
\text{and}\\
\text{Aut}_1(X,\mu)&:=\{T\in\text{Aut}(X,\mu)\mid \frac{d\mu\circ T}{d\mu}-1\in L^1(\mu)\}.
\end{aligned}
$$
These groups appear naturally in the study of nonsingular Poisson suspensions (see~(\ref{NPT})).
Define a topology $d_2$ on $\text{Aut}_2(X,\mu)$ by setting that a sequence $(T_n)_{n=1}^\infty$
converges to $T$ in $d_2$  if 
$T_n\to T$ weakly and $\|\sqrt{\frac{d\mu\circ T_n}{d\mu}}-\sqrt{\frac{d\mu\circ T}{d\mu}}\|_2\to 0$
as $n\to\infty$.
In a similar way one can define a topology $d_1$ on $\text{Aut}_1(X,\mu)$: 
a sequence $(T_n)_{n=1}^\infty$
converges to $T$ in $d_1$  if $T_n\to T$ weakly and $\|{\frac{d\mu\circ T_n}{d\mu}}-{\frac{d\mu\circ T}{d\mu}}\|_1\to 0$
as $n\to\infty$.

\begin{theorem}\cite{DaKoRo1}
\begin{itemize}
\item $(\text{Aut}_2(X,\mu), d_2)$ is a Polish group.
\item $(\text{Aut}_1(X,\mu), d_1)$ is a Polish group.
\item The mapping $\chi:\text{Aut}_1(X,\mu)\ni T\mapsto\int_X(\frac{d\mu\circ T}{d\mu}-1)d\mu\in\Bbb R$ is a 
 continuous onto homomorphism.
 \item
 The short exact sequence
 $$
 \{1\}\to\ker\chi\to\text{Aut}_1(X,\mu) \overset{\chi}\to\Bbb R\to\{1\}
 $$
 splits.
 \item $\text{Aut}_1(X,\mu)$ is a dense and meager  subgroup of $(\text{Aut}_2(X,\mu),d_2)$.
 \item 
 $\chi$ is not $d_2$-continuous.
 \item 
 the group $\{T_*\mid T\in \text{Aut}_2(X,\mu)\}$ is weakly closed in $\text{Aut}(X^*,\mu^*)$.
 \item
 the groups $ (\text{Aut}_2(X,\mu),d_2)$ and $(\ker\chi,d_1)$ have the Rokhlin property (i.e. there exists a dense conjugacy class in each of these groups).
 \item
 The group  $ (\text{Aut}_1(X,\mu),d_1)$ does not have the Rokhlin property.
\end{itemize}
\end{theorem}

\section{Orbit theory }\label{orbits}
Orbit theory is, in a sense, the most complete part of  nonsingular ergodic theory.  We present here the  seminal Krieger's  theorem on orbit classification of ergodic nonsingular transformations in terms of  ratio sets and  associated flows. Examples of transformations of various types $III_\lambda$, $0\le \lambda\le 1$ are  given.
``Almost continuous'' refinement of the orbit equivalence  is also considered here.
Next, we consider the outer conjugacy problem for automorphisms of the orbit equivalence relations. This problem is solved  in terms of a simple complete system of invariants. We discuss also a general theory of cocycles (of nonsingular systems) taking values in locally compact Polish groups and present an important orbit classification theorem for cocycles. This theorem is an analogue of the aforementioned result of Krieger.
We complete the section by considering ITPFI-systems and their relation to AT-flows.

\subsection{Full groups.  Ratio set and types $III_\lambda$,
$0\le \lambda\le 1$
}
Let $T$ be a nonsingular transformation of a standard probability space
$(X,\mathcal B,\mu)$. Denote by Orb$_T(x)$ the $T$-orbit of $x$, i.e.,
Orb$_T(x)=\{T^nx\mid n\in\mathbb Z\}$. The {\it full group} $[T]$ of $T$
consists of all transformations $S\in\text{Aut}(X,\mu)$ such that
$Sx\in\text{Orb}_T(x)$ for a.a. $x$. If $T$ is ergodic  then $[T]$ is topologically simple (or even algebraically simple if $T$ is not of type $II_\infty$) \cite{Eig}. It is easy to see
that $[T]$ endowed with the uniform topology $d_u$ is a Polish group. If $T$ is ergodic then $([T],d_u)$ is
contractible \cite{Dan}.

The {\it ratio set\,} $r(T)$ of $T$ was  defined by Krieger  [Kr70] and as we shall see below it is the key concept in the orbit
classification (see Definition~\ref{orbit equivalent}). The ratio set is a subset of $[0,+\infty)$ defined as follows: $t\in
r(T)$ if and only if for every $A\in\mathcal B$ of positive measure and each
$\epsilon>0$ there is a subset $B\subset A$ of positive measure and an
integer $k\ne 0$ such that $T^kB\subset A$ and $| \omega_k^\mu(x)-t|<\epsilon$ for all $x\in B$. It is easy to verify that
$r(T)$ depends only on the equivalence class of $\mu$ and not on $\mu$
itself.
A basic
fact is that $1\in r(T)$ if and only if  $T$ is conservative. Assume now  $T$ to be
conservative and ergodic. Then $r(T)\cap(0,+\infty)$ is a closed subgroup
of the multiplicative group $(0,+\infty)$. Hence $r(T)$ is one of the
following sets:
\begin{enumerate}
\item[(i)]
$\{1\}$;
\item[(ii)]
$\{0,1\}$; in this case we say that $T$ is of {\it type $III_0$},
\item[(iii)]
$\{\lambda^n\mid n\in\mathbb Z\}\cup\{0\}$ for $0<\lambda<1$; then we say that  $T$ is  of {\it type $III_\lambda$},
\item[(iv)]
$[0,+\infty)$; then we say that  $T$ is  of {\it type III$_1$}.
\end{enumerate}
Krieger showed that $r(T)=\{1\}$ if and only if $T$ is of type $II$.
 Hence we obtain a further subdivision    of  type $III$ into subtypes  $III_0$, $III_\lambda$, $0<\lambda<1$, and $III_1$.

\begin{ex}\label{type III}
(i)  Fix $\lambda\in(0,1)$. Let
$\nu_n(0):=1/(1+\lambda)$ and $\nu_n(1):=\lambda/(1+\lambda)$ for all
$n=1,2,\dots$. Let $T$ be the nonsingular product odometer associated with the sequence $(2,\nu_n)_{n=1}^\infty$  (see \S\,\ref{SS:odometerdef}). We claim that $T$ is of type $III_\lambda$.
Indeed, the group $\Sigma$ of finite permutations of $\mathbb N$ acts on $X$
by $(\sigma x)_n=x_{\sigma^{-1}(n)}$, for all $n\in\mathbb N$,
$\sigma\in\Sigma$ and $x=(x_n)_{n=1}^\infty\in X$. This action preserves
$\mu$. Moreover, it is ergodic by the Hewitt-Savage 0-1 law. It remains to
notice that $(d\mu\circ T/d\mu)(x)=\lambda$ on the cylinder $[0]$ which is of positive measure.

(ii)  Fix positive reals $\rho_1$ and $\rho_2$ such
that $\log\rho_1$ and $\log\rho_2$ are rationally independent. Let
$\nu_n(0):=1/(1+\rho_1+\rho_2)$, $\nu_n(1):=\rho_1/(1+\rho_1+\rho_2)$ and
$\nu_n(2):=\rho_2/(1+\rho_1+\rho_2)$ for all $n=1,2,\dots$. Then the
nonsingular product odometer associated with the sequence
$(3,\nu_n)_{n=1}^\infty$ is of type $III_1$. This can be shown in a similar way as (i).

(iii) Partition $\Bbb N$ into two infinite subsets $A$ and $B$.
Fix a sequence $(\epsilon_n)_{n\in B}$ of positive reals such that $\epsilon_n<0.4$ for all $n\in B$ and $ \sum_{n\in B}\epsilon_n<\infty$.
We now let $\nu_n(0):=0.5$ if $n\in A$ and $\mu_n(0):=1-\epsilon_n$ if $n\in B$.
Then the
nonsingular product odometer associated with the sequence
$(2,\nu_n)_{n=1}^\infty$ is of type $II_\infty$.
This follows from \cite{Moo}.
\end{ex}

Non-singular product odometer of type $III_0$ will be constructed in Example~\ref{IIIzero}
below.

\subsection{Maharam extension, associated flow and orbit classification of
type $III$ systems
}\label{S:maharam}
On $X\times\mathbb R$ with the $\sigma$-finite measure $\mu\times\kappa$,
where $d\kappa(y)=\exp(y)dy$, consider the transformation
$$
\widetilde T(x,y):=(Tx,y-\log\frac{d\mu\circ T}{d\mu}(x)).
$$
We call it the {\it Maharam extension}  of $T$ (see \cite{Mah}, where these
transformations were introduced). It is measure-preserving and it commutes
with the flow $S_t(x,y):=(x,y+t)$, $t\in\mathbb R$. 
It is conservative if and
only if $T$ is conservative \cite{Mah}.
 However $\widetilde T$ is not
necessarily ergodic when $T$ is ergodic.
 Let $(Z,\nu)$ denote the space of $\widetilde T$-ergodic
components. Then $(S_t)_{t\in\mathbb R}$ acts nonsingularly
on this space. The restriction of $(S_t)_{t\in\mathbb R}$ to $(Z,\nu)$ is
called the {\it associated flow} of $T$. The associated flow is ergodic
if and only if  $T$ is ergodic. It is easy to verify that the isomorphism class of
the associated flow is an invariant of the orbit equivalence of the
underlying system.

\begin{prop}[\cite{HO81}]
\begin{enumerate}
\item[(i)]
$T$ is of type $II$ if and only if its associated flow is the translation
on $\mathbb R$, i.e., $x\mapsto x+t$, $x,t\in\mathbb R$,
\item[(ii)]
$T$ is of type $III_\lambda$, $0\le\lambda<1$ if and only if its associated
flow is the periodic flow on the interval $[0,-\log\lambda)$, i.e.,
$x\mapsto x+t\mod(-\log\lambda)$,
\item[(iii)]
$T$ is of type $III_1$ if and only if its associated flow is the trivial
flow on a singleton or, equivalently, $\widetilde T$ is ergodic,
\item[(iv)]
$T$ is of type $III_0$ if and only if its associated flow is nontransitive.
\end{enumerate}
\end{prop}

\begin{ex}\label{IIIzero}  Let $A_n=\{0,1,\dots,2^{2^n}\}$ and $\nu_n(0)=0.5$
and $\nu_n(i)=0.5\cdot 2^{-2^n}$ for all $0<i\le 2^n$. Let $T$ be the
nonsingular product odometer associated with $(2^{2^n}+1,\nu_n)_{n=0}^\infty$.
 It is straightforward that the associated flow of $T$ is the
flow built under the constant function $1$ with the probability preserving 2-adic product odometer (associated with $(2,\kappa_n)_{n=1}^\infty$, $\kappa_n(0)=\kappa_n(1)=0.5$) as the
base transformation. In particular, $T$ is of type $III_0$.
\end{ex}

A natural problem arises: to compute Krieger's type (or the ratio set) for the nonsingular product odometers---the simplest class of nonsingular systems. Some partial progress was achieved in \cite{AW}, \cite{Moo}, \cite{Os2},  \cite{DKQ}, etc. However in the general setting this problem  remains open.

The map $\Psi:\text{Aut}(X,\mu)\ni T\mapsto\widetilde
T\in\text{Aut}(X\times\mathbb R,\mu\times\kappa)$ is a continuous group
homomorphism. Since the set $\mathcal E$ of ergodic
transformations on $X\times\mathbb R$ is a $G_\delta$ in Aut$(X\times\mathbb
R,\mu\times\kappa)$ (See \S\,\ref{topological group}), the subset $\Psi^{-1}(\mathcal E)$
of type $III_1$ ergodic transformations on $X$ is also $G_\delta$. The
latter subset is non-empty in view of  Example~\ref{type III}(ii). Since it is invariant under
conjugacy, we deduce from Lemma~\ref{conjugacy} that the set of ergodic transformations of type $III_1$ is a dense $G_\delta$
in (Aut$(X,\mu),d_w)$  (\cite{PaS}, \cite{CHP}).

Now we state the main result of this section---Krieger's theorem on orbit
classification for ergodic transformations of type $III$. It is a far
reaching generalization of the basic result by H. Dye: any two ergodic
probability preserving transformations are orbit equivalent \cite{dye}.

\begin{theorem}[Orbit equivalence for type $III$ systems
\cite{Kr1}---\cite{Kr76}]\label{T:2.4}  Two ergodic transformations of type $III$ are
orbit equivalent if and only if their associated flows are isomorphic. In
particular, for a fixed
$0<\lambda\le 1$, any two ergodic transformations of type $III_\lambda$  are orbit equivalent.
\end{theorem}

The original proof of this theorem is rather complicated. Simpler treatment
of it can be found in \cite{HO81} and \cite{KW}.

We also note that every free  ergodic flow can be realized as the
associated flow of a type $III_0$ transformation. 
However it is somewhat
easier to construct a $\mathbb Z^2$-action of type $III_0$ whose associated
flow is the given one. For this, we take an ergodic nonsingular
transformation $Q$ on a probability space $(Z,\mathcal B,\lambda)$ and a
measure-preserving transformation $R$ of an infinite $\sigma$-finite
measure space $(Y,\mathcal F,\nu)$ such that there is a continuous
homomorphism $\pi:\mathbb R\to C(R)$ with
$(d\nu\circ\pi(t)/d\nu)(y)=\exp(t)$ for a.a. $y$ (for instance, take a type
$III_1$ transformation $T$ and put $R:=\widetilde T$ and $\pi(t):=S_t$).
Let $\phi:Z\to \mathbb R$ be a Borel map with $\inf_Z\phi>0$. Define two
transformations $R_0$ and $Q_0$ of $(Z\times Y,\lambda\times\nu)$ by
setting:
$$
R_0(x,y):=(x,Ry), \ \ Q_0(x,y)=(Qx, U_xy),
$$
where $U_x=\pi(\phi(x)-\log(d\mu\circ Q/d\mu)(x))$. Notice that $R_0$ and
$Q_0$ commute. The corresponding $\mathbb Z^2$-action generated by these
transformations is ergodic. Take any transformation $V\in
\text{Aut}(Z\times Y,\lambda\times\nu)$ whose orbits coincide with the
orbits of the $\mathbb Z^2$-action. (According to \cite{CFW}, any ergodic
nonsingular action of any countable amenable group is orbit equivalent to
a single transformation.)  It is now easy to
verify that the associated flow of $V$ is the special flow built under
$\phi\circ Q^{-1}$ with the base transformation $Q^{-1}$. 
Then $V$ is of type $III_0$.
Since $Q$ and
$\phi$ are arbitrary, we deduce the following from Theorem~\ref{ambrose}.

\begin{theorem} \label{T:2.5}  Every nontransitive ergodic flow is an associated flow of
an ergodic transformation of type $III_0$.
\end{theorem}

In \cite{Kr76} Krieger introduced a map $\Phi$ as follows. Let $T$ be an
ergodic transformation of type $III_0$. Then the associated flow of $T$ is
a flow built under function with a base transformation $\Phi(T)$. We note
that the orbit equivalence class of $\Phi(T)$ is well defined by the orbit
equivalent class of $T$. If $\Phi^n(T)$ fails to be of type $III_0$ for
some $1\le n<\infty$ then $T$ is said to {\it belong to Krieger's
hierarchy}. For instance, the transformation constructed in Example~\ref{IIIzero} belongs to Krieger's hierarchy. Connes gave in \cite{Co} an example of $T$ such that $\Phi(T)$
is orbit equivalent to $T$ (see also \cite{HO81} and \cite{GiS}). Hence $T$
is not in
 Krieger's hierarchy.

 \subsection{Almost continuous orbit equivalence}
 In this subsection, by a dynamical system we mean a quadruple $(X,\tau,\mu,T)$, where $(X,\tau)$ is a Polish space, $\mu$ is a non-atomic Borel measure of  full support, $T$ is a nonsingular  ergodic homeomorphism of $X$ such that  the function $\omega_1:X\to\Bbb R$ is continuous (has a continuous version).

\begin{defn} Two dynamical systems $(X,\tau,\mu,T)$ and $(X',\tau',\mu',T')$ are {\it almost continuously orbit equivalent} if there are dense invariant $G_\delta$ subsets $X_0\subset X$ and $X_0'\subset X'$ of full measure and a homeomorphism $\phi:X_0\to X_0'$ such that
\begin{itemize}
\item $\phi(\{T^nx\mid n\in\Bbb Z\})=\{(T')^n\phi(x)\mid n\in\Bbb Z\}$ at every $x\in X_0$,
\item
$\mu\circ\phi^{-1}\sim\mu'$ and the Radon-Nikodym derivative
$\frac{d\mu\circ\phi^{-1}}{d\mu'}$ is (can be chosen) continuous,
\item
letting $S:=\phi^{-1}T'\phi$ we have $Tx=S^{n(x)}x$ and $Sx=T^{m(x)}x$, where $n$ and $m$ are continuous on $X_0$.
\end{itemize}
 \end{defn}
 
We note that in the case where $X$ and $X'$ are infinite product spaces, $T$ and $T'$ preserve $\mu$ and $\mu'$ respectively and we omit the requirement that $X_0$ and $X_0$ are $G_\delta$ then the above definition of $\phi$ is equivalent to the ``finitary'' equivalence from the celebrated work of Keane and Smorodinsky \cite{KS}.
It was shown by del Junco and  {\c S}ahin \cite{dJS}  that any two ergodic probability preserving homeomorphisms
of Polish spaces are almost continuously orbit equivalent. 
 The same is
true for
any ergodic homeomorphisms preserving infinite
$\sigma$-finite local measures \cite{dJS}.
In \cite{DadJ}, a topological analogue
$r_{\text{top}}(T)$ of $r(T)$ was introduced.
 It is a closed subgroup of
$\Bbb R$
which contains $r(T)$ and it is invariant under the almost continuous orbit equivalence. 
In   \cite{DadJ}, two  type
$III$
homeomorphisms    were constructed which are measure-theoretically orbit equivalent but not almost continuously orbit equivalent (their $r_{\text{top}}$-invariants are different).

\begin{theorem}[\cite{DadJ}]
Let
$(X,\tau,\mu,T)$ and $(X',\tau',\mu',T')$
be ergodic non-singular homeomorphisms of Polish spaces. 
If the two systems are either
\begin{itemize}
\item[{\rm (i)}]
of type
$III_\lambda$ with $0<\lambda<1$ 
and $r_{\text{top}}(T)=r_{\text{top}}(T')=\log\lambda\cdot\Bbb Z$
or
\item[{\rm (ii)}]
of type
$III_1$
\end{itemize}
then they are almost continuously orbit equivalent.
\end{theorem}

 Characterization  of almost continuous orbit equivalence for  homeomorphisms of type $III_0$  remains an open problem.

 \subsection{Normalizer of the full group. Outer conjugacy problem
}
Let
$$
N[T]=\{R\in\text{Aut}(X,\mu)\mid R[T]R^{-1}=[T]\},
$$
i.e., $N[T]$ is the {\it normalizer} of the full group $[T]$ in
Aut$(X,\mu)$. We note that a transformation $R$ belongs to $N[T]$ if and
only if $R(\text{Orb}_T(x))=\text{Orb}_T(Rx)$ for a.a. $x$. To define a
topology on $N[T]$ consider the $T$-orbit equivalence relation $\mathcal
R_T\subset X\times X$ and a $\sigma$-finite measure $\mu_\mathcal R$ on $\mathcal
R_T$ given by $\mu_{\mathcal
R_T}=\int_X\sum_{y\in\text{Orb}_T(x)}\delta_{(x,y)}d\mu(x)$. For $R\in
N[T]$, we define a transformation $i(R)\in\text{Aut}(\mathcal R_T,\mu_{\mathcal
R_T})$ by setting $i(R)(x,y):=(Rx,Ry)$. Then the map $R\mapsto i(R)$ is an
embedding of $N[T]$ into Aut$(\mathcal R_T,\mu_{\mathcal R_T})$. Denote by $\tau$
the topology on $N[T]$ induced by the weak topology on Aut$(\mathcal
R_T,\mu_{\mathcal R_T})$ via $i$ \cite{Dan}. Then $(N[T],\tau)$ is a Polish group. A
sequence $R_n$  converges to $R$ in $(N[T],\tau)$ if $R_n\to R$ weakly (in
Aut$(X,\mu)$) and $R_nTR_n^{-1}\to RTR^{-1}$ uniformly (in $[T]$).

Given $R\in N[T]$, denote by $\widetilde R$ the Maharam extension of $R$.
Then $\widetilde R\in N[\widetilde T]$ and it commutes with $(S_t)_{t\in
\mathbb R}$. Hence it defines a nonsingular transformation mod\,$R$ on the
space $(Z,\nu)$ of the associated flow $W=(W_t)_{t\in\mathbb R}$ of $T$.
Moreover, mod\,$R$ belongs to the centralizer $C(W)$ of $W$ in
Aut$(Z,\nu)$. Note that $C(W)$ is a closed subgroup of
(Aut$(Z,\nu),d_w)$.

Let $T$ be of type $II_\infty$ and let $\mu'$ be the invariant
$\sigma$-finite measure equivalent to $\mu$. If $R\in N[T]$ then it is easy
to see that the Radon-Nikodym derivative $d\mu'\circ R/d\mu'$ is invariant
under $T$. Hence it is constant, say $c$. Then mod$\,R=\log c$.

\begin{theorem}[\cite{HO81}, \cite{Ham}]\label{T:2.6} If $T$ is of type III then
the map $\mod:N[T]\to C(W)$ is a continuous onto homomorphism. The kernel
of this homomorphism is the $\tau$-closure of $[T]$. Hence the quotient
group $N[T]/\overline{[T]}^\tau$ is (topologically) isomorphic to $C(W)$.
In particular, $\overline{[T]}^\tau$ is co-compact in $N[T]$ if and only if
 $W$ is a finite measure-preserving flow with a pure point spectrum.
\end{theorem}

The following theorem describes the homotopical structure of normalizers.

\begin{theorem}[\cite{Dan}]\label{T:2.7}  Let $T$ be of type $II$ or
$III_\lambda$, $0\le\lambda<1$. The group $\overline{[T]}^\tau$ is
contractible. $N[T]$ is homotopically equivalent to $C(W)$. In particular,
$N[T]$ is contractible if $T$ is of type $II$. 
If $T$ is of type
$III_\lambda$ with $0<\lambda<1$ then $\pi_1(N[T])=\mathbb Z$.
\end{theorem}

The {\it outer period} $p(R)$ of $R\in N[T]$ is the smallest positive
integer $n$ such that $R^n\in[T]$. We write $p(R)=0$ if no such $n$ exists.

Two transformations $R$ and $R'$ in $N[T]$ are called {\it outer conjugate}
if there are transformations $V\in N[T]$ and $S\in[T]$ such that
$VRV^{-1}=R'S$. The following theorem provides convenient (for
verification) necessary and sufficient conditions for the outer conjugacy.

\begin{theorem}[\cite{CK} for type $II$ and \cite{BG} for type
$III$]\label{T:2.8}  Transformations $R,R'\in N[T]$ are outer conjugate if and only if
$p(R)=p(R')$ and $\mod R$ is conjugate to $\mod R'$ in the centralizer of
the associated flow of $T$.
\end{theorem}

We note that in the case $T$ is of type $II$, the second condition in the
theorem is just mod\,$R=\text{mod\,} R'$. 
It is always satisfied when $T$ is of
type $II_1$.

\subsection{Cocycles of dynamical systems. Weak equivalence of cocycles
}\label{cocycles}
Let $G$ be a locally compact Polish group and $\lambda_G$ a left Haar
measure on $G$. A Borel map $\phi:X\to G$ is called a {\it cocycle} of $T$.
 Two cocycles $\phi$ and $\phi'$ are
{\it cohomologous} if there is a Borel map $b:X\to G$ such that
$$
\phi'(x)=b(Tx)^{-1}\phi(x)b(x)
$$
for a.a. $x\in X$. A cocycle cohomologous to the trivial one is called a {\it
coboundary}. Given a dense subgroup $G'\subset G$, then every cocycle is
cohomologous to a cocycle with values in $G'$ \cite{GS}. Each cocycle
$\phi$ extends to a (unique) map $\alpha_\phi:\mathcal R_T\to G$ such that
$\alpha_\phi(Tx,x)=\phi(x)$ for a.a. $x$ and
$\alpha_\phi(x,y)\alpha_\phi(y,z)=\alpha_\phi(x,z)$ for a.a.
$(x,y),(y,z)\in\mathcal R_T$. $\alpha_\phi$ is called the {\it cocycle of $\mathcal
R_T$ generated by $\phi$}. Moreover, $\phi$ and $\phi'$ are cohomologous
via $b$ as above if and only if $\alpha_\phi$ and $\alpha_{\phi'}$ are {\it
cohomologous} via $b$, i.e., $
\alpha_\phi(x,y)=b(x)^{-1}\alpha_{\phi'}(x,y)b(y) $ for $\mu_{\mathcal
R_T}$-a.a. $(x,y)\in \mathcal R_T$. 
The following notion was introduced by Golodets and Sinelshchikov \cite{GS3}, \cite{GS}: two cocycles $\phi$ and $\phi'$ are {\it
weakly equivalent} if there is a transformation $R\in N[T]$ such that the
cocycles $\alpha_\phi$ and $\alpha_\phi'\circ (R\times R)$ of $\mathcal R_T$ are
cohomologous. 
Let $\mathcal M(X,G)$ denote the set of Borel maps from $X$ to
$G$. It is a Polish group when endowed with the topology of convergence in
measure. Since $T$ is ergodic, it is easy to deduce from Rokhlin's lemma
that the cohomology class of any cocycle is dense in $\mathcal M(X,G)$. Given
$\phi\in\mathcal M(X,G)$, we define the $\phi$-{\it skew product extension}
$T_\phi$ of $T$ acting on $(X\times G,\mu\times\lambda_G)$ by setting
$T_\phi(x,g):=(Tx,\phi(x)g)$. Thus Maharam extension is (isomorphic to) the
Radon-Nikodym cocycle-skew product extension.
We now specify some basic classes of
cocycles \cite{Sc},  \cite{BG2}, \cite{GS}, \cite{Dan2}:
\begin{enumerate}
\item[(i)]  $\phi$ is called {\it transient} if $T_\phi$ is totally disipative,
\item[(ii)]  $\phi$ is called {\it recurrent} if $T_\phi$ is conservative
 (equivalently, $\phi$ is not transient),
\item[(iii)]  $\phi$ {\it has dense range in $G$} if $T_\phi$ is
ergodic.
\item[(iv)]  $\phi$ is called {\it regular} if $\phi$ cobounds with dense
range into a closed subgroup $H$ of $G$ (then $H$ is defined up to
conjugacy).
\end{enumerate}
These properties are invariant under the cohomology and the weak
equivalence. The Radon-Nikodym cocycle $\omega_1$ is a
coboundary if and only if $T$ is of type $II$. It is regular if and only if
$T$ is of type $II$ or $III_\lambda$, $0<\lambda\le 1$. It has dense range
(in the multiplicative group $\mathbb R_+^*$) if and only if $T$ is of type
$III_1$. Notice that $\omega_1$
 is never transient (since $T$ is conservative).

In case $G$ is Abelian, Schmidt introduced in \cite{Sc5} an invariant $R(\phi):=\{g\in G\mid \phi-g\text{\ is recurrent}\}$. He showed in particular that
\begin{itemize}
\item[(i)]
$R(\phi)$ is a cohomology invariant,
\item[(ii)]
 $R(\phi)$ is a Borel set in $G$,
\item[(iii)] $R(\log\omega_1)=\{0\}$ for each aperiodic conservative $T$,
\item[(iv)]  there are cocycles $\phi$ such that $R(\phi)$ and $G\setminus R(\phi)$ are dense in $G$,
\item[(v)] if $\mu(X)=1$, $\mu\circ T=\mu$ and $\phi:X\to\mathbb R$ is integrable then $R(\phi)=\{\int\phi\, d\mu\}$.
\end{itemize}
We note that (v) follows from Atkinson theorem \cite{Atk}. A nonsingular version of this theorem was established in \cite{Ull}: if $T$ is ergodic and $\mu$-nonsingular and $f\in L^1(\mu)$ then $$
\liminf_{n\to\infty}\bigg|\sum_{j=0}^{n-1}f(T^jx)\omega_j(x)\bigg|=0\text{ \ for a.a. }x
$$
if and only if
$\int f\,d\mu=0$.

Since $T_\phi$ commutes with the action of $G$ on $X\times G$ by inverted
right translations along the second coordinate, this action induces an ergodic $G$-action
$W_\phi=(W_\phi(g))_{g\in G}$ on the space $(Z,\nu)$ of
$T_\phi$-ergodic components. It is called the {\it Mackey range (or Poincar{\'e} flow)} of $\phi$
\cite{Mac}, \cite{FM}, \cite{Sc}, \cite{Zim3}. We note that $\phi$ is regular (and cobounds with
dense range into $H\subset G$) if and only if $W_\phi$ is transitive (and
$H$ is the stabilizer of a point $z\in Z$, i.e., $H=\{g\in G\mid
W_\phi(g)z=z\}$). Hence every cocycle taking values in a compact
group is regular.

It is often useful to consider the {\it double cocycle} $\phi_0:=\phi\times
\omega_1$ instead of $\phi$. It takes values in the group $G\times\mathbb
R^*_+$. Since $T_{\phi_0}$ is exactly the Maharam extension of $T_\phi$, it
follows from \cite{Mah} that $\phi_0$ is transient or recurrent if and only
if $\phi$ is transient or recurrent respectively.

\begin{theorem}[Orbit classification of cocycles \cite{GS}]\label {T:2.9} Let
$\phi,\phi':X\to G$ be two recurrent cocycles of an ergodic transformation
$T$. They are weakly equivalent if and only if their Mackey ranges
$W_{\phi_0}$ and $W_{\phi_0'}$ are isomorphic.
\end{theorem}

 Another proof of this theorem was  presented in \cite{Fed}.

\begin{theorem}\label{T:2.10} Let $T$ be an ergodic nonsingular transformation.
Then there is a cocycle of $T$ with dense range in $G$ if and only if $G$ is
amenable.
\end{theorem}

It follows that if $G$ is amenable then the subset of cocycles of $T$ with dense
range in $G$ is a dense $G_\delta$ in $\mathcal M(X,G)$ (just adapt the
argument following Example~\ref{IIIzero}). The `only if' part of Theorem~\ref{T:2.10} was
established in \cite{Zim}. The `if' part was considered by many authors in
particular cases: $G$ is compact \cite{Zim2}, $G$ is solvable or amenable
almost connected \cite{GS2}, etc. The
general case was proved in \cite{GS3} and \cite{Her} (see also a recent
treatment in \cite{AaW}).

We note that the ``if'' part in Theorem~\ref{T:2.10} can be refined in the case where $G$ is a compactly generated Abelian group.

\begin{theorem}[\cite{Dan21}]
Let $T$ be an ergodic nonsingular transformation.
If $G$ is a compactly generated [FIA]-group then there is a bounded ergodic cocycle $\varphi$ of $T$ with values in $G$.
\end{theorem}

We recall that a locally compact Polish group $G$ is called a {\it [FIA]-group} if 
the group of inner automorphisms of $G$ is relatively compact in the group of all automorphisms of $G$ furnished with the natural topology \cite{GrMo}.
Of course, each Abelian group is [FIA].
The cocycle $\phi$ is {\it bounded} if there is a compact subset $K$ in $G$ such that $\varphi$ takes values in $K$.

  Theorem~\ref{T:2.5} is a particular case of the following result.

\begin{theorem}[\cite{GS4}, \cite{Fed}, \cite{AEG}]\label{T:2.11} Let $G$ be amenable.
Let $V$ be an ergodic nonsingular action of $G\times\mathbb R_+^*$.
Then there is an ergodic nonsingular transformation $T$ and a recurrent cocycle
$\phi$ of $T$ with values in $G$ such that $V$ is isomorphic to the Mackey
range of the double cocycle $\phi_0$.
\end{theorem}

Given a cocycle $\phi\in\mathcal M(X,G)$ of $T$, we say that a transformation
$R\in N[T]$ is {\it compatible with} $\phi$ if the cocycles $\alpha_\phi$
and $\alpha_\phi\circ(R\times R)$ of $\mathcal R_T$ are cohomologous. Denote by
$D(T,\phi)$ the group of all such $R$. It has a natural Polish topology
which is stronger than $\tau$ \cite{DaG}. Since $[T]$ is a normal subgroup
in $D(T,\phi)$, one can consider the outer conjugacy equivalence relation
inside $D(T,\phi)$. It is called {\it $\phi$-outer conjugacy}. Suppose that
$G$ is Abelian. Then an analogue of Theorem~\ref{T:2.8} for the $\phi$-outer
conjugacy is established in \cite{DaG}.
 Also, the cocycles $\phi$ with $D(T,\phi)=N[T]$ are
described there.

\subsection{ITPFI transformations and AT-flows}\label{AT-f}
A nonsingular transformation $T$ is called {\it ITPFI\footnote{This
abbreviates `infinite tensor product of factors of type $I$' (came from the
theory of von Neumann algebras).}} if it is orbit equivalent to a nonsingular product odometer
(associated to a sequence $(m_n,\nu_n)_{n=1}^\infty$, see \S\,\ref{SS:odometerdef}). If the sequence
$m_n$ can be chosen bounded then $T$ is called {\it ITPFI of bounded type}. If $m_n=2$ for all $n$
 then $T$ is called {\it ITPFI$_2$}. By \cite{GiS2}, every ITPFI-transformation of bounded type is
ITPFI$_2$. 
In view of Theorem~\ref{T:2.4} and Example~\ref{type III}, 
every ergodic transformation of type $II$ or $III_\lambda$ with $0<\lambda\le 1$
 is ITPFI$_2$.

A remarkable characterization of ITPFI transformations in terms
of their associated flows was obtained by Connes and Woods \cite{CoW}. We
first single out a class of ergodic flows. A nonsingular flow
$V=(V_t)_{t\in\mathbb R}$ on a space $(\Omega,\nu)$ is called {\it approximate
transitive (AT)} if given $\epsilon>0$ and $f_1,\dots,f_n\in L^1_+(X,\mu)$,
there exists $f\in L^1_+(X,\mu)$ and $\lambda_1,\dots,\lambda_n\in
L^1_+(\mathbb R,dt)$ such that
$$
\left|\left|f_j-\int_\mathbb R f\circ V_t\frac{d\nu\circ
V_t}{d\nu}\lambda_j(t)dt\right|\right|_1<\epsilon
$$
for all $1\le j\le n$. A flow built under a constant ceiling function with
a funny rank-one \cite{Fer} probability preserving base transformation is AT \cite{CoW}. In
particular, each ergodic finite measure-preserving flow with a pure point
spectrum is AT.

\begin{theorem}[\cite{CoW}] \label{2.12}
 An ergodic nonsingular transformation
is ITPFI if and only if its associated flow is AT.
\end{theorem}

The original proof of this theorem was given in the framework of von
Neumann algebras theory. 
A simpler, purely measure theoretical proof was
given later in \cite{Haw} (the `only if' part) and \cite{Ham2} (the `if'
part). 
It follows from Theorem~\ref{2.12} that every ergodic flow with pure point
spectrum is the associated flow of an ITPFI transformation. 
This was refined recently 
 in \cite{BeVa2}:  every ergodic flow with a pure point
spectrum is the associated flow of an ITPFI$_2$ transformation. 
This fact was proved earlier in  \cite{HO2} only for 
flows whose spectrum  is $\theta\Gamma$, where $\Gamma$ is a subgroup of $\mathbb Q$
and $\theta\in\mathbb R$.
The existence of  ITPFI transformations which are not of bounded type was shown in
\cite{Kr5}. 

Krieger  introduced   an invariant for the orbit equivalence, called {\it property A},
and showed that each product odometer satisfies  property A.
He also constructed an ergodic  nonsingular  transformation which does not satisfy this property \cite{Kr5}. 
Hence this transformation is not ITPFI.
Though not every ergodic transformation is orbit equivalent to a nonsingular product odometer,  a ``weaker'' form of this statement holds.

 \begin{theorem}[\cite{DoH1}]\label{2.13} Each ergodic nonsingular transformation
is orbit equivalent to a Markov odometer (see \S\ref{MarkOd}).
\end{theorem}

In \cite{DoH2}, an explicit example of a non-ITPFI ergodic Markov odometer (not satisfying property A) was
constructed.
Later Munteanu  in  \cite{Mu12} exhibited  an ergodic non-ITPFI transformation satisfying property A.
In \cite{JM}, it was constructed an explicit example of a non-AT nonsingular flow $W$ built under a function and over a nonsingular product odometer.
Hence every nonsingular ergodic transformation whose associated flow is isomorphic to $W$
is non-ITPFI.

\section{Mixing notions and multiple recurrence}\label{mixing}
The study of mixing and multiple recurrence are central topics in  classical ergodic theory \cite{CFS}, \cite{Fu}. Unfortunately, these notions are considerably less `smooth' in the world of nonsingular systems. The very concepts of any kind of mixing and multiple recurrence are not well understood in view of their ambiguity. Below we discuss nonsingular systems possessing a surprising diversity of such properties that seem equivalent but are different indeed.

\subsection{Weak mixing}\label{S:weakmixing}
Let $T$ be an ergodic conservative nonsingular transformation. A number
$\lambda\in\mathbb C$ is an {\it $ L^\infty $-eigenvalue} for $T$ if there
exists a nonzero $f\in L^\infty$ so that $f\circ T=\lambda f$ a.e. It
follows that $|\lambda|=1$ and $f$ has constant modulus, which we assume to
be $1$. Denote by $e(T)$ the set of all $L^\infty$-eigenvalues of $T$. $T$
is said to be {\it weakly mixing} if $e(T)=\{1\}$. We refer to \cite[
Theorem~2.7.1]{Aa} for proof of the following Keane's ergodic
multiplier theorem: given an
ergodic probability preserving transformation $S$, the product
transformation $T\times S$ is ergodic if and only if $\sigma_S(e(T))=0$,
where $\sigma_S$ denotes the measure of (reduced) maximal spectral type of
the unitary $U_S$ (see (\ref{koopman oper})). It follows that $T$ is weakly
mixing if and only $T\times S$ is ergodic for every ergodic probability
preserving $S$.
 While in the finite measure-preserving case this implies that $T\times T$ is ergodic, it was shown in \cite{ALW} that there exits a weakly mixing nonsingular $T$  with $T\times T$ not conservative, hence not ergodic.  In  \cite{AFS},   a weakly mixing $T$ was constructed with $T\times T$ conservative but not ergodic.  A nonsingular transformation $T$ is said to be {\it weakly doubly ergodic} (originally called {\it doubly ergodic} in \cite{BFMS01}) if for all sets
of positive measure $A$ and $B$ there exists an integer $n>0$ such that $\mu(A\cap T^{-n}A)>0$ and $\mu(A\cap T^{-n}B)>0$. 
Furstenberg \cite{Fu} showed that for finite measure-preserving transformations weak double ergodicity is equivalent to weak mixing. 
In \cite{AFS,BFMS01} it is shown that for nonsingular transformations weak mixing does not imply weak double ergodicity and weak double ergodicity does not imply that $T\times T$ is ergodic.

 $T$ is said to have {\it ergodic index $k$} if the Cartesian product of $k$ copies of $T$ is ergodic but the product of $k+1$ copies of $T$ is not ergodic. If all finite Cartesian products of $T$ are ergodic  then $T$ is said to have {\it infinite ergodic index}.
 In a similar way one can define the {\it conservative index} of $T$.
  Parry and Kakutani  \cite{KP} constructed for each $k\in\mathbb N\cup\{\infty\}$, an infinite  Markov shift of ergodic index $k$.   
  We note that for each infinite Markov shift, the ergodic index coincides with the conservative index.
 Infinite measure preserving rank-one transformations of an arbitrary ergodic index $k$ and infinite conservative index we constructed in \cite{AS15} and \cite{Dan2016}.
 A stronger property is
{\it power weak mixing}, which requires that for all nonzero integers $k_1,\ldots,k_r$ the product
$T^{k_1}\times\cdots\times T^{k_r}$ is ergodic \cite{DGMS}. The following examples were constructed  in \cite{afs01}, \cite{d01}, \cite{d04}, \cite{DGMS}:
\begin{itemize}
\item[(i)] power weakly mixing rank-one transformations,
\item[(ii)] non-power weakly mixing rank-one transformations with infinite ergodic index,
\item[(iii)] non-power weakly mixing rank-one transformations with infinite ergodic index and such that $T^{k_1}\times\cdots\times T^{k_r}$ are all conservative, $k_1,\ldots,k_r\in\mathbb Z$,
\end{itemize}
of types $II_\infty$ and $III$ (and various subtypes of $III$, see Section~\ref{orbits}).
Thus we have the following scale of properties (equivalent to  weak mixing in the probability preserving case), where every next property is strictly stronger than the previous ones:
$$
\begin{aligned}
\text{
$T$ is weakly mixing } &\text{$\Leftarrow$ $T$ is weakly doubly ergodic $\Leftarrow$
 $T\times T$ is ergodic} \\&\text{$\Leftarrow$  $T\times T\times T$ is ergodic }\Leftarrow\cdots\\ &\text{$\Leftarrow$ $T$ has infinite ergodic index $\Leftarrow$ $ T$ is power weakly mixing}.
\end{aligned}
$$
There is a rank-one  weakly doubly ergodic $T$ such that $T\times T$ is 
nonconservative \cite{BFMS01}
and there is a rank-one  weakly doubly ergodic $T$ such that $T\times T$ 
conservative but not ergodic \cite{LS}.
We also mention an example of a power weakly mixing transformation of type $II_\infty$ which embeds into a rank-one flow \cite{DanSol}.
This result was sharpened in \cite{DanPa}: there is an infinite measure preserving rank-one flow $(R_t)_{t\in\Bbb R}$ such that for each $t\ne 0$, the transformation $T_t$ has infinite ergodic index.
Several of these notions have been studied in the context of nonsingular actions of locally compact groups by Glasner and Weiss \cite{GlWe16}; we mention one condition that has not yet been discussed though only in the context of transformations. A nonsingular transformation $T$ on a probability space is said to be 
{\it ergodic with isometric coefficients} if every factor map
onto an isometry of 
 a (separable) metric space 
  is constant a.e. Glasner and Weiss show that if $T\times T$ is ergodic (i.e., $T$ is doubly ergodic), then $T$ is  ergodic with isometric coefficients, and that if $T$ is ergodic with isometric coefficients, then it is weakly mixing. In \cite{LS} it is shown that if $T$ is weakly doubly ergodic, then it is 
ergodic with isometric coefficients. In \cite{GlWe16} there is an example of a $T$ that is
ergodic with isometric coefficients but $T\times T$ is not conservative, hence not ergodic.
Further conditions related to weak mixing (in the case
of infinite measure-preserving transformations) are discussed in the survey \cite{AdSi18}.

\subsection{Rational ergodicity and rational  weak mixing}

Let $T$ be a conservative, ergodic  measure-preserving transformation of a $\sigma$-finite measure space 
$(X,\fL, \mu)$.  
For a function $f:X\to X$, let $S_n(f)=\sum_{k=0}^{n-1}f\circ T^k.$
$T$ is  called {\it rationally ergodic} \cite{Aa77} if there is a subset $F\in\fL$, $0<\mu(F)<\infty$, satisfying a {\it Renyi inequality},
i.e.,  there exists a constant $M>0$ such that for all $n\geq 1$,

\[\int_F(S_n(\mathbb I_F))^2\ d\mu\leq M\left(\int_F S_n(\mathbb I_F)\ d\mu\right)^2.\]
We now set
\[u_k(F):=\frac{\mu(F\cap T^{-k}F)}{\mu(F)^2}\ \text{ and }\ a_n(F):=\sum_{k=0}^{n-1}u_k(F).\]

\begin{theorem} [\cite{Aa77}]\label{rat} If $T$ is rationally ergodic  and $F$ satisfies the Renyi inequality then for all measurable sets $A$ and $B$ contained in $F$,
$$
\lim_{n\to\infty}\frac{1}{a_n(F)}\sum_{k=0}^{n-1}\mu(A\cap T^{-k}B)=\mu(A)\mu(B).
$$
\end{theorem}

An ergodic conservative transformation $T$ is {\it rationally weakly mixing} \cite{Aa13} if there exists a measurable set $F$ of positive finite measure such that for all measurable sets $A$ and $B$ contained in $F$ we have 
\begin{align}\label{raterg}
\lim_{n\to\infty}\frac{1}{a_n(F)}\sum_{k=0}^{n-1}|\mu(A\cap T^{-k}B)-\mu(A)\mu(B)u_k(F)|=0,\end{align}
where 
$u_k(F)$  and  $a_n(F)$ are defined as above.
When $\mu(X)=1$  and we let $F=X$, then $a_n(F)=n$ and \eqref{raterg} becomes the 
condition equivalent to the weak mixing property for a finite measure-preserving transformation.  In infinite measure, however, the rational weak mixing condition is not equivalent to weak mixing as we shall see.  If in \eqref{raterg} we drop the absolute values then this condition defines the notion of {\it weak rational ergodicity}  \cite{Aa80}. 
Then Theorem~\ref{rat} claims that rational ergodicity implies  weak rational ergodicity.
 If in \eqref{raterg} the sequence $(n)$ is replaced by a subsequence $(n_i)$ we say $T$
is {\it subsequence rational weak mixing}. {\it Subsequence weak rationally ergodic} is defined in a similar way.  
The transformation $T$ is {\it boundedly rationally ergodic} \cite{Aa80} if 
\[\sup_{n\geq 1}{\bigg\|}\frac{a_n(F)}{S_n(\mathbb {I}_F)}{\bigg\|}_\infty<\infty,\]
The notions of {\it subsequence boundedly rationally ergodic} and {\it subsequence rationally ergodic} are defined when the 
sequence $(n)$ is replaced by a subsequence $(n_i)$. It can be seen from the definition that bounded rational ergodicity implies 
rational ergodicity (and similarly for the subsequence versions). Aaronson \cite{Aa80} showed that rational ergodicity does not imply bounded rational ergodicity, and more recently it was shown by Adams and Silva \cite{AS16} that weak rational ergodicity does not imply rational ergodicity.  

The following theorem was proved in \cite{Aa80}  for the  weakly rationally ergodic transformations.
The same argument works in the more general case.

\begin{theorem}\label{Squa} Each subsequence weakly rationally ergodic transformation  $T$  of $(X,\mu)$ is non-squashable, i.e.,
each each non-singular transformation commuting with $T$ preserves $\mu$.
\end{theorem}

There are several examples of rationally ergodic transformations which are infinite Markov shifts, see \cite{Aa}. 
More recently, it was shown in \cite{DGPS, Bozgan}
that rank-one (infinite measure-preserving) transformations are subsequence boundedly rationally ergodic. The first version of \cite{Bozgan} has a proof that the 
rank-one transformations are subsequence weakly rationally ergodic; a simpler proof was found in \cite{D16}, where this property is also established  for the class of funny rank one transformations and the class of ergodic  transformations of {\it balanced finite rank}.
(A transformation is called of balanced finite rank if 
if it is of finite rank and the
bases of the Rokhlin towers on the $n$-th step of the cutting-and-stacking inductive construction have asymptotically comparable measures as $n\to\infty$.)
Therefore all these transformations are non-squashable in view of Theorem~\ref{Squa}.

The
rank-one transformations for which the sequence of cuts $(r_n)_{n=1}^\infty$ is bounded are boundedly rationally ergodic \cite{DGPS, Bozgan}; a 
stronger condition was shown in \cite{AKW}. 
As for the examples of rationally weakly mixing transformations,
Aaronson \cite{Aa13} shows that Markov shifts with certain conditions on their associated renewal sequences are rationally weakly mixing, and Dai et al \cite{DGPS} give rank-one examples. 
Subsequence rational weak mixing and rational weak mixing for products of powers have been studied in 
 \cite{Aa13}  and \cite{AdRigid}.

We have the following implications for rational weak mixing.

\begin{theorem} [\cite{Aa13}]  If  a transformation  is sequentially rationally weakly mixing, then it is weakly mixing.
\end{theorem}

\begin  {theorem} [\cite{Bozgan}] If a transformation  is  rationally weakly mixing, then it is weakly doubly ergodic.
\end{theorem}

It is an open problem whether weak double ergodicity implies rational weak mixing.
 Aaronson \cite{Aa13}  asked if weak rational ergodicity plus weak mixing imply
  rational weak mixing.
This was answered in negative in \cite{DGPS}, where it was constructed an example of a weakly mixing    rationally ergodic rank-one transformation that is not rationally weakly mixing.
We also mention an example of   a weakly mixing, rationally ergodic and  Koopman mixing (or zero type, see \S\ref{Mix} for the definition)  transformation that is not subsequence rationally weakly mixing \cite{A16}. 

The set of transformations that are subsequence rationally weakly mixing is residual   \cite{Aa13}, while the set
of rationally weakly mixing transformations is meagre  \cite{Aa13}. 
Since the set of power weakly mixing rank-one transformations is residual,
and the rank-one transformations are subsequence boundedly rationally ergodic, there exist rank one transformations
that are power weakly mixing and subsequence boundedly rationally ergodic but not rationally weakly mixing.

\subsection{Mixing, zero type}\label{Mix}
We now consider several attempts to define (strong) mixing for nonsingular  maps.
Probably the first notion of mixing for infinite measure preserving systems was proposed by Hopf in \cite{Ho37}. The idea was to show an asymptotic rate for the sequence $\mu(A\cap T^{-n}B)$ for a large class of finite measure sets $A, B$. 
More precisely, a transformation $T$ is {\it mixing for a ring $\mathcal R$}  (called now {\it Krickeberg mixing}), 
where  $\mathcal R$ is a ring of sets of finite measure that is invariant under $T$ and generates the entire $\sigma$-algebra, 
if there is  a sequence $(\rho_n)_{n=1}^\infty$ such that for all $A,B\in\mathcal R$ we have 
\[\lim_{n\to\infty}\rho_n\mu(A\cap T^{-n}B)=\mu(A)\mu(B).\]
Hopf proved such a property for an infinite measure-preserving transformation  defined on $\mathbb R^+\times [0,1]$ that is now called an infinite random walk; with  $\mathcal R$ being the ring of Riemann measurable subsets. 
If $\mathcal R$ is the ring of all subsets of finite measure then there are no
 Krickerberg $\mathcal R$-mixing transformations 
 because of the existence of weakly wandering sets.
 We note that the above (purely measure theoretical) definition $\mathcal R$-mixing is due to Friedman \cite{F76} who
 extended 
Krickeberg's one \cite{Kr67} given for continuous transformations of topological spaces endowed with a measure.
Recently there have been several works showing this version of mixing and computing mixing rates for several transformations. Melbourne and Terhesiu \cite{MT12} have verified mixing for a large class of maps including AFN maps with indifferent fixed points; these methods were extended to invertible transformations  by Melbourne \cite{M15} and to additional maps by Gou\"ezel \cite{Go11}. Recently, Dolgopyat and N\'{a}ndori \cite{DN19} have shown Krickerberg mixing for a class of special flows; other recent work appeared in \cite{BMT19}.

Another approach to mixing was proposed by Krengel and Sucheston \cite{KS69}
for nonsingular maps.
Given a sequence of measurable sets $\{A_n\}$ let $\sigma_k(\{A_n\})$ denote the $\sigma$-algebra generated by $A_k, A_{k+1},\ldots$. A sequence $\{A_n\}$  is said to be {\it remotely trivial}
if $\bigcap_{k=0}^\infty \sigma_k(\{A_n\}) = \{\emptyset, X\}$ mod $\mu$, and it is
{\it semi-remotely trivial} if every subsequence contains a further subsequence that is remotely trivial.
A nonsingular transformation $T$
of a $\sigma$-finite measure space is called {\it mixing} if for every set $A$ of finite measure the sequence
$\{T^{-n}A\}$ is semi-remotely trivial, and {\it completely mixing} if $\{T^{-n}A\}$ is semi-remotely trivial
for all measurable sets $A$.  
Krengel and Sucheston show that $T$ is
completely mixing if and only if it is type $II_1$ and mixing for the equivalent finite invariant measure.
Thus there are no type $III$ and $II_\infty$ completely mixing  nonsingular transformations on probability spaces. We note that this definition of mixing in infinite measure spaces depends on the choice of measure inside the equivalence class (but it is independent if we replace the measure by an equivalent measure with the same collection of sets of finite measure).

Hajian and Kakutani showed \cite{HK} that
an ergodic infinite measure-preserving  transformation $T$ is either of {\it zero type}:
$\lim_{n\to\infty}\mu(T^{-n}A\cap A)=0$ for all sets $A$ of finite measure, or of
{\it positive type}: $\limsup_{n\to\infty}\mu(T^{-n}A\cap A)>0$ for all subsets $A$ of finite  positive measure.  
It appears that $T$ is mixing if and only if it is of zero type \cite{KS69}.
We note that in infinite measure, mixing implies mixing of all orders, i.e., if a measure preserving $T$ is of zero type
then $\mu( T^{n_1}A_1\cap\cdots\cap T^{n_k}A_k)\to 0$ for  each $k$ and all subsets $A_1,\dots,A_k$ of finite measure whenever $|n_i-n_j|\to\infty$ if $i\ne j$.

 For $0\le\alpha\le 1$ Kakutani suggested a related definition of $\alpha$-type: an infinite measure preserving transformation is {\it of $\alpha$-type} if $\lim\sup_{n\to\infty} \mu(A\cap T^nA)=\alpha\mu(A)$ for every subset $A$ of finite measure. In \cite{OH}
examples of ergodic transformations of any $\alpha$-type and a transformation of not any type were constructed.

It was shown in \cite{Dan2016} and \cite{loh} that for each  pair $k\le n$, there exists a mixing rank-one infinite measure preserving transformation of ergodic index $k$
and conservative index $n$.
Rigid infinite measure preserving rank-one transformations of arbitrary ergodic index were constructed in \cite{Dan2016}.
Of course,  rigidity implies infinite conservative index.

We now isolate an important class of concrete rank-one transformations and examine mixing properties within this class.
Let $T$ be a rank-one transformation associated with a sequence $(r_n,w_n,s_n)_{n=1}^\infty$.
If $w_n(0)=w_n(1)=\cdots= w_n(r_n-1)$ and  $s_n(j)=z_n+j$ for $j=0,\dots,r_{n-1}$
then $T$ is called {\it a high staircase} (called also tower staircase in \cite{BFMS01}).
It was shown in \cite{BFMS01} that each high staircase is weakly doubly ergodic and hence weak mixing.
However there exist high staircases whose Cartesian square is not ergodic \cite{BFMS01}.
As for the mixing of the high staircases, the following theorem was proved in \cite{DanR2}.
It is an infinite analogue of the Adams solution \cite{Ada} of the Smorodinsky conjecture.

\begin{theorem} If $\lim_{n\to\infty}\frac{r_n^2}{r_1\cdots r_{n-1}}=0$ and $\sum_{n=1}^\infty\frac{z_n}{h_n}=\infty$ then the associated high staircase is infinite measure preserving and mixing.
\end{theorem}

Mixing  high staircases which are power weakly mixing
were constructed in \cite{DanR2}.

We note that mixing (zero type) does not imply either ergodicity or conservativity in the category of infinite measure-preserving transformations.
Indeed, a translation on $\Bbb R$ endowed with the Lebesgue measure
is non-ergodic, totally dissipative but of zero type.
It may seem that mixing plus ergodicity together are stronger than any kind of nonsingular weak mixing considered above. 
However, it is not the case: if $T$ is a weakly mixing infinite measure-preserving transformation of zero type and $S$ is an ergodic probability preserving transformation then $T\times S$ is ergodic and of zero type. On the other hand, the $L^\infty$-spectrum $e(T\times S)$ is nontrivial, i.e., $T\times S$ is not weakly mixing, whenever $S$ is not weakly mixing. We also note that
there exist rank-one infinite measure-preserving transformations $T$ of zero type such that $T\times T$ is not conservative (hence not ergodic) \cite{AFS}. 
In contrast to that,
if $T$ is of positive type then  all of its finite Cartesian products are conservative \cite{AN00}. Another result that suggests that there is no good definition of mixing in the  infinite measure-preserving case was proved  in
\cite{james}. It is shown there that while the mixing finite measure-preserving transformations are
 measurably sensitive, there exists no  infinite measure-preserving system that
is  measurably sensitive. (Measurable sensitivity is a measurable version of the strong sensitive dependence on initial conditions---a concept from topological theory of chaos.)

The Krengel-Sucheston concept of  mixing (or the Hajian-Kakutani zero type) considered above for infinite measure-preserving systems extends  naturally to nonsingular transformations $T$
 of $(X,\mathcal B,\mu)$ without finite absolutely continuous  invariant measure in the following way: we say that
 $T$ is {\it Koopman mixing} (or {\it of zero type}) if the maximal spectral type of the Koopman operator $U_T$ generated by $T$ (see  (\ref{koopman oper})) is a Rajchman measure, i.e., 
$\int_Xf\cdot U_T^nf\,d\mu\to 0$ for each $f\in L^2(X,\mu)$.
It is easy to see that this definition of mixing  will not be affected if we replace $\mu$ with an equivalent measure. 
Examples of Koopman mixing rank-one transformations of type $III$  were constructed in \cite{d01}.

More recently, Lenci \cite{Le1} introduced a new notion of mixing for infinite measure-preserving maps 
that is motivated by statistical mechanics and uses global observables. The definition is with respect 
to  a collection of sets, global observables and local observables. We choose a family $\mathcal V$ of measurable sets of finite measure so that it contains sets $V_1\subset V_2\subset \cdots $ such that 
$\bigcup_i V_i =X$. We also have a subspace $\mathcal G$ of $L^\infty$ functions (called global observables),
and a subspace $\mathcal L$ of $L^1$ functions (called local observables). There is also a condition on the growth rate of the measure of  $\mathcal V$-elements under iteration by $T$. Then Lenci defines an {\it infinite volume average} for elements $F$ of $\mathcal G$ by
\[\bar{\mu}(F)=\lim_{V\to X}\int_V F \ d\mu.\]
By this limit we mean that for every neighborhood of  $\bar{\mu}(F)$, there is a number $M>0$ so that 
when $\mu(V)>M$ for a set $V$ in $\mathcal V$, then  $\int_VFd\mu$ is in the neighborhood. 
He shows that under the above conditions,
 $\bar{\mu}(F\circ T^n)=\bar{\mu}(F)$. 
Then he defines several notions of what he calls infinite volume mixing  \cite{Le1}; we mention three here. 
The transformation $T$ is said to be {\it global-local mixing-1} if for all $F$ in $\mathcal G$ and all $g$ in $\mathcal L$  with   $\int g\ d\mu=0$, we have
\[\lim_{n\to\infty} \int (F\circ T^n)g\ d\mu= 0.\]  The transformation $T$ is said to be {\it global-local mixing-2} if for all $F$ in $\mathcal G$ and all $g$ in $\mathcal L$ we have
\[\lim_{n\to\infty} \int (F\circ T^n)g\ d\mu= \bar{\mu}(F) \int g\ d\mu.\]
The transformation is said to be 
{\it global-global mixing} if 
for all $F, G$ in $\mathcal G$   we have
\[\lim_{n\to\infty} \bar{\mu}(F\circ T^nG)\ d\mu= \bar{\mu}(F) \bar{\mu}(G)\]
Lenci proves in \cite{Le3} that if $T$ is an infinite measure-preserving K-automorphism then 
$T$  is global-local mixing-1 for any choice
of $\mathcal V$ satisfying the measure growth condition, for $\mathcal L=L^1$, and for $\mathcal C$ 
that is  the closure in $L^1$ of $T^n\mathcal F $, where $\mathcal F$ is as in the definition 
of the K-automorphism (see Subsection~\ref{Kprop}). 
 Infinite mixing has been shown for other examples \cite{Le1}, in particular 
for uniformly expanding maps of the interval \cite{Le17}, and for one-dimensional maps with 
an indifferent fixed point \cite{BGL18}.

\subsection{$K$-property}\label{Kprop}
A nonsingular transformation $T$ of $(X,\mathcal B,\mu)$ is
called {\it K-automorp\-hism} \cite{ST2} if there exists a
sub-$\sigma$-algebra $\mathcal F\subset\mathcal B$ such that $T^{-1}\mathcal
F\subset\mathcal F$, $\bigcap_{k\ge 0}T^{-k}\mathcal F=\{\emptyset,X\}$,
$\bigvee_{k=0}^{+\infty}T^k\mathcal F=\mathcal B$ and the Radon-Nikodym
derivative $\omega_1^\mu$ is $\mathcal F$-measu\-rable (see also
\cite{wP65} for the case when $T$ is of type $II_\infty$; the authors in \cite{ST2} required $T$ to be conservative).
If $R$ is a nonsingular endomorphism on $(X,\mathcal B,\mu)$ then the natural extension of $R$ is a $K$-automorphism if and only if $R$ is exact, i.e., $\bigwedge_{n=1}^\infty R^{-n}\mathcal B=\{\emptyset, X\}$ mod 0.
It follows from the Kolmogorov 0-1 theorem that a  nonsingular Bernoulli shift from the generalized Krengel class (see \S\,\ref{S:hamachi}) is a
 $K$-automorphism.
Parry \cite{wP65} showed that a type $II_\infty$ $K$-automorphism is either dissipative or ergodic.
Krengel \cite{Kre70}  proved the same for the Krengel class of Bernoulli nonsingular shifts, and finally Silva and Thieullen extended  this result to  the nonsingular $K$-automorphisms \cite{ST2}.
 It is
  also shown in \cite{ST2} that if $T$ is a
nonsingular $K$-automorphism, for any ergodic  nonsingular transformation $S$, if $S\times T$ is
conservative, then it is ergodic. This property is called {\it sharply weak mixing} in \cite{DanL22}.
It follows  that a conservative nonsingular $K$-automorphism is weakly mixing. However, it does not necessarily have infinite ergodic index \cite{KP}.
 Krengel and Sucheston \cite{KS69} showed that an infinite measure-preserving conservative
$K$-automorphism is mixing. 
In \cite{EKRSS}  it is shown that if $T$ is a rank-one transformation with $(r_n)$ bounded and 
$s_{n-1}(r_{n-1}-1)\geq h_n/2$ for all $n\geq 1$  (so the space is of infinite measure, and an $T$ is the infinite Chac\'on map), there exists a rank-one transformation $S$ such that $T\times S$ is conservative but not ergodic, so it is not sharply weak mixing.

\subsection{Multiple and polynomial recurrence}\label{MPR}
Let $p$ be a positive integer.
A nonsingular transformation $T$
is called $p$-{\it recurrent} if for every subset
$B$ of positive measure there exists a positive integer $k$
such that
$$
\mu(B\cap T^{-k}B\cap \dots \cap T^{-kp}B)>0.
$$
 If $T$ is $p$-recurrent for any $p>0$, then it is called
{\it multiply recurrent}.
It is easy to see that $T$ is 1-recurrent if and only if it is conservative.  
 Clearly, if $T$ is  rigid then it is multiply recurrent.  
 Furnstenberg showed \cite{Fu} that  every finite measure-preserving
transformation is multiply recurrent. 
In contrast to that,
Eigen, Hajian and Halverson \cite{ehh} constructed for any $p\in\mathbb
N\cup\{\infty\}$, a nonsingular product odometer of type $II_\infty$ which is $p$-recurrent but  not $(p+1)$-recurrent.
Aaronson and Nakada showed in \cite{AN00} that an infinite measure
preserving Markov shift $T$ is $p$-recurrent if and only if the product
$T\times \cdots\times T$ ($p$ times) is conservative. It
follows from this and \cite{ALW} that in the class of ergodic Markov
shifts, infinite ergodic index implies multiple recurrence. However, in
general this is not true. It was shown in \cite{afs01}, \cite{Ketal03} and \cite{DanS} that for each $p\in\mathbb N\cup\{\infty\}$ there exist
\begin{itemize}
\item[(i)]
 power weakly mixing rank-one transformations and
\item[(ii)] non-power weakly mixing rank-one transformations  with infinite ergodic index
\end{itemize}
which are $p$-recurrent but not $(p+1)$-recurrent (the latter holds when $p\ne\infty$, of course).

A subset $A$ is called $p$-{\it wandering} if $\mu(A\cap T^kA\cap\cdots\cap T^{pk}A)=0$ for each $k$. Aaronson and Nakada established in \cite{AN00} a $p$-analogue of Hopf decomposition (see Theorem~\ref{T:hopf}).

\begin{prop} If $(X,\sS,\mu,T)$ is conservative aperiodic nonsingular dynamical system and $p\in\mathbb N$ then $X=C_p\sqcup D_p$, where $C_p$ and $D_p$ are $T$-invariant disjoint subsets, $D_p$ is a countable union of $p$-wandering sets, $T\restriction C_p$ is $p$-recurrent and $\sum_{k=1}^\infty\mu(B\cap T^{-k}B\cap\cdots\cap T^{-dk}B)=\infty$ for every $B\subset C_p$.
\end{prop}

Let $T$ be an infinite measure-preserving transformation and  let $\mathcal F$
be a $\sigma$-finite factor (i.e., invariant subalgebra) of $T$. Inoue \cite{Ino}
showed  that for each $p>0$, if $T\restriction \mathcal F$
is $p$-recurrent then so is $T$ provided that the extension
$T\to T \restriction \mathcal F$ is isometric.  It is unknown yet whether
the latter assumption can be dropped.  However, partial progress was
 achieved in \cite{Mey}: if  $T\restriction\mathcal F$ is multiply
recurrent then so is $T$.

Let $\mathcal P:=\{q\in\mathbb Q[t]\mid q(\mathbb Z)\subset\mathbb Z \text{ and }
q(0)=0\}$. An ergodic  conservative nonsingular transformation $T$  is called
 $p$-{\it polynomially recurrent} if for every
$q_1,\dots,q_p\in\mathcal
P$ and every subset $B$ of positive measure there exists $k\in\mathbb N$ with
$$
\mu(B\cap T^{q_1(k)}B\cap\dots\cap T^{q_p(k)}B)>0.
$$
If $T$ is $p$-polynomially recurrent
 for every $p\in\mathbb N$
then it is called {\it polynomially recurrent}. Furstenberg's theorem on multiple recurrence was significantly strengthened in \cite{BeL}, where it was shown that every finite measure-preserving transformation is polynomially recurrent.
However, Danilenko and Silva \cite{DanS} constructed
\begin{itemize}
\item[(i)]  type $II_\infty$ transformations $T$ which are $p$-polynomially recurrent but not $(p+1)$-polynomially recurrent (for each fixed $p\in\mathbb N$),
\item[(ii)]  polynomially recurrent   transformations $T$ of type $II_\infty$,
\item[(iii)] rigid (and hence multiply recurrent) type $II_\infty$ transformations $T$ which are not polynomially recurrent.
\end{itemize}
Moreover, such $T$ can be chosen inside the class of rank-one transformations with infinite ergodic index.

\section{Dynamical properties of IDPFT systems}
Let  $T_n$ be an ergodic  transformation of  a standard probability space $(X_n,\mu_n, T_n)$ and let $\nu_n$ be a $\mu_n$-equivalent invariant probability measure
for each $n\in\Bbb N$.
Assume that the transformation $T=\bigotimes_{n=1}^\infty T_n$ is nonsingular 
with respect to the infinite product measure
$\mu:=\bigotimes_{n=1}^\infty\mu_n$.
In other words, 
$(X,\mu, T)$ is an IDPFT dynamical system.

\begin{theorem}[\cite{DanL18}] The  transformation $T$ is either conservative or totally dissipative.
Moreover,
if $S$ is  an ergodic conservative nonsingular transformation then the direct product 
$T\times S$ is either conservative or totally dissipative.
\end{theorem}

\begin{theorem}[\cite{DanL18}]
Let $(X_n,\nu_n,T_n)$ 
be mildly mixing (see \S\ref{mild mixing} for the definition) for each $n>0$.
If $T$ is $\mu$-conservative then $T$ is sharply weak mixing.
\end{theorem}

Examples of rigid ergodic but not weakly mixing IDPFT transformations of Krieger type $III_\lambda$, for each $\lambda\in(0,1)$ were constructed in \cite{DanL18}.
Some families of 0-type IDPFT transformations of type $III_1$ appeared in \cite{DanL18}
and of all possible Krieger types
 in \cite{DaKo}.

\begin{theorem}[\cite{DaKo}]
 Let $K\in\{III_\lambda\mid 0\le \lambda\le 1\}\sqcup\{II_\infty\}$.
Then there is a 0-type  weakly mixing IDPFT transformation
 of type $K$.
\end{theorem}

\section{Dynamical properties of nonsingular Bernoulli and Markov shifts}
We will use below the notation introduced in \S\ref{S:hamachi}.
Thus,
 $T$ stands for the  nonsingular Bernoulli shift on the space $(X,\mu)=(A^\Bbb Z,\bigotimes_{n\in\Bbb Z}\mu_n)$ associated with a sequence of   probability measures $(\mu_n)_{n\in\Bbb Z}$ on $A$.
\subsection{Krengel class}
 Nonsingular Bernoulli shifts appeared  originally in Krengel's work 
 \cite{Kre70}.
  He introduced there a class of nonsingular shifts for which $A=\{0,1\}$ and $\mu_n$ is the equidistribution on $\{0,1\}$ for all $n\le 1$.
We  will call it {\it the Krengel class}.
Krengel showed that  this class contains totally dissipative transformations.
He also used an inductive procedure to construct the sequence $(\mu_n)_{n=1}^\infty$ in such a way that the corresponding  Bernoulli shift  is ergodic conservative and not of type $II_1$.
Krengel conjectured  that the shift is of type $III$ indeed.
In \cite{Ham3}, Hamachi showed that    Krengel's class contains ergodic conservative nonsingular Bernoulli shifts of type $III$.
This was further refined by Kosloff who constructed type $III_1$ ergodic conservative shifts belonging to Krengel's class \cite{Ko11}.
In \cite{Ko13}
Kosloff  constructed a nonsingular Bernoulli shift of type $III_1$ (and belonging to the Krengel class)  which is power weakly mixing.
Weiss asked about possible Krieger's types for the nonsingular Bernoulli shifts.
Answering his question, Kosloff proved in 
 a subsequent paper \cite{Ko14}  that each conservative Bernoulli shift from the Krengel class is ergodic and either of type $II_1$  or of type $III_1$.
In particular, the non-type-$II_1$ conservative Bernoulli shifts constructed in \cite{Kre70}, \cite{Ham3} and \cite{Ko11} are all of type $III_1$ indeed.

\subsection{Generalized Krengel class}
Kosloff's result from \cite{Ko14} was further extended in \cite{DanL18}.
We  say that a nonsingular Bernoulli shift belongs to {\it the generalized Krengel class} if $A=\{0,1\}$ and $\mu_n=\mu_1$ for each $n\le 0$.
We  note that these transformations
are the natural extension of the one-sided nonsingular Bernoulli shifts
 defined on $(A^\Bbb N, \bigotimes_{n>0}\mu_n)$.
Every shift from the generalized Krengel class  is a $K$-automorphism.

\begin{theorem}[On types of nonsingular Bernoulli shifts from the generalized Krengel class  \cite{Ko14}, \cite{DanL18}]\label{BernShiftType}
 Let $A=\{0,1\}$ and let $T$ be a nonsingular Bernoulli shift on $(A^{\Bbb Z},\bigotimes_{n\in\Bbb Z}\mu_n)$ from the generalized Krengel class.
\begin{enumerate}
\item[(i)] If $\sum_{n>0}(\mu_n(0)-\mu_1(0))^2<\infty$ then $\mu$ is equivalent to 
$\bigotimes_{n\in\Bbb Z}\mu_1$ and hence $T$ is of type $II_1$.
\item[(ii)] If $\sum_{n>0}(\mu_n(0)-\mu_1(0))^2=\infty$ and $T$ is conservative then $T$ is ergodic of type $III_1$.
Moreover, the Maharam extension of $ T$ is a  weakly mixing $K$-automorphism.
\end{enumerate}
\end{theorem}

Thus, Krieger's type of each nonsingular Bernoulli shift from the generalized Krengel class is never of type $III_\lambda$, $0\le\lambda<1$.
In  \cite{VaWa} 
Vaes and Wahl, answering a question from \cite{DanL18},  found a convenient condition 
for a nonsingular Bernoulli shift from the generalized  Krengel class to be conservative.
Utilizing that condition they
constructed, for each $\lambda\in (0,1)$, an explicit example of 
power weakly mixing nonsingular Bernoulli shift of type $III_1$  with $\mu_1(0)=\lambda$.
We note  that the previously known  Bernoulli shifts of type $III_1$
were constructed via involved inductive procedures. 
Vaes and Wahl also provided an example of type $III_1$ Bernoulli shift  with finite ergodic index (less than 73) in \cite{VaWa}.
Modifying their example Kosloff and Soo showed that there are nonsingular Bernoulli shifts
of type $III_1$ of arbitrary ergodic index.

\begin{theorem}[\cite{KoSo2}]
Let $A=\{0,1\}$ and $c >0$.
 Let
$
\mu_n^c(0):=\frac 12+\frac c{\sqrt n}1_{\{n\in\Bbb N\mid \sqrt n>2c\}}
$
for each $n\in\Bbb Z$.
There exists $D >\frac16$ such that the Bernoulli shift $T$ on
$\bigotimes_{n\in\Bbb Z}(A,\mu_n)$
is ergodic of type $III_1$ for all $c < D$ and totally dissipative for all $c > D$. 
In addition, if $k\in \Bbb N$ and $\frac D{\sqrt{k+1}}\le c< \frac D{\sqrt k}$ then $T$
is of ergodic index $k$.
\end{theorem}

\subsection{General nonsingular Bernoulli shifts}
The studying of general nonsingular Bernoulli shifts  was initiated by Kosloff
in \cite{Ko13}.

\begin{theorem}[Mixing of nonsingular Bernoulli shifts \cite{Ko13}]
If $\#A=2$ and  $(\mu_n)_{n\in\Bbb Z}$ is a sequence of probabilities on $A$ such that
{\rm (\ref{ffff})} holds.
Then
$T$ is
either of type $II_1$ and mixing (with respect to the equivalent invariant probability measure) or of zero type.
\end{theorem}

In  \cite{Ko18}, Kosloff noticed that under some natural conditions, conservativity of Bernoulli shifts implies ergodicity. 
His proof was based essentially on the Hurewicz ergodic theorem
and properties of the tail equivalence relation on $A^{\Bbb Z}$.
Danilenko \cite{Dan2018} refined his results  by exploiting the interplay  between $T$
and the measurable equivalence relation on $A^{\Bbb Z}$ generated by the finite permutations of coordinates.

\begin{theorem}[Weak mixing of conservative nonsingular Bernoulli shifts]\label{BernWeak} Let $A$ be finite. 
\begin{enumerate}
\item[(i)]
If $\inf_{n\in\Bbb Z}\min_{a\in A}\mu_n(a)>0$ and $T$ is conservative then
$T$ is weakly mixing $($see \cite{Ko18}, \cite{Dan2018}$)$.
\item[(ii)] If $\# A=2$ and $\inf_{n\in\Bbb Z}\min_{a\in A}\log|\frac{\mu_n(a)}{\mu_{n+1}(a)}|>0$ and $T$ is conservative then
$T$ is weakly mixing \cite{Dan2018}.
\item[(iii)]
Under condition
{\rm (i)} or 
{\rm (ii)}, if
$T\times\cdots\times T$ ($p$ times)
is conservative for some
$p \ge 1$
then
$T\times\cdots\times T$ ($p$ times)
is weakly mixing \cite{Dan2018}.
\end{enumerate}
\end{theorem}

The techniques from \cite{Ko18} and \cite{Dan2018} were further elaborated in \cite{BjKoVa}
to prove the following refinement of Theorem~\ref{BernWeak}.

\begin{theorem}\label{th7.4}
Let $A=\{0,1\}$ and let $T$ be a nonsingular Bernoulli shift on $(A^\Bbb Z,\bigotimes_{n\in\Bbb Z}\mu_n)$ and let $\mu=\bigotimes_{n\in\Bbb Z}\mu_n$ be nonatomic.
If $T$ is not totally dissipative then $T$ is weakly mixing and its type is given as follows.
\begin{enumerate}
\item[(i)]
If there is $\lambda\in(0,1)$ such that $\sum_{n\in\Bbb Z}(\mu_n(0)-\lambda)^2<+\infty$
then  $T$ is of type $II_1$.
\item[(ii)]
If there is $\lim_{n\to\infty}\mu_n(0)=\lambda\in(0,1)$ and $\sum_{n\in\Bbb Z} (\mu_n(0)-\lambda)^2=+\infty$ then $T$ is of type $III_1$.
\item[(iii)]
 If either $\lim_{n\to\infty}\mu_n(0)=0$ or $\lim_{n\to\infty}\mu_n(0)=1$ then $T$ is of type $III$.
 \item[(iv)]
If the sequence $(\mu_n(0))_{n}$ does not converge as $n\to\infty$
then $T$ is of type $III_1$.
\end{enumerate}
\end{theorem}

It follows, in particular, that the nonsingular Bernoulli shifts on $\{0,1\}^{\Bbb Z}$ are never of type  $II_\infty$.
It is still an open problem which subtypes of $III$ are realized in the case (iii). 
An example of  $T$ of  type $III_1$ with  $\lim_{n\to\infty}\mu_n(0)=0$ was constructed in \cite{BeVa}.
Is there $T$ of type $III_0$ satisfying~(iii)?

However, the situation is different if we consider $A=[0,1]$
and the probability measures $(\mu_n)_{n\in\Bbb Z}$ on $A$ with infinite support.
In \cite{KoSo}, for each $\lambda\in(0,1)$, an example of type $III_\lambda$ Bernoulli shift was constructed with $\mu_n\sim \text{Leb}$ for each $n$.
Examples of nonsingular Bernoulli shifts of each possible Krieger's type were given in a later paper \cite{BeVa}.
An alternative proof of this result  appeared in a recent work \cite{DaKo}.

\begin{theorem}[\cite{DaKo}]\label{I+Z}
 Let $K\in\{III_\lambda\mid 0\le \lambda\le 1\}\sqcup\{II_\infty\}$.
Then there is a sequence $(\mu_n)_{n\in\Bbb Z}$ of probability measures on $[0,1]$ such that $\mu_n\sim\text{{\rm Leb}}$ for each $n\in\Bbb N$, the Bernoulli shift $T$ on $([0,1]^\Bbb Z,\bigotimes_{n\in\Bbb Z}\mu_n)$
is weakly mixing, IDPFT, and of Krieger type $K$.
\end{theorem}

We also note that if  there is $C>1$ such that $C^{-1}\le \frac{d\mu_n}{d\mu_0}\le C$ for all $n\in\Bbb Z$ then
$T$ can  be neither of type $III_0$ nor of type $II_\infty$  \cite{BeVa}.

The following result is an analog of Theorem~\ref{th7.4} for Bernoulli shifts with 
general state spaces.

\begin{theorem}[\cite{BeVa2}]
Let $(\mu_n)_{n\in\Bbb Z}$ be a sequence of mutually equivalent probabilities on a 
standard Borel space $A$ and let $T$ be the  Bernoulli shift on
$(X,\mu)=\bigotimes_{n\in\Bbb Z}(A,\mu_n)$.
If $T$ is nonsingular and conservative then $T$ is weakly mixing and the following holds.
\begin{enumerate}
\item[(i)]
$T$ of type $II_1$ if and only if there exists a probability measure
$\nu\sim\mu_0$ such that
$\nu^{\otimes Z}\sim\mu$.
\item[(ii)]
$T$ of type $II_\infty$ if and only if 
if and only if there exists a $\sigma$-finite measure
$\nu\sim\mu_0$
and Borel subsets  $B_n\subset A$
such that
$\nu(B_n)<\infty$ for each $n\in\Bbb Z$ and
$$
\sum_{n\in\Bbb Z}\mu_n(A\setminus B_n)<\infty, \ 
\sum_{n\in\Bbb Z}\nu_n(A\setminus B_n)=\infty, \
\sum_{n\in\Bbb Z} H^2(\mu_n, \nu(B_n)^{-1}\nu\restriction B_n)
$$
\item[(iii)]
$T$ is of type $III$ in all other cases.
\end{enumerate}
\end{theorem}

The next very natural question to ask is: which ergodic flows arise as the associated flow of nonsingular Bernoulli shifts of type $III_0$?

\begin{defn} Let $V=(V_t)_{t\in\Bbb R}$ be a nonsingular flow on a standard probability space $(X,\nu)$.
Given $n>1$, consider two mutually commuting nonsingular flows $U=(U_t)_{t\in\Bbb R}$ and $D=(D_t)_{t\in\Bbb R}$ on the product space
$(X,\nu)^{\otimes n}$:
$$
U_t(x_1,\dots,x_n):=(V_tx_1,x_2,\dots, x_n), \ D_t(x_1,\dots,x_n):=(V_tx_1,V_tx_2,\dots, V_tx_n),
$$
for  each $t\in\Bbb R$ and $(x_1,\dots, x_n)\in X^n$.
The restriction of $U$ to the $\sigma$-algebra of $D$-invariant subsets in $X^n$
is a well defined nonsingular flow.
It is called the {\it joint flow of $n$ copies of $V$}.
The flow $V$ is called {\it infinitely divisible} if for each $n>1$,  there is a flow $W$ such
that $V$ is isomorphic to the joint flow of $n$ copies of $W$.
\end{defn}

If $V$ is an AT-flow associated with an ITPFI$_2$-transformation  (see Section~\ref{AT-f}) then $V$
is infinitely divisible.

\begin{theorem}[\cite{BeVa2}]\label{recent}
\begin{enumerate}
\item[(i)]
Let $(\mu_n)_{n\in\Bbb Z}$ be a sequence of mutually equivalent probabilities on a 
standard Borel space $A$ and let $T$ be the  Bernoulli shift on
$\bigotimes_{n\in\Bbb Z}(A,\mu_n)$.
If $T$ is nonsingular and conservative then the associated flow of $T$
is infinitely divisible.
For each ergodic probability preserving transformation $S$, the product $T\times S$ is ergodic and the associated flows of $T$ and $T\times S$ are isomorphic. 
\item[(ii)]
Conversely, let $V$ be an AT-flow associated with an 
ITPFI$_2$-transformation.
Then there is a weakly mixing nonsingular Bernoulli shift whose associated flow is isomorphic to $V$.
\end{enumerate}
\end{theorem}

We note that formally Berendschot and Vaes proved in \cite{BeVa2} a stronger claim than (ii).
They introduced a concept of {\it Poisson flow} as the tail boundary of certain nonhomogeneous random walks on $\Bbb R$ (see \S\ref{brw} below).
Every Poisson flow is infinitely divisible.
Every AT-flow associated with an ITPFI$_2$-transformation is Poisson.
It is shown in \cite{BeVa2} that Theorem~\ref{recent}(ii) holds for each Poisson flow $V$.
However, the following problems remain open:
\begin{itemize}
\item
Whether each Poisson flow is an AT-flow associated 
with an ITPFI$_2$-transformation?
\item
Whether each infinitely divisible flow is Poisson?
Whether each infinitely divisible AT-flow is Poisson?
\end{itemize}

\subsection{Bernoulli factors  of nonsingular Bernoulli shifts}
\begin{defn}Let $T$ be a nonsingular Bernoulli shift on  $(X,\mu):=\bigotimes_{n\in\Bbb Z}(A,\mu_n)$.
Suppose that there is another  sequence $(\nu_n)_{n\in\Bbb Z}$ of probabilities on $A$
such that $T$ is nonsingular on the product space $(X,\nu):=\bigotimes_{n\in\Bbb Z}(A,\nu_n)$.
If there is a map $\pi:A^\Bbb Z\to A^{\Bbb Z}$
such that $\pi T=T\pi$ and the measure $\mu\circ \pi^{-1}$ is equivalent to $\nu$ then $(X,\nu,T)$ is called {\it a factor} of $(X,\mu,T)$.
\end{defn}
We  assume that  the factors  are non-trivial.
What is the relation between Krieger's types of $(X,\nu,T)$ and $(X,\mu,T)$?

\begin{theorem}(Bernoulli factors of type $II_1$ \cite{KoSo2})
Let $A$ be finite and  
$$
1>\inf_{n\in\Bbb Z}\min_{a\in A}\mu_n(a)>0.
$$
 Then there is a probability measure $\rho$ on $A$ 
 such that $(X,\rho^{\otimes\Bbb Z}, T)$ is a factor of $(X,\mu,T)$.
\end{theorem}

\begin{theorem}[ \cite{KoSo2}]
Let $\lambda,\lambda'\in(0,1]$.
 There exists a type $III_\lambda$ Bernoulli shift
  which has a type $III_{\lambda'}$ Bernoulli shift as a factor
   in each of the two cases:
   \begin{enumerate}
    \item[(i)]
  $\lambda'<\lambda=1$,
  \item[(ii)]
$\lambda<\lambda'$.
\end{enumerate}
\end{theorem}

\subsection{Nonsingular Markov shifts}
An analog  of Theorem~\ref{BernWeak}(i) holds also for nonsingular Markov shifts (see \cite{Ko18}, \cite{Dan2018}).

\begin{theorem}[Weak mixing of conservative nonsingular Markov shifts \cite{Ko18}, \cite{Dan2018}]
Let $A$ be finite and let
 $M=(M(a,b))_{a,b\in A}$ be a 0-1-valued  $A\times A$-matrix.
 Suppose that $M$ is primitive, i.e., there is $n>0$ such that all the entries of $M^n$ are strictly positive.
 Let $(X_M,T,\mu)$ be a nonsingular Markov shift and let $\mu$ be generated by a sequence $(\pi_n,P_n)_{n\in\Bbb Z}$ as in {\rm \S\ref{S:MarkSh}.}
 Suppose that $\mu$ is nonatomic and that $\pi_n$ is fully supported on $A$ for each $n$.
If $\inf\{P_n(a,b)\mid n\in\Bbb Z,  M(a,b)=1\}>0$ and $T$ is conservative then
$T$ is weakly mixing. 
\end{theorem}

 We  isolate a class of nonsingular Markov shifts for which $P_n=P_1$ and $\pi_n=\pi_1$ for all $n\le 0$ and call it {\it the Markov-Krengel} class.
 Each shift from this class is the natural extension of the corresponding one-sided nonsingular Markov shift \cite{DanL18}.
There is  an analog of  Theorem~\ref{BernShiftType} for the Markov-Krengel shifts.

\begin{theorem}[\cite{DanL18}]
 Let $M^*=\begin{pmatrix}
1 & 1\\
1& 1
\end{pmatrix}$, 
$P_n$ be a bistochastic matrix for each $n\in\Bbb Z$, $P_k=\begin{pmatrix}
0.5 & 0.5\\
0.5& 0.5
\end{pmatrix}$ and $\pi_k=(0.5, 0.5)$ for each $k\le 0$.
Let the corresponding Markov-Krengel shift $(X_{M^*},T,\mu)$ be  nonsingular  and conservative.
Then either $T$ is of type $II_1$ (if $\sum_{n>0}|P_n(0,0)-0.5|<\infty$) or $III_1$ (otherwise).
In the latter case the Maharam extension of $ T$ is a  weakly mixing $K$-automorphism.
Moreover, if $\mu$ is equivalent to a Bernoulli (i.e., infinite product) measure then $T$ is of type $II_1$.
\end{theorem}

Concrete examples of Markov-Krengel shifts $(X_{M},\mu,T)$ of type $III_1$ 
such that $\mu$ is not equivalent to a Bernoulli measure were constructed in \cite{DanL18}
and \cite{Ko15}.

Recently Avraham-Re'em extended and refined the aforementioned results on Markov shifts
\cite{AR}.
To state his results, we  introduce some notation.
Given $n\in\Bbb Z$ and $a,b\in A$ we let
$$
\widehat P_n(a,b)=
\begin{cases}\frac{\pi_{n-1}(b)}{\pi_n(a)}P_{n-1}(b,a)&\text{if }\pi_n(a)\ne 0,\\
0&\text{otherwise.}
\end{cases}
$$
For a stochastic $A\times A$-matrix $Q$, let $\lambda$ be the distribution on $A$ such that
$\lambda Q=\lambda$.
We let $\widehat Q(a,b):=\frac{\lambda(b)}{\lambda(a)}Q(b,a)$.

\begin{theorem}[\cite{AR}]\label{ART}
Let $A$ be finite and let
 $M=(M(a,b))_{a,b\in A}$ be a primitive 0-1-valued  $A\times A$-matrix.
  Let a measure $\mu$ on $X_M$ be generated by a sequence $(\pi_n,P_n)_{n\in\Bbb Z}$ as in {\rm \S\ref{S:MarkSh}} and 
    $\inf\{P_n(a,b)\mid n\in\Bbb Z,  M(a,b)=1\}>0$.
If 
the Markov shift 
 $(X_M,T,\mu)$ is nonsingular and conservative and 
 
 \begin{enumerate}
 \item[(i)]
 If $\lim_{n\to\infty}P_n$ does not exist then
 $T$ is of type $III_1$. 
 \item[(ii)]
 If there exists the limit $P_+:=\lim_{n\to+\infty}P_n$ and $P_-:=\lim_{n\to-\infty}P_n$ then
 $P_+=P_-$.
 \item [(iii)]
 If $A=\{0,1\}$ and there exist the limit $Q:=\lim_{n\to\infty}P_n$ then $T$ is either of type $II_1$ or $II_\infty$.
More precisely,  $T$ is of type $II_1$ if and only if
$$
\sum_{n\in\Bbb Z}\sum_{a,b,a',b'\in A}\bigg(\sqrt{\widehat P_n(a,b)P_n(a',b')}-\sqrt{\widehat Q(a,b)Q(a',b')}\bigg)^2<\infty.
$$
The corresponding equivalent invariant probability measure (if exists) is the Markov measure
 defined by $Q$ and the distribution $\lambda$ on $A$ satisfying $\lambda Q=\lambda$.
 \end{enumerate}
\end{theorem}

It was also shown in \cite{AR} that an analogue of Theorem~\ref{ART}(iii) holds also for the golden mean Markov shift for which 
$A=\begin{pmatrix} 1 &0 &1\\
1&0&1\\
0&1&0
\end{pmatrix}$.

\section{Dynamical properties of nonsingular Poisson suspensions and nonsingular Gaussian
transformations}
\subsection{Nonsingular Poisson suspensions.} 
Let $(X,\mathcal B,\mu)$ be a $\sigma$-finite infinite standard measure space and let
$T\in\text{Aut}_2(X,\mu)$.
Then the Poisson suspension $(X^*,\mu^*,T_*)$ is well defined by Theorem~\ref{nonsing-P}.
The first problem to consider is to find out when $T_*$ admits an absolutely continuous
invariant probability measure. 
A satisfactory solution of this problem is obtained in \cite{DaKoRo1}.

\begin{theorem} The following are equivalent:
\begin{itemize}
\item
there exists a  $T_*$-invariant probability measure $\rho\prec\mu_*$,
\item   
$\sup_{n\in\Bbb Z}\|\sqrt{\frac{d\mu\circ T^n}{d\mu}}-1\|_2<\infty$,
\item
there is a measure $\kappa\prec\mu$ such that $\sqrt{\frac{d\kappa}{d\mu}}-1\in L^2(\mu)$
and $\kappa\circ T=\kappa$.
\end{itemize}
\end{theorem}

It follows that there exists ``properly'' nonsingular Poisson suspensions, i.e. suspensions that do not admit an absolutely continuous invariant measure.
The next problem is to find out when $T_*$ is  dissipative and when it is conservative.

\begin{theorem}
\begin{enumerate}
\item[(i)]
If 
$\sum_{n\in\Bbb Z}e^{-\frac12\|\sqrt{\frac{d\mu\circ T^n}{d\mu}}-1\|_2^2}<\infty$
then $T_*$ is totally dissipative  \cite{DaKoRo1}. 
\item[(ii)] 
Let  $T\in \text{Aut}_1(X,\mu)$, $\chi (T)= 0$ and
$\bigg(\frac{d\mu}{d\mu\circ T^{-n}}\bigg)^2-1\in L^1(X,\mu)$
for each $n\in\Bbb N$.
 If there is a sequence $(b_n)_{n=1}^\infty$
of positive reals such that
$\sum_{n=1}^\infty b_n=\infty$ but
$$
\sum_{n=1}^\infty b_n^2e^{\int_X\big(\big(\frac{d\mu}{d\mu\circ T^{-n}}\big)^2-1\big)d\mu}<\infty
$$
then $T_*$ is conservative  \cite{DaKoRo2}.
\end{enumerate}
\end{theorem}

We also note that if $T\in \text{Aut}_1(X,\mu)$ and $\chi (T)\ne 0$ then $T$ is not conservative
\cite{DaKoRo1}.
Moreover, $T_*$ is totally dissipative \cite{DaKoRo2}.

\begin{theorem}[\cite{DaKoRo1}] If $T$ is of 0-type and there is no $T$-invariant measure
$\kappa\prec\mu$ such that $\sqrt{\frac{d\kappa}{d\mu}}-1\in L^2(\mu)$ then $T_*$ is of 0-type.
\end{theorem}

The following theorem is proved via Baire category tools. 

\begin{theorem}[\cite{DaKoRo2}]
The set
$$
\{T\in \text{Aut}_2(X,\mu)\mid \text{$T$ and $T_*$ are both ergodic and of type $III_1$ }\}
$$
is a dense $G_\delta$ in $(\text{Aut}_2(X,\mu),d_2)$.
The set
$$
\{T\in \ker\chi \mid \text{$T$ and $T_*$ are both ergodic and of type $III_1$ }\}
$$
is a dense $G_\delta$ in $(\ker\chi,d_1)$.

\end{theorem}

It is easy to see that if $T\in \text{Aut}_2(X,\mu)$ then $T\in \text{Aut}_2(X,t\cdot\mu)$ for each $t>0$.
Hence $T_*$ is $(t\cdot\mu)^*$-nonsingular for each $t>0$.
However the dynamical properties of the systems $(X^*, T_*, (t\cdot\mu)^*)$ depend heavily on the choice of $t$.
This is illustrated by a concrete example constructed in \cite{DaKoRo2}.

\begin{ex}[Phase transition] There is a totally dissipative transformation $T\in\text{Aut}_2(X,\mu)$ and $t_0>0$ such that the dynamical system 
$(X^*, T_*, (t\cdot\mu)^*)$ is ergodic of type $III_1$ for each $t\in(0,t_0)$ and
and totally dissipative for each $t>t_0$.
\end{ex}

The point $t_0$ can be interpreted as a ``bifurcation point''.

We note that if $T\in\text{Aut}_2(X,\mu)$ and $T$ is totally dissipative
then $T_*$ is isomorphic to a nonsingular Bernoulli shift.
Therefore the following theorem was proved in \cite{DaKo}
simultaneously with~Theorem~\ref{I+Z}.

\begin{theorem} Let $K\in\{III_\lambda\mid 0\le \lambda\le 1\}\sqcup\{II_\infty\}$.
Then there is a totally dissipative 
transformation $T\in\text{Aut}_2(X,\mu)$ such that
$(X^*, T_*, \mu^*)$ is weakly mixing of type $K$.
\end{theorem}

This theorem is further refined in  type $III_0$ as follows (cf. Theorem~\ref{recent}).

\begin{theorem}[\cite{BeVa2}]
Let $V$ be an AT-flow associated with an 
ITPFI$_2$-transformation.
Then there is 
a totally dissipative 
transformation $T\in\text{Aut}_2(X,\mu)$ such that
$T_*$ is
 weakly mixing and  the associated flow of $T_*$ is isomorphic to $V$.
\end{theorem}

\subsection{Nonsingular Gaussian transformations.}
Let $\mathcal H$ be a separable infinite dimensional real Hilbert space.
Let $(h,O)\in\text{Aff}\,\mathcal H$.
For $n\in\Bbb Z$, define a vector  $h^{(n)}\in\mathcal H$ by setting $(h,O)^n=(h^{(n)},O^n)$.
We first consider when $G_{h,O}$ admits an equivalent
invariant probability measure.

\begin{theorem}[ \cite{AIM}, \cite{DanL22}] Let $(X,\mu)$ denote the space of $G_{h,O}$.
The following are equivalent:
\begin{itemize}
\item
there exists a  $G_{h,O}$-invariant probability measure $\rho\sim\mu$,
\item   
$h$ is an $O$-coboundary, i.e. there is $a\in\mathcal H$ such that $h=a-Oa$,
\item
the affine operator $(h,O)\in\text{{\rm Aff}}\,\mathcal H$  on $\mathcal H$ has a fixed point,
\item
the sequence $(h^{(n)})_{n\in\Bbb Z}$ is bounded in $\mathcal H$.
\end{itemize}
\end{theorem}

The following dichotomy for nonsingular Gaussian transformations was established in \cite{AIM} (see also \cite{DanL22}).

\begin{theorem} If $O$ has no nontrivial fixed vectors in $\mathcal H$ then
$G_{h,O}$ is either conservative or totally dissipative.
\end{theorem}

In \cite{AIM}, for each $(h,O)\in\text{{\rm Aff}}\,\mathcal H$, the authors consider a one-parametric family of nonsingular Gaussian transformations $G_{\theta h,O}$, $\theta\in(0,+\infty)$.

\begin{defn} The {\it Poincare exponent}  of  $(h,O)$ is
$$
\delta_{h,O}:=\inf\{\alpha>0\mid \sum_{n=1}^\infty e^{-\alpha \|h^{(n)}\|^2}<\infty\}\in[0,+\infty].
$$
\end{defn}

The following theorem demonstrates a phase transition from conservativity to total 
dissipativity.

\begin{theorem}[\cite{AIM}] 
Let $(h,O)\in\text{{\rm Aff}}\,\mathcal H$.
There exists $\theta_{diss}\in[0+\infty]$ such that $G_{\theta h,O}$ is conservative for all   $\theta<\theta_{diss}$ and  totally dissipative for all $\theta>\theta_{diss}$.
Moreover, 
$$
\sqrt{2\delta_{h,O}}\le \theta_{diss}\le2\sqrt{2\delta_{h,O}}.
$$
\end{theorem}

This result can be strengthened under  additional assumptions on $O$.

\begin{theorem}[\cite{AIM}, \cite{DanL22}]  Let $(h,O)\in\text{{\rm Aff}}\,\mathcal H$.
Suppose that  $f$ is not an $O$-coboundary
 and the  $\mu$-preserving transformation $G_{0,O}$ is mildly mixing.
Then for each $\theta<\theta_{diss}$, the nonsingular Gaussian transformation $G_{\theta h,O}$
is weakly mixing, IDPFT,  and  of Krieger type $III_1$.
\end{theorem}

This theorem is generalized  in a further work \cite{MaVa} for some cases in which 
$G_{0,O}$ is weakly mixing but not  mildly mixing.
In these cases  $G_{0,O}$ has  ``mixing  subsequences'' along which
the cocycle $(h^{(n)})_n$ is ``proper''.
In particular, an example of rigid weakly mixing Gaussian transformation $G_{0,O}$ is constructed
such that the nonsingular Gaussian transformation $G_{h,O}$ is ergodic of type $III_1$ for some $h\in\mathcal H$.

Currently, the only known examples of ergodic nonsingular Gaussian transformations are either of type $III_1$ or $II_1$.
Therefore,
though there is an obvious analogy between the theory of nonsingular Gaussian transformations and the theory of nonsingular Poisson suspensions, the latter theory looks more diverse. 
We illustrate this by the following theorem.

\begin{theorem}[\cite{Da23}]
For each  $K\in\{III_\lambda\mid 0\le \lambda\le 1\}\sqcup\{II_\infty\}$,
there is  transformation $T\in\text{{\rm Aut}}_2(Y,\nu)$ such that
the nonsingular Poisson suspension $T_*$ is 0-type, weakly mixing and of Krieger type $K$,
while the nonsingular Gaussian transformation $G_{\sqrt{\frac{d\nu\circ T}{d\nu}}-1, U_T}$
is 0-type, weakly mixing and  of Krieger type $III_1$.
The unitary Koopman operators of the two nonsingular transformations are unitarily equivalent.
\end{theorem}

\section {Spectral theory for nonsingular systems
 }
While the spectral theory for probability preserving systems is developed
in depth, the spectral theory of nonsingular systems is still in its
infancy. We discuss below some problems related to $L^\infty$-spectrum
which may be regarded as an analogue of the discrete spectrum. We also
include results on computation of the maximal spectral type of the
`nonsingular' Koopman operator for rank-one nonsingular transformations.

\subsection{$L^\infty$-spectrum and groups of quasi-invariance
}\label{9.1}
 Let $T$ be an ergodic nonsingular transformation of $(X,\mathcal
B,\mu)$. A number $\lambda\in\mathbb T$ belongs to the $L^\infty$-spectrum
$e(T)$ of $T$ if there is a function $f\in L^\infty(X,\mu)$ with $f\circ
T=\lambda f$. $f$ is called an $L^\infty$-{\it eigenfunction} of $T$
corresponding to $\lambda$. Denote by $\mathcal E(T)$ the group of all
$L^\infty$-eigenfunctions of absolute value $1$. It is a Polish group when
endowed with the topology of converges in measure. If $T$ is of type $II_1$
then the $L^\infty$-eigenfunctions are $L^2(\mu')$-eigenfuctions of $T$,
where $\mu'$ is an equivalent invariant probability measure. Hence $e(T)$
is countable. Osikawa constructed in \cite{Os} the first examples of ergodic
nonsingular transformations with uncountable $e(T)$.

We state now a nonsingular version of the von Neumann-Halmos discrete
spectrum theorem. Let $Q\subset \mathbb T$ be a countable infinite subgroup.
Let $K$ be a compact dual of $Q_d$, where $Q_d$ denotes $Q$ with the
discrete topology. Let $k_0\in K$ be the element defined by $k_0(q)=q$ for
all $q\in Q$. Let $R:K\to K$ be defined by $Rk=k+k_0$. The system $(K,R)$
is called a {\it compact group rotation}. The following theorem was proved in \cite{AN87}.

\begin{theorem}\label{T:4.1}  Assume that the $L^\infty$-eigenfunctions
of $T$ generate the entire $\sigma$-algebra $\mathcal B$. Then $T$ is
isomorphic to a compact group rotation equipped with an ergodic
quasi-invariant measure.
\end{theorem}

A natural question arises: which subgroups of $\mathbb T$ can appear as $e(T)$
for an ergodic $T$?

\begin{theorem}[\cite{MooS}, \cite{Aa}]\label{T:4.2}  $e(T)$ is a Borel subset of
$\mathbb T$ and carries a unique Polish topology which is stronger than the usual
topology on $\mathbb T$. The Borel structure of $e(T)$ under this topology
 agrees with the Borel structure inherited from $\mathbb T$. There is a Borel
 map $\psi:e(T)\ni\lambda\mapsto\psi_\lambda\in\mathcal E(T)$ such that
 $\psi_\lambda\circ T=\lambda\psi_\lambda$ for each
 $\lambda$. Moreover, $e(T)$ is of Lebesgue  measure $0$ and it can have an
 arbitrary Hausdorff dimension.
 \end{theorem}

A proper Borel subgroup $E$ of $\mathbb T$ is called
\begin{enumerate}
\item[(i)]
{\it weak Dirichlet} if $\limsup_{n\to\infty}\widehat\lambda(n)=1$ for each
finite complex measure $\lambda$ supported on $E$;
\item[(ii)]  {\it saturated} if $\limsup_{n\to\infty}|\widehat\lambda(n)|\ge
|\lambda(E)|$ for each finite complex measure $\lambda$ on $\mathbb T$,
\end{enumerate}
where $\widehat\lambda(n)$ denote the $n$-th Fourier coefficient of
$\lambda$.
Every countable subgroup of $\mathbb T$ is saturated.

\begin{theorem}\label{T:4.3}  $e(T)$ is $\sigma$-compact in the usual
topology on $\mathbb T$ \cite{HMP} and saturated (\cite{Mel}, \cite{HMP}).
\end{theorem}

It follows that $e(T)$ is weak Dirichlet (this fact was established earlier
in \cite{Sc3}).

It is not known if every Polish group continuously embedded in $\mathbb T$
as a $\sigma$-compact saturated group is the eigenvalue group of some
ergodic nonsingular transformation. This is the case for the so-called
$H_2$-groups and the groups of quasi-invariance of measures on $\mathbb T$
(see below). Given a sequence $n_j$ of positive integers and a sequence
$a_j\ge 0$, the set of all $z\in\mathbb T$ such that $\sum_{j=1}^\infty
a_j|1-z^{n_j}|^2<\infty$ is a group. It is called an $H_2$-{\it group}.
Every $H_2$-group is Polish in an intrinsic topology stronger than the
usual circle topology.

\begin{theorem}[\cite{HMP}]\label{T:4.4}

\begin{itemize}
\item[(i)] Every $H_2$-group is a saturated (and hence weak Dirichlet)
$\sigma$-compact subset of $\mathbb T$.
\item[(ii)]
If $\sum_{j=0}^\infty a_j=+\infty$ then the corresponding $H_2$-group is a
proper subgroup of $\mathbb T$.
\item[(iii)] If $\sum_{j=0}^\infty a_j(n_j/n_{j+1})^2<\infty$ then
the corresponding $H_2$-group is uncountable.
\item[(iv)] Any $H_2$-group is
$e(T)$ for an ergodic nonsingular compact group rotation $T$.
\end{itemize}
\end{theorem}

It is an open problem whether every eigenvalue group $e(T)$ is an
$H_2$-group.  It is known however that $e(T)$ is close `to be an $H_2$-group':
if a compact subset $L\subset\mathbb T$ is disjoint from $e(T)$ then there is an
$H_2$-group containing $e(T)$ and disjoint from $L$.

\begin{ex}[\cite{AN87}, see also \cite{Os}]\label{4.5} Let $(X,\mu,T)$ be the nonsingular 
product odometer associated to a sequence $(2,\nu_j)_{j=1}^\infty$
Let $n_j$ be a sequence of positive integers such that $n_j>\sum_{i<j}n_i$
for all $j$. For $x\in X$, we put $h(x):=n_{l(x)}-\sum_{j<l(x)}n_j$. Then
$h$ is a Borel map from $X$ to the positive integers. Let $S$ be the tower
over $T$ with height function $h$ (see \S\,\ref{tower}). Then $e(S)$ is the
$H_2$-group of all $z\in\mathbb T$ with $\sum_{j=1}^\infty
\nu_j(0)\nu_j(1)|1-z^{n_j}|^2<\infty$.
\end{ex}

It was later shown in \cite{HMP} that if $\sum_{j=1}^\infty
\nu_j(0)\nu_j(1)(n_j/n_{j+1})^2<\infty$ then the $L^\infty$-eigenfunctions of $S$
generate the entire $\sigma$-algebra, i.e., $S$ is isomorphic (measure
theoretically) to a nonsingular compact group rotation.

Let $\mu$ be a finite measure on $\mathbb T$. Let $H(\mu):=\{z\in\mathbb
Z\mid \delta_z*\mu\sim\mu\}$, where $*$ means the convolution of measures.
Then $H_\mu$ is a group called the {\it group of quasi-invariance of}
$\mu$. It has a Polish topology whose Borel sets agree with the Borel sets
which $H(\mu)$ inherits from $\mathbb T$ and the injection map of $H(\mu)$
into $\mathbb T$ is continuous. This topology is induced by the weak
operator topology on the unitary group in the Hilbert space $L^2(\mathbb
T,\mu)$ via the map $H(\mu)\ni z\mapsto U_z$, $(U_zf)(x)=\sqrt{(d(\delta_z
* \mu)/d\mu)(x)}f(xz)$ for $f\in L^2(\mathbb T,\mu)$.
Moreover, $H(\mu)$ is saturated \cite{HMP}. If $\mu(H(\mu))>0$ then either
$H(\mu)$ is countable or $\mu$ is equivalent to $\lambda_\mathbb T$
\cite{ManN}.

\begin{theorem}[\cite{AN87}]\label{T:4.6}
Let $\mu$ be an ergodic with respect to
the $H(\mu)$-action by translations on $\mathbb T$. Then there is a compact
group rotation $(K,R)$ and a finite measure on $K$ quasi-invariant and
ergodic under $R$ such that $e(R)=H(\mu)$. Moreover, there is a continuous
one-to-one homomorphism $\psi:e(R)\to E(R)$ such that $\psi_\lambda \circ
R=\lambda\psi_\lambda$ for all $\lambda\in e(R)$.
\end{theorem}

 It was shown by Aaronson and
Nadkarni \cite{AN87} that if $n_1=1$ and $n_j=a_ja_{j-1}\cdots a_1$ for
positive integers $a_j\ge 2$ with $\sum_{j=1}^\infty a_j^{-1}<\infty$ then the
transformation $S$ from Example~\ref{4.5} does not admit a continuous
homomorphism $\psi:e(S)\to E(S)$ with $\psi_\lambda \circ
T=\lambda\psi_\lambda$ for all $\lambda\in e(S)$. Hence $e(S)\ne H(\mu)$
for any measure $\mu$ satisfying the conditions of Theorem~\ref{T:4.6}.

Assume that $T$ is an ergodic nonsingular compact group rotation. Let
$\mathcal B_0$ be the $\sigma$-algebra generated by a sub-collection of
eigenfunctions. Then $\mathcal B_0$ is invariant under $T$ and hence a
factor (see \S\ref{S:joinings}) of $T$. It is not known if every factor of $T$ is of this form. It
is not even known whether every factor of $T$ must have non-trivial
eigenvalues.

\subsection{Koopman unitary operator for a nonsingular system
}
Let $(X,\mathcal B,\mu,T)$ be a nonsingular dynamical system. In this subsection we
consider spectral properties of the Koopman operator $U_T$ defined by~(\ref{koopman oper}). First, we note that the spectrum of
$T$ is the entire circle $\mathbb T$ \cite{Nad2}. Next, if $U_T$ has an
eigenvector then $T$ is of type $II_1$. Indeed, if there are
$\lambda\in\mathbb T$ and $0\ne f\in L^2(X,\mu)$ with $U_Tf=\lambda f$ then the
measure $\nu$, $d\nu(x):=|f(x)|^2d\mu(x)$, is finite, $T$-invariant and
equivalent to $\mu$. Hence if $T$ is of type $III$ or $II_\infty$ then the
maximal spectral type $\sigma_T$ of $U_T$ is continuous. Another `restriction' on $\sigma_T$ was  found in \cite{Roy2}: no Fo\"{\i}a\c{s}-Str\u{a}til\u{a} measure is absolutely continuous with respect to $\sigma_T$ if $T$ is of type  $II_\infty$. We recall that
a symmetric measure on $\mathbb T$ possesses {\it Fo\"{\i}a\c{s}-Str\u{a}til\u{a} property} if for each ergodic probability preserving system $(Y,\nu,S)$ and $f\in L^2(Y,\nu)$, if $\sigma$ is the spectral measure of $f$ then $f$ is a Gaussian random variable \cite{LPT}. For instance, measures supported on Kronecker sets possess this property.

As we have noted in \S\ref{mixing},  mixing (0-type) is an $L^2$-spectral property for  nonsingular transformations.
Also, if $T$ is infinite measure-preserving then $T$ is mixing if and only if
$n^{-1}\sum_{i=0}^{n-1}U_T^{k_i}\to 0$ in the strong operator topology for each strictly increasing sequence $k_1<k_2<\cdots$ \cite{KS69}. 
This generalizes a well known theorem of Blum and Hanson for probability preserving maps. For comparison, we note that ergodicity is not an $L^2$-spectral property of infinite measure-preserving systems.

Now let $T$ be a rank-one nonsingular transformation associated with a sequence $(r_n,w_n,s_n)_{n=1}^\infty$ as in \S\ref{S:rankone}.

\begin{theorem}[\cite{HMP}, \cite{ChN}]\label{T:4.7}  The spectral multiplicity
of $U_T$ is $1$ and the maximal spectral type $\sigma_T$ of $U_T$ (up to a
discrete measure in the case $T$ is of type II$_1$) is the weak limit of the measures $\rho_k$ defined as
follows:
$$
d\rho_k(z)=\prod_{j=1}^k w_j(0)|P_j(z)|^2\,dz,
$$
where $P_j(z):=1+\sqrt{w_j(1)/w_j(0)}z^{-R_{1,j}}+
\cdots+\sqrt{w_j(m_j-1)/w_j(0)}z^{-R_{r_j-1,j}}$, $z\in\mathbb T$,
$R_{i,j}:=ih_{j-1}+s_j(0)+\cdots+s_j(i)$, $1\le i\le r_k-1$ and $h_j$ is the hight of the $j$-th column.
\end{theorem}

Thus the maximal spectral type of $U_T$ is given by a so-called {\it
generalized Riesz product}. We refer the reader to \cite{HMP}, \cite{HMP1},
\cite{ChN}, \cite{Nad} for a detailed study of Riesz products:
their convergence, mutual singularity, singularity to $\lambda_\mathbb T$,
etc.

It was shown in \cite{AN87} that $H(\sigma_T)\supset e(T)$ for any ergodic
nonsingular transformation $T$. Moreover, $\sigma_T$ is ergodic under the
action of $e(T)$ by translations if $T$ is isomorphic to an ergodic
nonsingular compact group rotation. However it is not known:
\begin{enumerate}
\item[(i)]
Whether $H(\sigma_T)= e(T)$ for all ergodic  $T$?
\item[(ii)]
Whether ergodicity of $\sigma_T$ under $e(T)$ implies that $T$ is an
ergodic compact group rotation?
\end{enumerate}

The first claim of Theorem~\ref{T:4.7} extends to the rank $N$ nonsingular systems
as follows: if $T$ is an ergodic nonsingular transformation of rank $N$
 then the spectral multiplicity of
$U_T$ is bounded by $N$ (as in the finite measure-preserving case). It is
not known whether this claim is true  for a more general class of transformations which are defined as rank $N$ but without the assumption that the Radon-Nikodym cocycle is constant on the tower levels.

 Danilenko and Ryzhikov showed in \cite{DanR1} that
for  each subset $E\subset\Bbb N$, there is an ergodic conservative infinite measure-preserving transformation $T$  such that  the set of essential values of the multiplicity function of $U_T$ is $E$.
In a subsequent paper \cite{DanR2} they sharpened this result: for  each subset $E\subset\Bbb N\cup\{\infty\}$, there is a mixing  ergodic conservative infinite measure-preserving transformation $T$  such that  the set of essential values of the multiplicity function of $U_T$ is $E$.
We note that the analogous realization problem for  spectral multiplicities  of ergodic probability preserving transformations is still open \cite{Dan27S}.

In  \cite{DanR1}, a mixing rank-one infinite measure-preserving transformation $T$ was constructed such that the measures $\sigma_T,\sigma_T*\sigma_T, \sigma_T*\sigma_T*\sigma_T,\dots$ on $\Bbb T$ are mutually disjoint.
Hence the unitary operator $U_T\oplus U_T^{\odot 2}\oplus U_T^{\odot 3}\oplus\cdots$ has a simple spectrum.

El Abdalaoui and  Nadkarni constructed an  ergodic nonsingular transformation whose spectrum has Lebesgue component of multiplicity one \cite{AbNa}. 
The problem of existence of an ergodic nonsingular transformation  with a simple Lebesgue spectrum   is still open.

\subsection{  Koopman  unitary operators associated with nonsingular Poisson transformations}
We recall a very important structural feature of Poisson spaces (see \cite{Ner}):
 there is a canonical isometry between $L^2(\mu^*)$
 and the symmetric Fock space $F(L^2(\mu))$ over $L^2(\mu)$.
  We recall that
  $F(L^2(\mu)):=\bigoplus_{n=0}^\infty L^2(\mu)^{\odot n}$,
where $L^2(\mu)^{\odot 0}=\Bbb C$
and each factor $L^2(\mu)^{\odot n}$ is equipped with the normalized scalar product 
$n!\langle.,.\rangle_{L^2(\mu)^{\odot n}}$.
There is an exponential map $\mathcal E:L^2(\mu)\to F(L^2(\mu))$, given by
$$
\mathcal E(f):=\bigoplus_{n=0}^\infty\frac1{n!}f^{\odot n}.
$$
The family $\{\mathcal E(f)\mid f\in L^2(\mu)\}$ is total in $F(L^2(\mu))$.
Let  Aff$_\Bbb R( L^2(\mu))$ be the   subgroup of invertible affine operators in $L^2(\mu)$  that preserve invariant the $\Bbb R$-subspace $L^2_\Bbb R(\mu)$.  
Then Aff$_\Bbb R( L^2(\mu))=L^2_\Bbb R(\mu)\rtimes \mathcal U_{\Bbb R}(L^2(\mu))$, where 
$\mathcal U_{\Bbb R}(L^2(\mu))$ is the group of unitary operators that preserve invariant 
$L^2_\Bbb R(\mu)$.
The Weil unitary representation  $W$ of Aff$_\Bbb R( L^2(\mu))$ is well defined in $F(L^2(\mu))$ via the formula
$$
W(f,V)\mathcal E(h):=e^{-\frac12\|f\|_2-\langle f,Vh\rangle}\mathcal E(f+Vh),\qquad f\in L^2(\mu).
$$
Given $T\in\text{Aut}_2(X,\mu)$, we have that
$\Big(\sqrt{\frac{d\mu\circ T}{d\mu}}-1,U_T\Big)\in \text{Aff}_\Bbb R( L^2(\mu))$.

\begin{theorem}[Koopman operator associated with $T_*$ \cite{DaKoRo1}]\label{KoopmanP}
 If $T\in\text{Aut}_2(X,\mu)$ then under the canonical identification of  $L^2(\mu^*)$ and $F(L^2(\mu))$,
$$
U_{T_*}= W\bigg(\sqrt{\frac{d\mu\circ T}{d\mu}}-1,U_T\bigg).
$$
\end{theorem}

\subsection{  Koopman  unitary operators associated with nonsingular Gaussian transformations}
Let  $\mathcal H$ be a separable infinite dimensional real Hilbert space. 
Denote by $(X,\mu)$ the probability space where the nonsingular Gaussian transformations
$G_{h,O}$ for all $(h,O)\in\text{Aff} \,H$ are defined (see \S\ref{gaus}).  
It is well known that there is a canonical isometry between $L^2(X,\mu)$
and the symmetric Fock space $F(\mathcal H)$.

\begin{theorem}[Koopman operator associated with $G_{h,O}$ \cite{DanL22}]\label{KoopmanG}
 If $(h,O)\in\text{{\rm Aff}} \,H$ then under the canonical identification of  $L^2(X,\mu)$ and $F(\mathcal H)$,
$$
U_{G_{h,O}}= W\bigg(\frac12h,O\bigg).
$$
\end{theorem}

It follows from Theorems~\ref{KoopmanP} and \ref{KoopmanG} that  each nonsingular Poisson transformation $T_*$ is spectrally equivalent to the nonsingular Gaussian
transformation $G_{2\big(\sqrt{\frac{d\mu\circ T}{d\mu}}-1\big),U_T}$.

\section {Entropy and other invariants
 }
Let $T$ be an ergodic conservative nonsingular transformation of a
standard probability space $(X,\mathcal B,\mu)$. If $\mathcal P$ is a finite
partition of $X$, we define the entropy $H(\mathcal P)$ of $\mathcal P$ as $H(\mathcal
P)=-\sum_{P\in\mathcal P}\mu(P)\log\mu(P)$. In the study of measure-preserving
systems the classical (Kolmogorov-Sinai) entropy proved to be a very useful
invariant for isomorphism \cite{CFS}. The key fact of the theory is that if $\mu\circ
T=\mu$ then the limit $\lim_{n\to\infty}n^{-1}H(\bigvee_{i=1}^n
T^{-i}\mathcal P)$ exists for every $\mathcal P$. However if $T$ does not preserve
$\mu$, the limit may no longer exist. Some efforts have been made to extend
the use of entropy and similar invariants to the nonsingular domain.
These include Krengel's entropy of conservative  measure-preserving maps and its extension to nonsingular maps, Parry's entropy and Parry's nonsingular version of Shannon-McMillan-Breiman theorem, Poisson entropy, critical dimension by Mortiss and Dooley, etc. Unfortunately, these invariants are less informative than their classical counterparts and they are more difficult to compute.

\subsection{{Krengel's and  Parry's   entropies}\label{Krengel}
}
Let $S$ be a conservative measure-preserving transformation of a $\sigma$-finite measure space
$(Y,\mathcal E,\nu)$. The {\it Krengel entropy} \cite{Kre1} of $S$ is defined
by
$$
h_{\text{Kr}}(S)=\sup\{\nu(E)h(S_E)\mid 0<\nu(E)<+\infty\},
$$
where $h(S_E)$ is the Kolmogorov-Sinai entropy of $S_E$.
It follows from Abramov's formula for the entropy of induced transformation
that $h_{\text{Kr}}(S)=\mu(E)h(S_E)$ whenever $E$ {\it sweeps out}, i.e.,
$\bigcup_{i\ge 0}S^{-i}E=X$. A generic transformation from Aut$_0(X,\mu)$
has entropy 0. Krengel raised a question in \cite{Kre1}: does there exist a
zero entropy infinite measure-preserving $S$ and a zero entropy finite
measure-preserving $R$ such that $h_{\text{Kr}}(S\times R)>0$? This problem
was solved in \cite{DaR} (a special case was announced by Silva and Thieullen in
an October 1995 AMS conference (unpublished)):
\begin{enumerate}
\item[(i)]
if $h_{\text{Kr}}(S)=0$ and $R$ is distal then $h_{\text{Kr}}(S\times
R)=0$;
\item[(ii)]
if $R$ is not distal then there is a rank-one transformation $S$ with
$h_{\text{Kr}}(S\times R)=\infty$.
\end{enumerate}
We also note that if a conservative $S\in\text{Aut}_0(X,\mu)$ is squashable, i.e., it commutes with
another transformation $R$ such that $\nu\circ R=c\nu$ for a constant $c\ne
1$, then $h_{\text{Kr}}(S)$ is either 0 or $\infty$ \cite{ST2}.

Now let $T$
be a type $III$ ergodic transformation of
$(X,\mathcal B,\mu)$. Silva and Thieullen define an entropy $h^*(T)$ of $T$ by
setting $h^*(T):=h_{\text{Kr}}(\widetilde T)$, where $\widetilde T$ is the
Maharam extension of $T$ (see \S\,\ref{S:maharam}). Since $\widetilde T$ commutes
 with transformations which `multiply'  $\widetilde T$-invariant measure,
 it follows that $h^*(T)$ is either 0 or $\infty$.

  Let $T$ be the standard
 $III_\lambda$-odometer from Example~\ref{type III}(i). Then $h^*(T)=0$. The same is true
  for a so-called ternary product odometer associated with the sequence $(3,\nu_n)_{n=1}^\infty$, where $\nu_n(0)=\nu_n(2)
  =\lambda/(1+2\lambda)$ and $\nu_n(1)=\lambda/(1+\lambda)$
   \cite{ST2}. It is not known
  however whether every ergodic nonsingular product odometer  has
  zero entropy. On the other hand, it was shown in
\cite{ST2} that $h^*(T)=\infty$ for every $K$-automorphism.

The Parry entropy \cite{Pa3} of $S$ is defined by
$$
h_{\text{Pa}}(S):=\{H(S^{-1}\mathfrak F |\mathfrak F)\mid \mathfrak F\text{ is a $\sigma$-finite subalgebra of $\mathfrak B$ such that $\mathfrak F\subset S^{-1}\mathfrak F$}\}.
$$
Parry showed \cite{Pa3} that $h_{\text{Pa}}(S)\le h_{\text{Kr}}(S)$.
It is  still an open question whether the two entropies coincide. This is the case when  $S$ is of rank one (since $h_{\text{Kr}}(S)=0$) and when $S$ is quasi-finite  \cite{Pa3}. The transformation $S$ is called {\it quasi-finite} if there exists a subset of finite measure $A\subset Y$ such that the first return time partition $(A_n)_{n>0}$ of $A$ has finite entropy. We recall that $x\in A_n\iff n\text{ is the smallest positive integer such that }T^nx\in A$. An example of non-quasi-finite ergodic infinite measure-preserving transformation was  constructed  in \cite{AaP}.
A natural question is about existence of the maximal invariant $\sigma$-finite subalgebra of zero (Krengel or Parry) entropy.
Such an algebra is called {\it  the Krengel-Pinsker or the Parry-Pinsker} factor of $T$ respectively.
Existence of the Krengel-Pinsker factors was proved in \cite{AaP} for a special class
of quasi-finite transformations called {\it LLB}.
This result was extended in \cite{Jan} in the following way.

\begin{theorem} Let $T$ be an ergodic quasi-finite transformation.
Then either there is the Krengel-Pinsker factor of $T$ which is also the Parry-Pinsker and the Poisson-Pinsker  (see the next subsection below) factor of $T$ or $T$ is remotely infinite, i.e.,
there exists a sub-$\sigma$-algebra $\mathcal F\subset\mathcal B$ such that
$T^{-1}\mathcal F\subset\mathcal F$, $\bigvee_{n>0}T^n\mathcal F=\mathcal F$ and the subalgebra $\bigwedge_{n>0}T^{-n}\mathcal F$ does not contain  subsets of positive finite measure.
\end{theorem}

\subsection{Poisson entropy}\label{Poisson}
Poisson entropy for infinite measure-preserving transformations was introduced in \cite{Roy}.
Let $(X,\mu)$ be an infinite $\sigma$-finite space and let $T$ be a $\mu$-preserving invertible transformation of $X$.
The Poisson suspension $T_*$ of $T$ is well defined on a probability space $(X^*,\mu^*)$ and
$\mu^*\circ  T_*=\mu^*$ (see \S\ref{NPT}).
It is ergodic if and only if $T$ has no invariant sets of finite positive measure.
 It follows from Theorem~\ref{KoopmanP} that $U_{ T_*}$  is the `exponent' of $U_T$.
Hence, the maximal spectral type of $U_{ T_*}$ is $\sum_{n\ge 0}(n!)^{-1}(\sigma_T)^{*n}$, where $\sigma_T$ is a measure of the maximal spectral type of $U_T$.

 Now {\it the Poisson entropy} $h_{\text{Po}}(T)$ of $T$ is $h(T_*)$.
 The main question is: whether $h_{\text{Po}}(T)$  coincides with $h_{\text{Pa}}(T)$ or
 $h_{\text{kr}}(T)$?
 It was shown in \cite{Jan} that  $h_{\text{Pa}}(T)\le h_{\text{Po}}(T)$.
 If $T$ is quasi-finite or rank one then the three entropies of $T$ coincide \cite{Jan}. 
 If $T$ is the infinite Markov shift associated with a pair $(P,\pi)$ for  recurrent  and irreducible
 $P$  (see \S\ref{S:MarkSh}) then
 $$
 h_{\text{Kr}}(T)= h_{\text{Pa}}(T)= h_{\text{Po}}(T)=-\sum_{a\in A}\pi(a)\sum_{b\in A}P(a,b)\log P(a,b).
 $$
 If $\sigma_T$ is singular or $U_T$ has finite multiplicity then $h_{\text{Po}}(T)=0$ \cite{Jan}.
It was  also shown  in \cite{Jan} that given a nontrivial invariant $\sigma$-finite algebra $\mathcal F$ of $\mathcal B$, the  natural $\mathcal F$-relative version of Poisson entropy coincides with the relative (Krengel) entropy defined in \cite{DaR}.
Hence if Krengel's and the Poisson entropies  coinside on $T\restriction\mathcal F$ for some $\mathcal F$ then $h_{\text{Kr}}(T)=h_{\text{Po}}(T)$.
On the other hand, Janvresse and de la Rue constructed an ergodic conservative infinite measure-preserving transformation $T$ such that $h_{\text{Kr}}(T)=0$ but $h_{\text{Po}}(T)>0$ \cite{Jande}.

\begin{defn}
 An ergodic measure-preserving transformation $T$ of a $\sigma$-finite measure space
$(X,\mathcal B,\mu)$
  is said {\it to have totally positive Poisson entropy} if for  each $\sigma$-finite $T$-invariant sub-$\sigma$-algebra ${\mathcal F}\subset {\mathcal B}$, the Poisson entropy of the system $(X,\mathcal F,\mu\restriction\mathcal F,T)$ is strictly positive.
  \end{defn}
  
  We note that the Poisson suspension of the system $(X,\mathcal F,\mu\restriction\mathcal F,T)$ from the above definition 
  is canonically a factor of $(\widetilde X,\widetilde{\mathcal B},\widetilde\mu,\widetilde T)$.
  Such factors of $\widetilde T$ are called Poissonian. 
  Roy showed in \cite{Roy3} that
if $T$ has totally positive Poisson entropy then $T$ is of zero type.

 \begin{theorem}[Existence of the Poisson-Pinsker factor \cite{Roy3}]  Let $T$ be an ergodic  measure-preserving transformation of an infinite $\sigma$-finite measure space $(X,\mathcal B,\mu)$.
 Then either $T$ has totally positive entropy  and $\widetilde T$ is CPE or there is  a $\sigma$-finite $T$-invariant sub-$\sigma$-algebra ${\mathcal E}\subset {\mathcal B}$ such that the Poisson suspension of 
 $(X,\mathcal F,\mu\restriction\mathcal F,T)$ is the Pinsker factor of $\widetilde T$.
  \end{theorem}
 
 If $T$ has totally positive entropy then the maximal spectral type of $T$ is Lebesgue countable.
 If $h_{\text{Po}}(T)>0$  and $T$ possesses a Poisson-Pinsker factor then 
 the maximal spectral type of $T$ in the orthocomplement to the Poisson-Pinsker factor is Lebesgue countable \cite{Roy3}.

\subsection{Parry's generalization of Shannon-MacMillan-Breiman theorem}

Let $T$ be an ergodic  transformation of a standard non-atomic probability
space $(X,\sS,\mu)$. Suppose that $f\circ T\in L^1(X,\mu)$ if and only
if $f\in L^1(X,\mu)$. This means that there is $K>0$ such that
$K^{-1}<\frac{d\mu\circ T}{d\mu}(x)<K$ for a.a. $x$. Let $\mathcal P$ be a
finite partition of $X$. Denote by $C_n(x)$ the atom of $\bigvee_{i=0}^n
T^{-i}\mathcal P$ which contains $x$. We put $\omega_{-1}=0$. Parry shows in \cite{Par2} that
\begin{multline*}
 \frac{\sum_{j=0}^n\log\mu(C_{n-j}(T^jx))(\omega_j(x)-\omega_{j-1}(x))}
{\sum_{i=0}^n\omega_j(x)}\to\\
 H\bigg(P\mid\bigvee_{i=1}^\infty
T^{-1}P\bigg)-\int_X\log E\bigg(\frac{d\mu\circ
T}{d\mu}\mid\bigvee_{i=0}^\infty T^{-i}\mathcal P\bigg)\,d\mu
\end{multline*}
for a.a. $x$. Parry also shows that under the aforementioned conditions on
$T$,
\[
\frac 1n\bigg(\sum_{j=0}^nH\bigg(\bigvee_{i=0}^j T^{-j}\mathcal
P\bigg)-\sum_{j=0}^{n-1}H\bigg(\bigvee_{i=1}^{j+1} T^{-j}\mathcal
P\bigg)\bigg)\to H\bigg(\mathcal P\mid\bigvee_{i=1}^\infty T^{-i}\mathcal P\bigg).
\]

\subsection{{Critical dimension}} The critical dimension introduced
by Mortiss \cite{Mor} measures the order of growth for sums of
Radon-Nikodym derivatives. Let $(X,\sS,\mu,T)$ be an ergodic nonsingular dynamical system.
Given $\delta>0$, let
\begin{align}
 X_\delta &:=\{x\in
X\mid\liminf_{n\to\infty}\frac{\sum_{i=0}^{n-1}\omega_i(x)}{n^\delta}>0\}
\text{ and}\\
X^\delta &:=\{x\in
X\mid\liminf_{n\to\infty}\frac{\sum_{i=0}^{n-1}\omega_i(x)}{n^\delta}=0\}.
\end{align}
Then $X_\delta$ and $X^\delta$ are $T$-invariant subsets.

\begin{defn} [\cite{Mor}, \cite{DooM}]
The {\it lower critical dimension} $\alpha(T)$ of $T$ is
$\sup\{\delta\mid\mu(X_\delta)=1\}$. The {\it upper critical dimension}
$\beta(T)$ of $T$ is $\inf\{\delta\mid\mu(X^\delta)=1\}$.
\end{defn}

It was shown in \cite{DooM} that the lower and upper critical dimensions
are invariants for isomorphism of nonsingular systems. Notice also that
$$
\alpha(T)=\liminf_{n\to\infty}\frac{\log(\sum_{i=1}^n\omega_i(x))}{\log
n}\text{ and
}\beta(T)=\limsup_{n\to\infty}\frac{\log(\sum_{i=1}^n\omega_i(x))}{\log n}.
$$
Moreover, $0\le\alpha(T)\le\beta(T)\le 1$. If $T$ is of type $II_1$ then
$\alpha(T)=\beta(T)=1$. If $T$ is the standard $III_\lambda$-odometer from
Example~\ref{type III} then
$\alpha(T)=\beta(T)=\log(1+\lambda)-\frac{\lambda}{1+\lambda}\log\lambda$.

\begin{theorem}\label{5.2}
\begin{enumerate}
\item[(i)]
For every $\lambda\in[0,1]$ and every $c\in[0,1]$ there exists a
nonsingular product odometer of type $III_\lambda$ with critical dimension equal
to $c$ \cite{Mor1}.
\item[(ii)]
For every $c\in[0,1]$ there exists a nonsingular product odometer of type
$II_\infty$ with critical dimension equal to $c$ \cite{DooM}.
\end{enumerate}
\end{theorem}

Let $T$ be the nonsingular product odometer associated with a sequence $(m_n,\nu_n)_{n=1}^\infty$.  Let $s(n)=m_1\cdots m_n$
and let $H(\mathcal P_n)$ denote the entropy of the partition of the first $n$
coordinates with respect to $\mu$. We now state a nonsingular version of
Shannon-MacMillan-Breiman theorem for $T$ from \cite{DooM}.

\begin{theorem}\label{5.3}  Let $m_i$ be bounded from above. Then
\begin{enumerate}
\item[(i)]
$\alpha(T)=\liminf_{n\to\infty}\inf\frac{-\sum_{i=1}^n\log m_i(x_i)}{\log
s(n)}=\liminf_{n\to\infty}\frac{H(\mathcal P_n)}{\log s(n)}$ \ and
\item[(ii)]
$\beta(T)=\limsup_{n\to\infty}\inf\frac{-\sum_{i=1}^n\log m_i(x_i)}{\log
s(n)}=\limsup_{n\to\infty}\frac{H(\mathcal P_n)}{\log s(n)}$
\end{enumerate}
for a.a. $x=(x_i)_{i\ge 1}\in X$.
\end{theorem}

It follows that in the case when $\alpha(T)=\beta(T)$, the critical
dimension coincides with $\lim_{n\to\infty}\frac{H(\mathcal P_n)}{\log s(n)}$.
In \cite{Mor1} this expression (when it exists) was called {\it AC-entropy}
(average coordinate). It also follows from Theorem~\ref{5.3} that if $T$ is a product 
odometer of bounded type then $\alpha(T^{-1})=\alpha(T)$
 and $\beta(T^{-1})=\beta(T)$. In \cite{DooM2}, Theorem~\ref{5.3} was extended  to a subclass of Markov
 odometers. 
 Those results were further extended to  so-called $G$-measures on product spaces \cite{MDoo} and a class of Bratteli-Vershik systems with multiple edges \cite{DooH}.
 The critical dimensions for nonsingular Bernoulli shifts  (see \S\,\ref{S:hamachi}) were investigated in \cite{DooM3}:

\begin{theorem} For any $\epsilon>0$, there exists a nonsingular Bernoulli shift  $S$ from the Krengel class
 with $\alpha(S)<\epsilon$ and $\beta(S)>1-\epsilon$.
\end{theorem}

\subsection{Nonsingular restricted orbit equivalence}
In \cite{Mor2} Mortiss initiated study of a nonsingular version of Rudolph's restricted orbit equivalence \cite{Ru1}. This work is still in its early stages and does not yet deal with any form of entropy. However she introduced nonsingular orderings of orbits,  defined sizes and showed  that much of the basic machinery  still works in the nonsingular setting.

\section{Nonsingular joinings and factors}\label{S:joinings}
The theory of joinings is a powerful tool to study probability preserving
systems and to construct striking counterexamples. 
It is interesting to study  what part
of this machinery can be extended to the nonsingular case.
 However, there are some principal obstacles for such extensions:
 \begin{itemize}
 \item there are   too many quasi-invariant measures in view of the Glimm-Effros theorem (see Theorem~\ref{T:Glimm});
 \item 
 ergodic components of a non-ergodic joining need not be  joinings of the original systems.
 \end{itemize}
There are several ways to bypass  these obstacles.
The principal idea is to select always an appropriate  (rather narrow) class of quasiinvariant measures under consideration or impose some restrictions on the structure of joinings.
 This approach led  to some
progress  in understanding $2$-fold joinings and constructing
prime systems of any Krieger type. As far as we know the higher-fold
nonsingular joinings have not been considered so far. It turned out however
that an alternative coding technique, predating joinings in studying the
centralizer and factors of the classical measure-preserving Chac{\'o}n maps, can be used as well to
classify factors of Cartesian products of some nonsingular Chac{\'o}n maps.

\subsection{Joinings, nonsingular MSJ and simplicity}   In this subsection all measures are probability measures. A {\it nonsingular joining } of two nonsingular systems $(X_1,\sS_1,\mu_1, T_1)$ and $ (X_2,\sS_2,\mu_2, T_2)$   is a  measure $\hmu$ on the product $\sS_1\times \sS_2$ that is nonsingular for $T_1\times T_2$ and satisfies:
  $\hmu(A\times X_2)=\mu_1(A)$ and $\hmu(X_1\times B)=\mu_2(B)$ for all $A\in\sS_1$ and $B\in\sS_2$. Clearly, the product $\mu_1\times\mu_2$ is a nonsingular joining. Given a transformation $S\in C(T)$, the measure $\mu_S$ given by $\mu_S(A\times B):=\mu(A\cap S^{-1}B)$ is a nonsingular joining of $(X,\mu,T)$ and $(X,\mu\circ S^{-1},T)$. It is called
 a {\it graph}-joining since it is supported on the graph of $S$. Another important kind of joinings that we are going to define now is related to factors of dynamical systems. Recall that
 given a  nonsingular  system
$  (X,{\cB} , \mu ,T)$, a  sub-$\sigma$-algebra $\mathcal A$ of $\sS$
 such that $T^{-1}({\mathcal A}) = {\mathcal A} \mod\mu$ is called a {\it factor} of $T$. There is another, equivalent, definition. A nonsingular
dynamical system $(Y,\mathcal C,\nu,S)$ is called a factor of $T$ if
there exists a  measure-preserving map $\phi : X \to
Y$,  called a  {\it factor map}, with $\phi
T = S\phi$ a.e.  (If $\phi$ is only nonsingular,
$\nu$ may be replaced  with the equivalent measure $\mu \circ \phi^{-1}$, for which $\phi$ is measure-preserving.) Indeed, the sub-$\sigma$-algebra
$\phi^{-1}(\mathcal C)\subset\sS$ is $T$-invariant and, conversely, any $T$-invariant
sub-$\sigma$-algebra of $\sS$ defines a factor map by immanent properties of standard probability spaces, see e.g.  \cite{Aa}.
If $\phi$ is a factor map as above,  then $\mu$
has a disintegration with  respect to $\phi$, i.e.,
$\mu \ = \ \int \mu_y d\nu (y)$ for
 a measurable map
$y\mapsto \mu_{y}$ from $Y$ to the probability measures on $X$ so that  $\mu_y(\phi^{-1}( y))=1$, the measure
 $\mu_{S\phi(x)}\circ T $ is equivalent to
 $\mu_{\phi(x)}$ and
\begin{equation}\label{10-1}
\frac{d\mu \circ
T}{d\mu} (x)
= \frac{d\nu \circ
S}{d\nu} (\phi(x)) {\frac{d\mu_{S\phi(x)}\circ
T}{d\mu_{\phi(x)}}} (x)
\end{equation}
 for a.e. $x\in X. $
Define now
the {\it relative product}
  $\hmu = \mu \times_{\phi} \mu$ on $X \times X$ by setting
$
\hmu  = \int \mu_y \times \mu_y \, d\nu (y).
$
Then it is easy to deduce from (\ref{10-1}) that
 $\hmu$ is a nonsingular self-joining of $T$.

We note however that the above definition of  joining is  too general to be satisfactory (as we noted in the introduction to this section). 
It does not reduce to the classical definition when we consider probability preserving systems.
 Indeed, the following result was proved in \cite{RS89}.

\begin{theorem}\label{T:njext} Let  $(X_1,\sS_1,\mu_1, T_1)$ and $ (X_2,\sS_2,\mu_2, T_2)$    be two finite measure-preserving systems such that $T_1 \times T_2$ is ergodic. Then  for every $\lambda, 0 < \lambda < 1$, there exists a nonsingular joining   $\hmu$  of $\mu_1$ and $\mu_2$ such that $(T_1 \times T_2, \hmu)$ is ergodic and of type $III_\lambda$.
\end{theorem}

It is not known however if the nonsingular joining $\hmu$ can be chosen in every orbit equivalence class. In view of the above, Rudolph and Silva \cite{RS89}
isolate an important subclass of joining.
It is used in the definition of a nonsingular version of minimal self-joinings.

\begin{defn} \label{E:rationaleq}
\begin{itemize}
\item[(i)] A nonsingular joining $\hmu$ of $(X_1,\mu_1,T_1)$ and
$(X_2,\mu_2,T_2)$ is    {\it rational}  if there exit
measurable
functions $c^{1}:X_1\to \mathbb R_+$ and $c^{2}:X_2\to\mathbb R_+$ such that
$$
\home^{\hmu}_1 (x_1, x_2)  =  \omega^{\mu_1}_1 (x_1 )
\omega^{\mu_2}_1 (x_2 ) c^{1} (x_1 )  =  \omega^{\mu_1}_1 (x_1 ) \omega^{\mu_2}_1 (x_2 )
c^{2}
(x_2 )\quad \hmu \ a.e.
$$
\item[(ii)]
A nonsingular dynamical system $ (X,{\cB} , \mu, T) $ has {\it minimal self-joinings (MSJ)}  over a class $\mathcal M$ of probability measures equivalent to $\mu$, if for every $\mu_1, \mu_2\in\mathcal M$, for every rational joining  $\hmu$ of $\mu_1, \mu_2$,  a.e. ergodic component of $\hmu$ is
either the product of its marginals or is the graph-joining supported on $T^j$ for some $j\in\mathbb Z$.
\end{itemize}
\end{defn}

Clearly, product measure, graph-joinings and the relative products are all rational joinings.   Moreover, a rational joining of finite measure-preserving systems
is measure-preserving and a rational joining of type $II_1$'s is of type $II_1$ \cite{RS89}. Thus we obtain the finite measure-preserving theory as a special case.   As for the definition of MSJ, it depends on a class $\mathcal M$ of equivalent measures. In the finite measure-preserving case $\mathcal M=\{\mu\}$. However, in the nonsingular case no particular measure is distinguished. We note also that Definition~\ref{E:rationaleq}(ii) involves some restrictions on all rational joinings and not only ergodic ones as in the finite measure-preserving case. The reason is that an ergodic component of a nonsingular joining needs not be a joining of measures equivalent to the original ones \cite{a87}.  For finite measure-preserving transformations, MSJ over $\{\mu\}$ is the same as the usual 2-fold MSJ \cite{JR87}.

A nonsingular transformation $T$ on $(X,\sS,\mu)$ is called {\it prime} if its only factors are $\sS$ and $\{X,\emptyset\}$ mod\,$\mu$.
A (nonempty) class  $\mathcal M$  of probability measures equivalent to $\mu$ is said to be
 {\it  centralizer stable} if
for each $S\in C(T)$ and $\mu_1\in\mathcal M$, the measure  $\mu_1\circ S$ is in $\mathcal M$.

\begin{theorem}[\cite{RS89}] \label{T:msj}Let $ (X,{\cB} , \mu, T) $ be a ergodic non-atomic dynamical system such that  $T$ has MSJ over a class $\mathcal M$
that is centralizer stable. Then $T$ is prime and  the centralizer of $T$  consists of the powers of $T$.
\end{theorem}

  A question that arises is whether  such nonsingular
dynamical system (not of type $II_1$) exist. Expanding on Ornstein's
original construction from \cite{Ornst},
  Rudolph and Silva
   construct in \cite{RS89}, for each $0\leq \lambda\leq 1$,
a nonsingular rank-one transformation
$T_\lambda$ that is of type $III_\lambda$    and that has MSJ over a class
$\mathcal M$ that is
centralizer stable. Type $II_\infty$ examples with analogues properties were also constructed there. In this connection it is worth to mention the example by
Aaronson and Nadkarni \cite{AN87} of $II_\infty$ ergodic transformations
that have no factor algebras on which the invariant measure is $\sigma$-finite (except for the entire ones); however these transformations are not prime.

A more general notion than MSJ called {\it graph self-joinings (GSJ)}, was
introduced \cite{SW92}: just replace the the words ``on $T^j$ for some
$j\in\mathbb Z$'' in Definition~\ref{E:rationaleq}(ii) with ``on $S$ for some
element $S\in C(T)$''. For finite measure-preserving transformations, GSJ
over $\{\mu\}$ is the same as the usual 2-fold simplicity \cite{JR87}. The
famous Veech theorem on factors of 2-fold simple maps (see \cite{JR87})
 was
extended to nonsingular systems in \cite{SW92} as follows: if a system
$(X,\sS,\mu,T)$ has GSJ then for every non-trivial factor $\mathcal A$ of
$T$ there exists a locally compact subgroup $H$ in $C(T)$ (equipped with
the weak topology) which acts smoothly (i.e., the partition into $H$-orbits
is measurable) and such that $\mathcal A=\{B\in\sS\mid \mu(hB\triangle B)=0
\text{ for all } h\in H\}$. It follows that there is a cocycle $\phi$ from
$(X,\mathcal A,\mu\restriction\mathcal A)$ to $H$ such that $T$ is
isomorphic to the $\phi$-skew product extension $(T\restriction A)_\phi$
(see \S\,6.4). Of course, the ergodic nonsingular product odometers and, more
generally, ergodic nonsingular compact group rotation (see \S\,\ref{9.1}) have
GSJ. However, except for this trivial case (the Cartesian square is
non-ergodic) plus the systems with MSJ from \cite{RS89}, no examples of
type $III$ systems with GSJ are known. In particular, no smooth examples have
been constructed so far. This is in sharp contrast with the finite measure
preserving case where abundance of simple (or close to simple) systems are
known (see \cite{JR87}, \cite{Tho}, \cite{Dan26}).

\subsection{Nonsingular coding and factors of Cartesian products of nonsingular maps}
As we have already noticed above, the nonsingular MSJ theory was developed in \cite{RS89} only for 2-fold self-joinings. The reasons for this were technical problems with extending the notion of rational joinings form $2$-fold to $n$-fold self-joinings. However while the $2$-fold nonsingular MSJ or GSJ properties of $T$ are sufficient to control the centralizer and the factors of $T$, it is not clear whether
 it implies anything about the factors or centralizer of $T\times T$. Indeed, to control them  one needs to know the $4$-fold joinings of $T$.
However even in the finite measure-preserving case it is a long standing open question whether $2$-fold MSJ implies $n$-fold MSJ. That is why del Junco and Silva \cite{JS03} apply an alternative---nonsingular coding---techniques to classify the factors of Cartesian products of nonsingular Chac\'on maps.
The techniques were originally used in \cite{j78} to show that the classical Chac\'on map is prime and has trivial centralizer. They were extended to nonsingular systems in \cite{JS}.

For each $0<\lambda<1$ we denote by $T_{\lambda}$
the Chac\'on map (see \S\,\ref{S:rankone}) corresponding the sequence
of probability
vectors  $w_n=
({\lambda}/({1+2\lambda}),
 {1}/({1+2\lambda}),
 {\lambda} /({1+2\lambda}))$ for all $n>0$. One can verify that the maps $T_{\lambda}$
are of type $III_\lambda$. (The classical Chac\'on map corresponds to $\lambda=1$.) All of these transformations are defined on the same standard Borel space $(X,\sS)$.
These transformations  were shown to be  power weakly mixing in \cite{afs01}.
 The centralizer of any finite Cartesian product of
nonsingular Chac\'on maps is computed in the following theorem.

\begin {theorem}[\cite{JS03}]{\label{secondthm}}  Let
$0<\lambda_{1}<\ldots <\lambda_{k} \leq 1$ and $n_{1},\ldots, n_{k}$
be positive integers.   Then the centralizer of the Cartesian product
$T_{\lambda_1}^{\otimes n_1} \times \ldots\times
T_{\lambda_k}^{\otimes n_k} $
is generated by maps of the form $U_1 \times \ldots \times U_k$,
where each $U_i$, acting on the $n_{i}$-dimensional
product space $X^{n_i}$, is a Cartesian
product of
powers of $T_{\lambda_{i}}$ or a co-ordinate permutation on
$X^{n_i}$.
\end{theorem}

Let $\pi$ denote the permutation on $X\times X$ defined by $\pi(x,y)=(y,x)$ and let $\sB$ denote the symmetric factor, i.e.,
$\sB=\{A\in\sS\otimes\sS\mid \pi(A)=A\}$.
The following theorem
  classifies  the factors of the Cartesian product of any
two nonsingular type $III_{\lambda}$, $0<\lambda < 1$, or type
$II_{1}$ Chac\'on maps.

\begin{theorem}[\cite{JS03}]{\label{firstthm} }  Let
$ T_{\lambda_{1}}$
and $T_{\lambda_{2}}$
be two nonsingular Chac\'on systems.  Let $\cF$ be a factor algebra of
$T_{\lambda_{1}} \times T_{\lambda_{2}}$.
\begin{itemize}
\item[(i)] If $\lambda_{1}\neq \lambda_{2}$ then
$\cF$ is equal $\mod 0$ to one of the four
algebras
 ${\cB}\otimes \cB$, ${\cB}\otimes \cN$, $\cN
 \otimes \cB$,   or $\cN\otimes \cN$, where $\cN=\{\emptyset, X\}$.
 \item[(ii)] If $\lambda_{1} = \lambda_{2}$ then
$\cF$ is equal $\mod 0$ to one of the
following algebras
 ${\cB}\otimes \cB$, ${\cB}\otimes \cN$, $\cN\otimes \cB$,
${\cN}\otimes{\cN}$,
  or $(T^m\times Id)\sB$ for some integer $m$.
\end{itemize}
\end{theorem}

It is not hard to obtain  type $III_{1}$ examples of Chac\'on maps for which the previous two theorems hold.
 However the construction of type
$II_{\infty}$ and type $III_{0}$ nonsingular Chac\'on transformations is
more subtle as it needs the choice of $\omega_n$ to vary with $n$. In
\cite{hs00}, Hamachi and Silva construct type $III_0$ and type $II_\infty$
examples, however the only property proved for these maps is ergodicity of
their Cartesian square. More recently, Danilenko \cite{d04} has shown that
all of them (in fact, a wider class of nonsingular Chac\'on maps of all
types) are power weakly mixing.

In \cite{cep},
  Choksi, Eigen and Prasad asked whether there exists a zero
 entropy, finite measure-preserving mixing automorphism $S$, and a
 nonsingular type $III$ automorphism $T$, such that $T\times S$ has no
 Bernoulli factors.
   Theorem~\ref{firstthm} provides a partial answer  (with a mildly mixing only instead of mixing) to this question:
    if $S$ is the finite measure-preserving  Chac\'on
    map and $T$ is a  nonsingular Chac\'on map as
   above, the factors of $T\times S$ are only the trivial
ones,
    so $T\times S$ has no Bernoulli factors.

\subsection{Joinings and MSJ for infinite measure-preserving systems}
Adams, Friedman and Silva introduced in \cite{AFS} an infinite version of Chac\'on map  $T$ as a rank-one transformation associated with $(r_n,\omega_n,s_n)_{n=1}^\infty$ such that
 $r_n=3$, $\omega_n(0)=\omega_n(1)=\omega_n(2)$, $s_n(0)=0$, $s_n(1)=1$ and $s_n(2)=3h_n+1$ for each $n>0$. This is called the {\it infinite Chac\'on} transformation. 
 Let $(X,\mu)$ be the space of  $T$.
 Of course, $\mu(X)=\infty$. This  transformation has infinite ergodic index \cite{AFS}, is not power weakly mixing and not multiply recurrent \cite{Ketal03}, and has trivial centralizer \cite{JanRR}.
 For each $d>0$, Janvresse, de la Rue and Roy investigated  $T^{\times d}$-invariant measures on $X^d$ which are
 {\it boundedly finite}.
 This means that  for each $d$   levels of every tower of the inductive construction, the measure of the Cartesian product of these levels is finite.
 The product $\bigotimes_{n=1}^d\mu$ and graph-joinings, i.e., measures of the form $(A_1,\dots, A_d)\mapsto \mu(S_1^{-1}A_1,\dots \cap S_d^{-1}A_d)$ for some transformations $S_1,\dots, S_d\in C(T)$, are boundedly finite.
Moreover,  $T$ itself  is uniquely ergodic in the sense that there is only one (up to scaling) boundedly finite $T$-invariant measure.
It was shown in \cite{JanRR} that each ergodic $T^{\times d}$-invariant boundedly finite measure is a direct product of  so-called {\it diagonal measures}.
 Unlike the finite measure-preserving case, the class of diagonal measures does not reduce to the graph-joinings (with $S_1,\dots, S_d$ being the powers of $T$).
 It contains  so-called {\it weird} measures whose  marginals are singular to $\mu$. 
 As a corollary, it was proved that $C(T)=\{T^n\mid n\in\Bbb Z\}$ \cite{JanRR}.
 Some of the weird measures are totally dissipative (supported on a single orbit) are some of them are conservative.
Danilenko showed in \cite{Dan*} that there is a conservative $T\times T$-invariant boundedly finite measure with absolutely continuous marginals whose ergodic components are all weird.
This phenomenon is impossible for another  infinite version $T$ of Chacon map constructed in \cite{JanRR2}.
Its construction mimics the construction of the classical Chacon map so much that it gives a $\mu$-conull subset $X_\infty$ such that for each $d\ge 1$, each ergodic $T^{\times d}$-invariant measure supported on $X_\infty^d$ is the direct product of several copies of $\mu$ and the graph-joinings generated by powers of $T$.
As a corollary, we obtain that   each boundedly finite $d$-fold self-joining  of $T$  (the marginals of a joining are absolutely continuous) is a convex combination of countably many ergodic joinings.

In \cite{Dan*}, the problems studied in \cite{JanRR} are considered from a different point of view.
Let $T$ be a homeomorphism of  a locally compact Cantor space $X$.
We assume that $T$ is Radon uniquely ergodic, i.e., there is only one (up to scaling) Radon $T$-invariant measure $\mu$ on $X$.
A {\it $d$-fold Radon self-joining of $T$} is a Radon measure on $X^d$ whose marginals (which may be non-sigma-finite) are equivalent to $\mu$.
We consider only Radon invariant measures, define Radon $d$-fold MSJ and Radon disjointness.
Of course, each ergodic component of nonergodic Radon joinning is Radon. However it needs not to be a joining.
Then the  $(C,F)$-construction  (see \cite{d01} and \cite{Dan26}) is used  to produce a number of rank-one homeomorphisms of $X$ whose ergodic joinings are explicitly described.
The weird measures from \cite{JanRR} appear now as a quasi-graph Radon measures, i.e., they are graphs of equivariant maps whose domain and range are meager (and of zero measure) subsets of $X$.
It is constructed an uncountable family of  pairwise Radon disjoint infinite Chacon like Radon uniquely ergodic homeomorphisms  with Radon MSJ.
Moreover, every transformation of this family is Radon disjoint with its inverse \cite{Dan*}.

 \section {Smooth nonsingular transformations
 }
Diffeomorphisms of smooth manifolds equipped with smooth measures are commonly considered as physically natural examples of dynamical systems. Therefore the construction of smooth models for various dynamical properties is a well established problem of the  modern (probability preserving) ergodic theory. Unfortunately, the corresponding `nonsingular' counterpart of this problem is almost unexplored. We survey here several interesting facts related to the topic.

For $r\in\mathbb N\cup\{\infty\}$, denote by Diff$^r_+(\mathbb T)$ the group of
orientation preserving $C^r$-diffeomorphisms of the circle $\mathbb T$. Endow
this set with the natural Polish topology. Fix $T\in\text{Diff}_+^r(\mathbb
T)$. Since $\mathbb T=\mathbb R/\mathbb Z$, there exists a $C^1$-function
$f:\mathbb R\to\mathbb R$ such that $T(x+\mathbb Z)=f(x)+\mathbb Z$ for all $x\in\mathbb R$.
The {\it rotation number} $\rho(T)$ of $T$ is the limit
$\lim_{n\to\infty}(\underbrace{f\circ\cdots\circ f}_{n\text{
times}})(x)\pmod 1$. The limit exists and does not depend on the choice of
$x$ and $f$. It is obvious that $T$ is nonsingular with respect to
Lebesgue measure $\lambda_\mathbb T$. Moreover, if $T\in \text{Diff}_+^r(\mathbb
T)$ and $\rho(T)$ is irrational then the dynamical system $(\mathbb
T,\lambda_\mathbb T,T)$ is ergodic \cite{CFS}. It is interesting to ask: which
Krieger's type can such systems have?

 Katznelson  showed in \cite{Kat} that the subset of type $III$
 $C^\infty$-diffeomorphisms
 and the subset of type $II_\infty$  $C^\infty$-diffeomorphisms are dense
 in Diff$_+^\infty(\mathbb T)$.
 Hawkins and Schmidt refined the idea of Katznelson  from \cite{Kat} to
 construct, for every irrational number $\alpha\in[0,1)$ which is not of
constant type (i.e., in whose continued fraction expansion the denominators
are not bounded) a transformation $T\in \text{Diff}_+^2(\mathbb T)$ which is
of type $III_1$ and $\rho(T)=\alpha$ \cite{HaS}. It should be
mentioned that class $C^2$ in the construction is essential, since it
follows from a remarkable result of Herman that if
$T\in\text{Diff}_+^3(\mathbb T)$ then under some condition on $\alpha$ (which
determines a set of full Lebesgue measure), $T$ is measure theoretically
(and topologically) conjugate to a rotation by $\rho(T)$ \cite{Her2}. Hence $T$ is of
type $II_1$.

In \cite{Haw3}, Hawkins shows that every smooth paracompact manifold of
dimension~$\geq 3$ admits a type $III_\lambda$ diffeomorphism for every
$\lambda\in[0,1]$. This extends a result of Herman \cite{Her} on the
existence of type $III_1$ diffeomorphisms in the same circumstances.

It is also of interest to ask: which free ergodic flows are
associated with smooth dynamical systems of type $III_0$? Hawkins proved
that any free ergodic $C^\infty$-flow on a smooth,
connected, paracompact manifold is the associated flow for a
$C^\infty$-diffeomorphism on another manifold (of higher dimension)
\cite{Haw4}.

A nice result was obtained in \cite{Kat2}: if $T\in\text{Diff}_+^2(\mathbb T)$
and the rotation number of $T$ has unbounded continued fraction
coefficients then $(\mathbb T, \lambda_\mathbb T, T)$ is ITPFI. 
Moreover, a
converse also holds: given a nonsingular product odometer $R$, the set of
orientation-preserving $C^\infty$-diffeomorphisms of the circle which are
orbit equivalent to $R$ is $C^\infty$-dense in the Polish set of all
$C^\infty$-orientation-preserving diffeomorphisms with irrational rotation
numbers. 
In contrast to that, Hawkins constructs in \cite{Haw2} a type
$III_0$ $C^\infty$-diffeomorphism of the 4-dimensional torus which is not
ITPFI.

Examples of $n$-to-1 conservative ergodic nonsingular C$^\infty$-endomorphisms on the 2-torus, not admitting an equivalent $\sigma$-finite invariant measure,  were constructed in \cite{HS91}.  
 In \cite{AB07} it is shown that a C$^1$ generic expanding map of $\mathbb T$ has no absolutely continuous $\sigma$-finite invariant measure.

Kosloff in \cite{Ko15}  showed that $\Bbb T^2$ admits a $C^1$ Anosov diffeomorphism of   type $III_1$ with respect to Lebesgue measure.
We recall that this phenomenon is impossible in the class of conservative $C^{1+\alpha}$
Anosov diffeomorphisms because by a theorem of Gurevich and Oseledets, every such transformation is of type $II_1$ (with respect to Lebesgue measure).
In a later work \cite{Ko20}, he extended this result to  $\Bbb T^d$ for every  $d>3$.
The case $d=3$ remains open.

\section{Miscellaneous topics}
Let $T$ be an ergodic  measure-preserving transformation of an infinite $\sigma$-finite nonatomic measure space
$(X,\mathcal B,\mu)$.

\subsection{On normalizing constants for ergodic theorem}\label{NC}
Replacing $\mu$ with an equivalent probability measure one can deduce from the Hurewicz ergodic theorem  that the average
$\frac{1}{n}\sum_{i=0}^{n-1}f(T^ix)$ converges to $0$ a.e. 
for each  function $f\in L^1(X,\mu)$.
In view of that, a natural question arises: is there  a sequence of positive numbers $(a_n)_{n=1}^\infty$ such that

\begin{align}\label{aarn}
\frac{1}{a_n}\sum_{i=0}^{n-1}f(T^ix) \to \int f d\mu
\  \text{a.e.}
\end{align}
 for each $f\in L^1(X,\mu)$?
Aaronson answered this question negatively in \cite{Aa77b}.
He showed that if 
 there is a sequence of positive numbers $(a_n)_{n=1}^\infty$ and a single integrable function $f\geq 0$ with $ \int f d\mu>0$ such that (\ref{aarn}) holds then
$\mu(X)<\infty$.
Thus,  no normalizing constants in the ergodic theorem for infinite measure-preserving transformations exist.
 Since then other forms of convergence of ergodic averages for a given sequence have been studied, for which  the reader may refer to \cite{Aa81}, \cite{Aa}, \cite{ThZw06}, \cite{AZ14} and  the references therein.

\subsection{Around King's weak closure theorem} We recall that if $S$ is probability preserving rank-one map then $C(S)$ is the weak closure of the set $\{S^n\mid n\in\Bbb Z\}$ (King theorem, \cite{King}).
It is still unclear whether this theorem extends to   the infinite measure-preserving rank-one transformations.
However, 
there are some classes of infinite rank-one maps for which it is true: zero type maps and partially bounded maps.

\begin{defn} A {\it $\sigma$-finite self-joining (of order 2)} of $T$ is a $\sigma$-finite $T\times T$-invariant measure $\lambda$ on $(X\times X,\mathcal B\otimes\mathcal B)$ such that $\lambda(A\times X)=\lambda(X\times A)=\lambda (A)$ for all $A\in\mathcal B$ of finite measure.
If for each ergodic $\sigma$-finite self-joining $\lambda$ of $T$, there is $n\in\Bbb Z$ such that
$\lambda(A\times B)=\mu(A\times T^{-n}B)$ then  $T$ is said to have {\it minimal $\sigma$-finite self-joinings (of order 2)}. 
\end{defn}

The above concept  of MSJ permits to control the $\mu$-preserving centralizer $C_0(T)$ of $T$:
if $T$ has MSJ then $C_0(T)=\{T^n\mid n\in\Bbb Z\}$.
It was shown by Ryzhikov and Thouvenot  \cite{RyT} that each zero type transformation of rank one has $\sigma$-finite MSJ.
Since the rank-one transformations are non-squashable, it follows that the centralizer of each zero type rank-one transformation is just its powers.

\begin{defn} Let $T$ be a rank-one transformation  associated with $(r_n,\omega_n,s_n)_{n=1}^\infty$.
It is called {\it partially bounded} if there is $L>0$ such that $r_n\le L$, $\omega_n(0)=\cdots=\omega_n(r_n-1)$, $\max_{0\le i<j<r_n-1}|s_n(i)-s_n(j)|<L$, $s_n(r_n-1)=0$ and $\min_{0\le i<r_n-1}s_n(i)\ge h_n$ for each $n>0$
 \cite{Gaeb}.
\end{defn}

It was shown in  \cite{Gaeb} that for each partially bounded transformation, the centralizer consists of just the powers. 
Of course, the family of partially bounded transformations does not intersect the set of zero type rank-one maps.

\subsection{Asymmetry and Bergelson's question}
 We say that $T$ is {\it asymmetric} if $T$ is not isomorphic to $T^{-1}$.
 Explicit examples of asymmetric infinite rank-one transformations are constructed in \cite{Ry2} and
 \cite{DanR3} (see there for asymmetric maps which embed into a flow).
It was shown in \cite{Gaeb} that if $T$ is a partially bounded rank-one transformation
then  $T$ is isomorphic to $T^{-1}$ if and only if $s_n(i)=s_n(r_n-2-i)$ for all $i=0,\dots, r_n-2$ eventually in $n$.

Bergelson asked: is there $T$ of infinite ergodic index such that  
 $T\times T^{-1}$ is not ergodic? 
 Of course, such a $T$ is asymmetric. 
 The question is still  open.
 However some partial progress was achieved in \cite{clancy} and \cite{Dan2016}. 
 In \cite{clancy}, an example of a rank-one $T$ was constructed such that $T\times T$ is ergodic but $T\times T^{-1}$ is not. 
 Similar examples appeared  also in \cite{Dan2016}.
 However they do not answer Bergelson's question because $T$ has ergodic index 2 in these examples.
 It was also shown in \cite{Dan2016} that within the class of infinite Markov shifts, the answer on Bergelson's question is negative.
 As for the rank-one transformations, it was shown in \cite{clancy} that $T\times T^{-1}$ is always conservative.

\subsection{Ergodicity of powers}
Let $T$ be a rank-one infinite measure-preserving transformation associated with
$(r_n,\omega_n,s_n)_{n=1}^\infty$.  (Hence $\omega_n(0)=\cdots=\omega_n(r_n-1)$.)
For each $n>0$, we denote by $C_n$ the set of bottom levels of all copies of $(n-1)$-th tower in the $n$-th tower.
Thus, formally, $C_n:=\{0\}\cup \{jh_{n-1}+\sum_{i=0}^{j-1}s_n(i)\mid j=1,\dots, r_n-1\}$.
The following was proved in \cite{DanAnti}.

\begin{theorem}[Ergodicity of powers of rank-one transformations]
\begin{itemize}
\item[(i)]
If $T^d$ is ergodic then for each divisor $p$ of $d$ there are infinitely many $n$ such that some $c\in C_n$ is not divisible by $p$.
\item[(ii)]
If $(r_n)_{n=1}^\infty$ is bounded and for each divisor $p$ of a positive integer $d$,
there are infinitely many
$n$
such that
$p$
does not divide some $c\in C_n$ then $T^d$ is ergodic.
\item[(iii)]
If the sequence $(r_n)_{n=1}^\infty$ is bounded
 then
$T$
is totally ergodic if and only if
for each
$d >
1$,
 there are infinitely many
$n >
0$
such that some element $c\in C_n$
is not divisible by $d$.
\end{itemize}
\end{theorem}

\subsection{Rigidity sequences}
Adams proved in \cite{AdRigid} that if $S$ is an ergodic   probability preserving transformation such that $S^{n_k}\to\text{Id}$ weakly for some sequence $(n_k)_{k=1}^\infty$ of positive integers  then there exists a {\it power rationally weakly 
mixing} infinite measure-preserving  transformation $T$ such that $T^{n_k}\to\text{Id}$ weakly.
The converse assertion is obvious---just pass to the Poisson suspension.
(An infinite measure-preserving $T$ is called power rationally weakly 
mixing if for each finite sequence of non-zero integers $l_1,\dots,l_k$, the product transformation $T^{l_1}\times\cdots\times T^{l_k}$ is rationally weakly mixing.)
Thus the class of rigidity sequences for the ergodic probability preserving transformations equals to the class of rigidity sequences for the ergodic infinite measure-preserving ones.

\subsection{Directional recurrence} Given an ergodic infinite measure-preserving  $\Bbb Z^d$-action
$T=(T_g)_{g\in\Bbb Z^d}$, it seems natural to study the dynamics of individual transformations $T_g$ when $g$ runs $\Bbb Z^d$.
For instance, one of the natural questions is to describe the set 
of those $g\in G$
 such that $T_g$ is recurrent.
Since $T_g$ is recurrent if and only if $T_{ng}$ is recurrent for every $n\in\Bbb Z\setminus\{0\}$, it make sense to talk about {\it recurrent directions} in $\Bbb Z^d$ or {\it rational recurrent directions} in the group $\Bbb R^d$ containg $\Bbb Z^d$  as a standard lattice.
Following Milnor's general idea of directional dynamics \cite{Mil}, Johnson and {\c S}ahin 
introduced in \cite{JS} a concept of directional recurrence of $T$ along an arbitrary (including irrational) direction in $\Bbb R^d$.
They showed that the set $R(T)$ of recurrent directions is a $G_\delta$-subset, produced 
examples of $T$ with  trivial and non-trivial $R(T)$ and asked about description of all possible $R(T)$
when $T$ runs the ergodic $\Bbb Z^d$-actions.
Some partial answers were obtained 
in  \cite{Dan2017}: given a $G_\delta$-subset $\Delta$ of the real projective space $P(\Bbb R^d)$ and a countable subset $D\subset\Delta$, there is a rank-one action $T$ with $D\subset R(T)\subset\Delta$.
However, in general, the problem remains open.

\section{Applications. Connections with other fields}\label{applic}
In this---final---section we shed light on numerous mathematical sources of nonsingular systems. They come from  the theory of stochastic processes, random walks, locally compact Cantor systems, horocycle flows on hyperbolic surfaces, von Neumann algebras, statistical mechanics, representation theory for  groups and anticommutation relations, etc. We also note that such systems sometimes appear in the context of probability preserving dynamics (see also a criterium of distality in terms of the Krengel entropy in  \S\,\ref{Krengel}).

\subsection{Mild mixing}\label{mild mixing}
An ergodic  finite measure-preserving dynamical system $(X,\sS,\mu,T)$ is called {\it mildly mixing} if for each non-trivial  $T$-invariant $\sigma$-algebra $\mathcal A\subset\sS$, the restriction $T\restriction\mathcal A$ is not rigid.
For equivalent definitions and extensions to actions of locally compact groups we refer to \cite{Aa} and \cite{sw}.
 There is an interesting criterium of the mild mixing that involves nonsingular systems: $T$ is mildly mixing if and only if for each ergodic nonsingular transformation $S$, the product $T\times S$ is ergodic \cite{fw}. 
Furthermore,  $T$ is mildly mixing if and only if for each ergodic nonsingular transformation $S$, the product
 $T\times S$ is  orbit equivalent to $S$ \cite{HS97}.
 In particular, the associated flows of $S$ and $T\times S$ are isomorphic.
  Moreover, if  $R$ is a nonsingular transformation such that $R\times S$ is ergodic for any ergodic nonsingular $S$ then $R$ is of type $II_1$ (and mildly mixing) \cite{sw}. 
  In this context we note that 
for every ergodic infinite measure-preserving transformation $T$
there is an ergodic Markov shift $S$ such that $T\times S $ is not conservative, hence not ergodic \cite{ALW}; also 
that $S$ can be chosen to be  rank-one and rigid \cite{EKRSS}.

\subsection{Ergodicity of Gaussian cocycles}
Let $T$ be an ergodic (equivalently, weakly mixing) measure preserving Gaussian transformation on a standard probability space  $(X,\frak B,\mu)$ and let $H$ be the corresponding invariant Gaussian subspace of the real Hilbert space $L^2_0(X,\mu):=L^2(X,\mu)\ominus\Bbb R$.
The following conjecture was stated in  \cite{LeLeSk}:
for each function $f\in H$, either $f$ is a $T$-coboundary (equivalently, a Gaussian coboundary, i.e. the transfer function belongs to $H$) or the skew product transformation $T_f$ acting on the product space $(X\times\Bbb R,\mu\times\text{Leb})$ is ergodic.
We note that  $T_f$ preserves the infinite measure $\mu\times\text{Leb}$.
The affirmative answer was obtained in \cite{DanL22} (and independently in \cite{MaVa}) under an assumption that $T$ is mildly mixing.
Some examples of ergodic $T_f$ with rigid weakly mixing $T$ were constructed in \cite{MaVa}.
However, in the general setting, the problem remains open.

\subsection{Disjointness and Furstenberg's class $\mathcal W^\perp$} Two probability preserving systems $(X,\mu,T)$ and $(Y,\nu,S)$ are called
{\it disjoint} if $\mu\times\nu$ is the only $T\times S$-invariant
probability measure on $X\times Y$ whose coordinate projections
 are $\mu$ and $\nu$ respectively. Furstenberg in \cite{Fur}
initiated studying  the class $\mathcal W^\perp$ of transformations
disjoint from all weakly mixing ones.
Let $\mathcal D$ denote the class of distal transformations
and $\mathcal M(\mathcal W^\perp)$ the class of multipliers
of $\mathcal W^\perp$ (for the definitions see \cite{Glas}).
Then $\mathcal D\subset\mathcal M(\mathcal W^\perp)\subset \mathcal W^\perp$.
 In \cite{LemP} and \cite{DanL} it was shown by constructing
explicit examples that these inclusions are strict. We record this
 fact here because nonsingular ergodic theory
 was the key ingredient of the arguments in the two papers pertaining to
the theory of probability preserving systems.
The examples
are of the form $T_{\phi,S}(x,y)=(Tx,S_{\phi(x)}y)$, where $T$ is an
ergodic rotation on $(X,\mu)$, $(S_g)_{g\in G}$ a mildly mixing action of a
locally compact group $G$ on $Y$ and $\phi:X\to G$ a measurable map.
Let $W_\phi$ denote the Mackey action of  $G$ associated with $\phi$ and let $(Z,\kappa)$ be the space of this action.
The key observation is that there exists an affine isomorphism of the simplex
of $T_{\phi,S}$-invariant  probability measures whose pullback on $X$ is $\mu$ and the simplex of $W_\phi\times S$ quasi-invariant probability measures whose pullback on $Z$ is $\kappa$ and whose Radon-Nikodym cocycle is measurable with respect to $Z$. This is a far reaching generalization of Furstenberg theorem on relative unique ergodicity of ergodic compact group extensions.

\subsection{Symmetric stable and infinitely divisible stationary processes}
Rosinsky in \cite{Ros} established a remarkable connection between structural studies of stationary stochastic processes and ergodic theory of nonsingular transformations (and flows). For simplicity we consider only real processes in discrete time. Let $X=(X_n)_{n\in\mathbb Z}$ be a measurable stationary symmetric $\alpha$-stable (S$\alpha$S) process, $0<\alpha<2$. This means that any linear combination $\sum_{k=1}^n a_kX_{j_k}$, $j_k\in\mathbb Z$, $a_k\in\mathbb R$ has  an S$\alpha$S-distribution. (The case $\alpha=2$ corresponds to Gaussian processes.) Then the process admits a spectral representation
\begin{equation}\label{integral}
X_n=\int_Yf_n(y)\,M(dy), \ n\in\mathbb Z,
\end{equation}
 where $f_n\in L^\alpha(Y,\mu)$ for a standard $\sigma$-finite measure space $(Y,\sS,\mu)$ and $M$ is an independently scattered random measure on $\sS$ such that $E\exp{(iuM(A))}=\exp{(-|u|^\alpha\mu(A))}$ for every $A\in\sS$ of finite measure. By \cite{Ros}, one can choose the kernel $(f_n)_{n\in\mathbb Z}$ in a special way: there are a $\mu$-nonsingular transformation $T$ and  measurable maps $\phi:X\to\{-1,1\}$ and $f\in L^\alpha(Y,\mu)$ such that
$f_n=U^nf$, $n\in\mathbb Z$, where $U$ is the isometry of $L^\alpha(X,\mu)$ given by $Ug=\phi\cdot({d\mu\circ T}/{d\mu})^{1/\alpha}\cdot g\circ T$.
If, in addition, the smallest $T$-invariant $\sigma$-algebra containing $f^{-1}(\sS_{\mathbb R})$ coincides with $\sS$ and Supp$\{f\circ T^n:n\in\mathbb Z\}=Y$ then the pair $(T,\phi)$ is called minimal. It turns out that minimal pairs always exist. Moreover, two minimal pairs $(T,\phi)$ and $(T',\phi')$ representing the same S$\alpha$S process are equivalent in some natural sense \cite{Ros}. Then one can relate ergodic-theoretical properties of $(T,\phi)$ to  probabilistic properties of $(X_n)_{n\in\mathbb Z}$. For instance, let $Y=C\sqcup D$ be the Hopf decomposition of $Y$ (see Theorem~\ref{T:hopf}). We let $X_n^D:=\int_Df_n(y)\,M(dy)$ and $X_n^C:=\int_Cf_n(y)\,M(dy)$.
Then we obtain a unique (in distribution) decomposition of $X$ into the sum
$X^D+X^C$ of two independent stationary S$\alpha$S-processes.

Another kind of decomposition was considered in \cite{Sam}. Let $P$ be the largest invariant subset of $Y$ such that $T\restriction P$ has a finite invariant measure. Partitioning $Y$ into $P$ and $N:=Y\setminus N$ and restricting the integration in (\ref{integral}) to $P$ and $N$ we obtain a unique (in distribution) decomposition of $X$ into the sum $X^P+X^N$ of two independent stationary S$\alpha$S-processes. Notice that the process $X$ is ergodic if and only if $\mu(P)=0$.

 Roy considered a more general class of {\it infinitely divisible (ID)} stationary processes \cite{Roy1}. Using Maruyama's representation of the characteristic function of an ID process $X$ without Gaussian part he  singled out the L\'evy measure $Q$ of $X$. Then $Q$ is a shift invariant $\sigma$-finite measure on $\mathbb R^\mathbb Z$. Decomposing the dynamical system $(\mathbb R^\mathbb Z,\tau,Q)$ in various natural ways (Hopf decomposition, 0-type and positive type, so-called `rigidity free' part and its complement) he obtains corresponding decompositions for the process $X$. Here $\tau$ stands for the shift on $\mathbb R^\mathbb Z$.

\subsection{Poisson suspensions of infinite measure preserving transformations}
Poisson suspensions over infinite measure-preserving transformations (see \S\,\ref{NPT} or  \cite{Roy}, \cite{Roy2})
 are widely used in statistical mechanics to model ideal gas, Lorentz gas, etc (see \cite{CFS}).
Together with the Gaussian dynamical systems
they are also an important source of examples and counterexamples in ergodic theory.
Due to a close
 similarity with the well studied  Gaussian systems, a natural question arises: are there ergodic Poisson suspensions whose ergodic self-joinings are Poissonian?
 Such suspensions are called PAP. 
 They are analogue of GAG in the theory of Gaussian systems \cite{LPT}.
Janvresse, de la Rue and Roy  constructed  PAP suspension in \cite{Jan2017} (see also \cite{ParR}).
The example of an infinite measure-preserving $T$ with ``minimal self-joinnings''  from \cite{JanRR2} plays a crucial role in their construction.
We also mention a result of Meyerovitch \cite{Mey1} related to weak mixing of infinite measure-preserving systems and Poisson suspensions: if $T$ is a conservative ergodic infinite measure-preserving transformation then the direct product of $T$ with the Poisson suspension
$T_*$ of $T$ is ergodic.

\subsection{Recurrence of random walks with non-stationary increments}
Using nonsingular ergodic theory one can introduce the notion of recurrence for random walks obtained from certain non-stationary processes.  Let $T$ be an ergodic  nonsingular transformation of a standard probability space $(X,\sS,\mu)$ and let $f:X\to\mathbb R^n$ a measurable function. Define for $m\ge 1$, $Y_m:X\to\mathbb R^n$ by $Y_m:=\sum_{n=0}^{m-1} f\circ T^n$. In other words, $(Y_m)_{m\ge 1}$ is the random walk associated with the (non-stationary) process $(f\circ T^n)_{n\ge 0}$. Let us call this random walk {\it recurrent} if the cocycle $f$ of $T$ is recurrent (see \S\,\ref{cocycles}). It was shown in \cite{Sc5} that in the case $\mu\circ T=\mu$, i.e., the process is stationary, this definition is equivalent to the standard one.

\subsection{Boundaries of random walks}\label{brw}
Boundaries of random walks on groups retain valuable information on the underlying groups (amenability, entropy, etc.) and enable one to obtain integral representation for harmonic functions of the random walk \cite{Zim}, \cite{Zim2}, \cite{KV}.
 Let $G$ be a locally compact group and $\nu$ a probability measure on $G$. Let $T$ denote the (one-sided) shift on the probability space $(X,\sS_X,\mu):=(G,\sS_G,\nu)^{\mathbb Z_+}$ and $\phi:X\to G$ a measurable map defined by $(y_0,y_1,\dots)\mapsto y_0$. Let $T_\phi$ be the $\phi$-skew product extension of $T$ acting on the space $(X\times G,\mu\times\lambda_G)$ (for non-invertible transformations the skew product extension is defined in the very same way as for invertible ones, see \S\,\ref{cocycles}). Then $T_\phi$ is isomorphic to the {\it homogeneous random walk} on $G$ with jump probability $\nu$. Let $\mathcal I(T_\phi)$ denote the sub-$\sigma$-algebra of $T_\phi$-invariant sets and let $\mathcal F(T_\phi):=\bigcap_{n>0}T_\phi^{-n}(\sS_X\otimes\sS_G)$. The former is called the {\it Poisson boundary} of $T_\phi$ and the latter one is called the {\it tail boundary} of $T_\phi$.  Notice that a nonsingular action of $G$ by inverted right translations along the second coordinate is well defined on each of the two boundaries.  The two boundaries (or, more precisely, the $G$-actions on them) are ergodic. The Poisson boundary is the Mackey range of $\phi$ (as a cocycle of $T$). Hence the Poisson boundary is amenable \cite{Zim}. If the support of $\nu$ generates a dense subgroup of $G$ then the corresponding Poisson boundary is weakly mixing \cite{AaL}. As for the tail boundary, we first note that it can be defined for a wider family of {\it non-homogeneous}  random walks. This means that the jump probability $\nu$ is no longer fixed and a sequence $(\nu_n)_{n>0}$ of probability measures on $G$ is considered instead. Now let $(X,\sS_X,\mu):=\prod_{n>0}(G,\sS_G,\nu)$. The one-sided shift on $X$ may not be nonsingular now. Instead of it, we consider the tail equivalence relation $\mathcal R$ on $X$ and a cocycle $\alpha:\mathcal R\to G$ given by $\alpha(x,y)=x_1\cdots x_ny_n^{-1}\cdots y_1$, where $x=(x_i)_{i>0}$ and $y=(y_i)_{i>0}$ are $\mathcal R$-equivalent and $n$ in the smallest integer such that $x_i=y_i$ for all $i>n$. The tail boundary of the random walk on $G$ with time dependent jump probabilities $(\nu_n)_{n>0}$ is  the Mackey $G$-action associated with $\alpha$. In the case of homogeneous random walks this definition is equivalent to the initial one. Connes and Woods showed \cite{CW} that the tail boundary is always amenable and AT. It is unknown whether the converse holds for general $G$. However it is true for $G=\mathbb R$ and $G=\mathbb Z$: the class of  AT-flows coincides with the class of tail boundaries of the random walks on $\mathbb R$ and a similar statement holds for $\mathbb Z$ \cite{CW}.
Jaworski showed \cite{Ja} that if $G$ is countable and a random walk is homogeneous then the tail boundary of the random walk possesses a so-called SAT-property (which is stronger than AT).

\subsection{Stationary actions} Let $T=(T_g)_{g\in G}$ be a continuous action of a countable group $G$ on a compact metrizable space $X$.
By Markov-Kakutani theorem, if $G$ is amenable then there is an invariant probability Borel measure on $K$.
If $G$ is non-amenable such a measure does not necessarily exist.
However, if $\kappa$ is a probability measure on $G$ whose support generates $G$ as a semigroup then there is always a $T$-quasiinvariant probability measure $\mu$ on $X$ such that
$$
\sum_{g\in G}\kappa(g)\frac{d\mu\circ T_g^{-1}}{d\mu}(x)=1\  \ \text{for a.e. }x\in X.
$$
$\mu$ is called  a {\it $\kappa$-stationary} measure.
For a deep theory of stationary actions and its applications we refer to
 \cite{FurG}, \cite{Furman} and references therein.
 However, if
 $G$ is Abelian or, more generally, nilpotent then each $\kappa$-stationary action is invariant under $T$.
 Thus, there are no stationary $\Bbb Z$-actions except for the probability preserving ones.
 That is why we do not discuss them in this survey.

\subsection{Classifying $\sigma$-finite ergodic invariant measures}
The description of ergodic finite invariant measures for topological (or, more generally, standard Borel) systems is a well established  problem in the classical ergodic theory \cite{CFS}. On the other hand, it seems impossible to obtain any useful information about the system by analyzing the set of all ergodic quasi-invariant (or just  $\sigma$-finite invariant) measures because this set is wildly huge (see \S\,\ref{S:Glimm}). The situation changes if we impose some restrictions on the measures.
For instance, if the system under question is a homeomorphism (or a topological flow) defined on a locally compact Polish space then it is natural to consider the class of ($\sigma$-finite) invariant Radon measures, i.e., measures taking finite values on the compact subsets. We give two examples.

First, the seminal results of Giordano, Putnam and Skau on the topological orbit equivalence of compact Cantor minimal systems were extended to locally compact Cantor minimal (l.c.c.m.) systems in \cite{Dan27} and \cite{Mat}. Given a l.c.c.m. system $X$, we denote by $\mathcal M(X)$ and $\mathcal M_1(X)$ the set of invariant Radon measures and the set of invariant probability measures on $X$. Notice that $\mathcal M_1(X)$ may be empty \cite{Dan27}. It was shown in \cite{Mat} that two systems $X$ and $X'$ are topologically orbit equivalent if and only if there is a homeomorphism of $X$ onto $X'$ which maps bijectively $\mathcal M(X)$ onto $\mathcal M(X')$ and $\mathcal M_1(X)$ onto $\mathcal M_1(X')$.
Thus $\mathcal M(X)$ retains an important information on the system---it is `responsible' for the topological orbit equivalence of the underlying systems. Uniquely ergodic l.c.c.m. systems (with unique up to scaling infinite invariant Radon measure)  were constructed in \cite{Dan27}.

The second example is related to study of the smooth horocycle  flows on tangent bundles of hyperbolic surfaces.
Let $\mathbb D$ be the open disk equipped with  the hyperbolic metric
$|dz|/(1-|z|^2)$ and let M\"{o}b$(\mathbb D)$ denote the group of M\"{o}bius
transformations of $\mathbb D$.
A hyperbolic surface can be written in the form
$M:=\Gamma\backslash\text{M\"{o}b}(\mathbb D)$,
where $\Gamma$ is a torsion free discrete subgroup of M\"{o}b$(\mathbb D)$.
Suppose that $\Gamma$ is a nontrivial normal subgroup of a lattice
$\Gamma_0$ in M\"{o}b$(\mathbb D)$.
Then $M$ is a regular cover of the finite volume surface
$M_0:=\Gamma_0\backslash\text{M\"{o}b}(\mathbb D)$.
 The group of deck transformations $G=\Gamma_0/\Gamma$ is finitely generated.
The horocycle flow $(h_t)_{t\in\mathbb R}$ and the geodesic flow
$(g_t)_{t\in\mathbb R}$ defined on the unit tangent bundle $T^1(\mathbb D)$
descend naturally to flows, say $h$ and $g$, on $T^1(M)$. We consider the
problem of classification of the $h$-invariant Radon measures on $M$.
According to Ratner, $h$ has no finite invariant measures on $M$ if $G$ is
infinite (except for measures supported on closed orbits). However there
are infinite invariant Radon measures, for instance the volume measure. In
the case when $G$ is free Abelian and $\Gamma_0$ is co-compact, every
homomorphism $\phi:G\to\mathbb R$ determines a unique up to scaling ergodic
invariant Radon measure (e.i.r.m.) $m$ on $T^1(M)$ such that $m\circ
dD=\exp(\phi(D))m$ for all $D\in G$ \cite{BabL} and every e.i.r.m. arises
this way \cite{Sar}. Moreover all these measures are quasi-invariant under
$g$. In the general case, an interesting bijection is established in
\cite{LedS} between the e.i.r.m. which are quasi-invariant under $g$ and
the `non-trivial minimal' positive eigenfunctions of the hyperbolic
Laplacian on $M$.

\subsection{Von Neumann algebras} There is a fascinating and productive interplay between nonsingular ergodic theory and von Neumann algebras. The two theories alternately influenced development of each other. Let $(X,\sS,\mu,T)$ be a nonsingular dynamical system. Given $\phi\in L^\infty(X,\mu)$ and $j\in \mathbb Z$, we define operators $A_\phi$ and $U_j$ on the Hilbert space $L^2(Z\times\mathbb Z,\mu\times\nu)$ by setting
$$
(A_\phi f)(x,i) := \phi(T^ix)f(x,i),\ (U_j f)(x,i) := f(x,i-j)
$$
Then $U_jA_\phi U^*_j=A_{\phi\circ T^j}$. Denote by $\mathcal M$
the von Neumann algebra (i.e., the weak closure of the $*$-algebra) generated by $A_\phi$, $\phi\in L^\infty(X,\mu)$ and $U_j$, $j\in\mathbb Z$. If $T$ is ergodic and aperiodic then $\mathcal M$ is a factor, i.e., $\mathcal M\cap\mathcal M'=\mathbb C 1$, where $\mathcal M'$ denotes the algebra of bounded operators commuting with $\mathcal M$. It is called a {\it Krieger's} factor.  Murray-von Neumann-Connes' type of $\mathcal M$ is exactly the Krieger's type of $T$. The flow of weights of $\mathcal M$ is isomorphic to the associated flow of $T$.
Two  Krieger's factors are isomorphic if and only if the underlying dynamical systems are orbit equivalent \cite{Kr76}. Moreover, a number of important problems in the theory of von Newmann algebras  such as classification of subfactors, computation of the flow of weights and Connes' invariants, outer conjugacy  for automorphisms, etc. are intimately related to the corresponding problems in  nonsingular orbit theory. We refer to \cite{Moo1}, \cite{FM}, \cite{GiS}, \cite{GiS2}, \cite{HamKos}, \cite{DanH} for details.

\subsection{Representations of CAR} Representations of canonical anticommutation relations (CAR) is one of  the most elegant and useful chapters of mathematical physics, providing a natural language for many body quantum physics and quantum field theory. By a representation of CAR we mean a sequence of bounded linear operators $a_1,a_2,\dots$ in a separable Hilbert space $\mathcal K$ such that
$
a_ja_k+a_ka_j=0 \ \text{and }a_ja_k^*+a_k^*a_j=\delta_{j,k}.
$

Consider $\{0,1\}$ as a group with addition mod\,2. Then $X=\{0,1\}^\mathbb N$ is a compact Abelian group. Let $\Gamma:=\{x=(x_1, x_2,\dots):\lim_{n\to\infty}x_n=0\}$. Then $\Gamma$ is a dense countable subgroup of $X$. It is generated by elements $\gamma_k$ whose $k$-coordinate is 1 and the other ones are 0. $\Gamma$ acts on $X$ by translations.  Let $\mu$ be an ergodic $\Gamma$-quasi-invariant measure on $X$. Let $(C_k)_{k\ge 1}$ be  Borel maps from $X$ to the group of unitary operators in a Hilbert space $\mathcal H$ satisfying $C_k^*(x)=C_k(x+\delta_k)$, $C_k(x)C_l(x+\delta_l)=C_l(x)C_k(x+\delta_k)$, $k\ne l$ for a.a. $x$.
In other words, $(C_k)_{k\ge 1}$ defines a cocycle of the $\Gamma$-action.
 We now put $\widetilde{\mathcal H}:=L^2(X,\mu)\otimes \mathcal H$ and
define operators $a_k$ in $\widetilde{\mathcal H}$ by setting
$$
(a_kf)(x)=(-1)^{x_1+\cdots+x_{k-1}}(1-x_k)C_k(x)\sqrt{\frac
{d\mu\circ\delta_k}{d\mu}(x)}f(x+\delta_k),
$$
where $f:X\to\mathcal H$ is an element of $\widetilde{\mathcal H}$ and $x=(x_1,x_2,\dots)\in X$. It is easy to verify that $a_1,a_2,\dots$ is a representation of CAR. The converse was established in \cite{GorW} and \cite{Golod}: every factor-representation (this means that the von Neumann algebra generated by all $a_k$ is a factor) of CAR can be represented as above for some ergodic measure $\mu$, Hilbert space $\mathcal H$ and a $\Gamma$-cocycle $(C_k)_{k\ge 1}$. Moreover, using nonsingular ergodic theory Golodets \cite{Golod} constructed for each $k=2,3,\dots,\infty$, an irreducible representation of CAR  such that dim\,$\mathcal H=k$. This answered a question of G{\aa}rding and Wightman \cite{GorW} who considered only the case $k=1$.

\subsection{Unitary representations of locally compact groups}
Nonsingular actions appear in a systematic way  in the theory of unitary representations of groups. Let $G$ be a locally compact second countable group and $H$ a closed normal subgroup of $G$. Suppose that $H$ is commutative (or, more generally, of type I, see \cite{Dix}).
 Then the natural action of $G$ by conjugation on $H$ induces a Borel $G$-action, say $\alpha$, on  the dual space $\widehat H$---the set of unitarily equivalent classes of irreducible unitary representations of $H$. If now $U=(U_g)_{g\in G}$ is a unitary representation of $G$ in a separable Hilbert space then by applying Stone decomposition theorem to $U\restriction H$ one can deduce that $\alpha$ is nonsingular with respect to a measure $\mu$ of the `maximal spectral type' for $U\restriction H$ on $\widehat H$. Moreover, if $U$ is irreducible then $\alpha$ is ergodic. Whenever $\mu$ is fixed, we obtain a one-to-one correspondence between the set of cohomology classes of irreducible cocycles for $\alpha$ with values in the unitary group on a Hilbert space
$\mathcal H$ and the subset of  $\widehat G$ consisting of classes of those unitary representations $V$ for which the measure associated to $V\restriction H$ is equivalent to $\mu$.  This correspondence is used in both directions. From information about cocycles we can deduce facts about representations and vise versa \cite{Kir}, \cite{Dix}.

\end{document}